\documentclass[12pt]{article}

\usepackage[letterpaper,margin =1in]{geometry}
\usepackage{amssymb}
\usepackage{amsmath}
\usepackage{amsthm}
\usepackage{mathtools}
\usepackage{enumerate}
\usepackage{fancyhdr}
\usepackage{braket}
\usepackage{titlesec}
\usepackage{xcolor}
\usepackage{hyperref}

\usepackage{tikz-cd} 
\usetikzlibrary{babel}

\usepackage{tikz}
\usetikzlibrary{shadows, arrows}

\usepackage{scalerel,stackengine}
\stackMath
\newcommand\reallywidehat[1]{%
\savestack{\tmpbox}{\stretchto{%
  \scaleto{%
    \scalerel*[\widthof{\ensuremath{#1}}]{\kern-.6pt\bigwedge\kern-.6pt}%
    {\rule[-\textheight/2]{1ex}{\textheight}}%WIDTH-LIMITED BIG WEDGE
  }{\textheight}% 
}{0.5ex}}%
\stackon[1pt]{#1}{\tmpbox}%
}

\usepackage{todonotes}

\definecolor{skyblue}{rgb}{0.85,0.85,1}

\titleformat{\section}
  {\normalfont\Large\bfseries}{\thesection}{1em}{}

\titleformat{\subsection}
  {\normalfont\large\bfseries}{\thesubsection}{1em}{}

\DeclareFontFamily{OT1}{pzc}{}
\DeclareFontShape{OT1}{pzc}{m}{it}%
             {<-> s * [1.200] pzcmi7t}{}
      
\newlength\tindent
     \newcommand{\bigO}{\mathcal{O}}

\DeclareMathAlphabet{\mathscr}{OT1}{pzc}{m}{it}

\providecommand{\norm}[1]{\lVert#1\rVert}

\providecommand{\del}{\partial}
\providecommand{\eps}{\varepsilon}

\newenvironment{frcseries}{\fontfamily{frc}\selectfont}{}

\newcommand{\ve}{\varepsilon}
\newcommand{\bupm}{\bar{U}^{\pm}}
\newcommand{\sgn}{\textrm{sgn}}

\newcommand{\cA}{\mathcal{A}}
\newcommand{\cL}{\mathcal{L}}
\newcommand{\cT}{\mathcal{T}}

\newcommand{\C}{\mathbb{C}}
\newcommand{\R}{\mathbb{R}}

\newtheorem{thr}{Theorem}[section]
\newtheorem{lem}[thr]{Lemma}
\newtheorem{coro}[thr]{Corollary}
\newtheorem{hyp}[thr]{Hypothesis}
\theoremstyle{definition}
\newtheorem{defi1}[thr]{Definition}
\newenvironment{defi}{\begin{defi1}\rm}{\hfill $\triangle$ \end{defi1}}

\theoremstyle{remk}
\newtheorem{remark1}[thr]{Remark}
\newenvironment{remk}{\begin{remark1}\rm}{\hfill $\triangle$ \end{remark1}}

\newcommand{\transv}{\mathrel{\text{\tpitchfork}}}
\makeatletter
\newcommand{\tpitchfork}{%
  \vbox{
    \baselineskip\z@skip
    \lineskip-.52ex
    \lineskiplimit\maxdimen
    \m@th
    \ialign{##\crcr\hidewidth\smash{$-$}\hidewidth\crcr$\pitchfork$\crcr}
  }%
}
\makeatother

\begin{document}

\title{Spectral stability of shock-fronted travelling waves under viscous relaxation}
\author{Ian Lizarraga\thanks{School of Mathematics and Statistics, The University of Sydney, {\tt ian.lizarraga@sydney.edu.au} (corresponding author)}  \: and Robert Marangell\thanks{School of Mathematics and Statistics, The University of Sydney, {\tt robert.marangell@sydney.edu.au}.}}
\maketitle

\begin{abstract}
Reaction-nonlinear diffusion partial differential equations can exhibit {\it shock-fronted} travelling wave solutions. Prior work by Li et al. (2021) has demonstrated the existence of such waves for two classes of regularisations, including viscous relaxation  (see \cite{li}). Their analysis uses geometric singular perturbation theory: for sufficiently small values of a parameter $\eps > 0$ characterising the `strength' of the regularisation, the waves are constructed as perturbations of a singular heteroclinic orbit. Here we show rigorously that these waves are spectrally stable for the case of viscous relaxation.  Our approach is to show that for sufficiently small $\eps>0$, the `full' eigenvalue problem of the regularised system is controlled by a reduced {\it slow eigenvalue problem} defined for $\eps = 0$. In the course of our proof, we examine the ways in which this geometric construction complements and differs from constructions of other reduced eigenvalue problems that are known in the wave stability literature.
\end{abstract}

\tableofcontents 

\section{Introduction}

Regularised PDEs exhibiting travelling wave solutions can frequently be regarded as {\it singularly perturbed} dynamical systems, since the regularisation is modelled by the inclusion of strictly higher-order partial derivatives multiplied by a small constant. The singular perturbation in turn has a physical/dynamical interpretation reflecting an inherent hierarchy of spatiotemporal scales, which we refer to more succinctly as a {\it fast-slow} structure. We consider the following reaction-nonlinear diffusion PDE with a viscous relaxation term, exhibiting (regularised) {\it shock-fronted} travelling waves:\footnote{Throughout the paper we use bar notation to refer to phase space variables. Later on we will use unbarred variables to denote variables defined on the linearized subspaces along the curve, and carets to denote projectivisations of these variables.}
\begin{align} \label{eq:master}
\frac{\del \bar{U}}{\del t} &= \frac{\del}{\del x} \left( D(\bar{U}) \frac{\del \bar{U}}{\del x} \right) + R(\bar{U}) + \eps \frac{\del^3 \bar{U}}{\del x^2 \del t}\end{align}

for $(x,t) \in \mathbb{R}\times \mathbb{R}$ and $\eps \geq 0$ the singular perturbation parameter characterising the strength of the the regularisation. The nonlinear diffusion and  reaction terms given by the quadratic resp. cubic polynomials $D(\bar{U})$ and $R(\bar{U})$. We suppose that $D(\bar{U}) < 0$ within an interval $(a,b) \subset (0,1)$; the {\it potential function} $F(\bar{U}) := \int D(\bar{U}) d\bar{U}$ is hence nonmonotone. Nonlinear diffusion processes with nonmonotone potential functions serve as prototypical models for the formation of coherent solutions with sharp fronts; see e.g. \cite{witelski1995,witelski1996}. Arguably the canonical system from a physical context is the Cahn-Hilliard model \cite{pego1989}. The system \eqref{eq:master} for $\eps = 0$ is derived as a particular continuum limit of stochastic agent-based models of invasion processes \cite{simpson1,simpson2}. The choice of viscous regularisation in \eqref{eq:master} appears in the shock regularisation literature \cite{novick,padron,witelski1996}, but there are other physically-motivated high-order regularising terms that give rise to (smoothed families of) shock-fronted travelling waves. We refer the reader to \cite{li} for a comparison of the existence problems under viscous resp. nonlocal regularising terms relative to system \eqref{eq:master}. \\

In the case of existence, the effectiveness of  {\it geometric singular perturbation theory} (GSPT) is now well established. When \eqref{eq:master} is written in a coordinate frame that follows the travelling wave, it takes the structure of a closed fast-slow system of ODEs. A one-parameter family of travelling waves for $0 < \eps \ll 1$ is then constructed rigorously as a perturbation of a singular heteroclinic orbit, composed of segments defined by the slow flow along a so-called critical manifold, concatenated with fast jumps along the fast fibres, according to the dynamics of the layer problem. Li et. al. follow this approach in \cite{li}, numerically demonstrating the existence of smooth travelling waves of \eqref{eq:master} for small values of $\eps>0$. These one-parameter families of waves limit to `genuine' (piecewise continuous) shock-fronted waves as $\eps \to 0$; furthermore, such singular limits are nonunique, strongly depending on the regularisation chosen. Li and coauthors have also studied smooth travelling-wave solutions of the unperturbed problem arising as so-called `hole in the wall' solutions; see \cite{li2} and the footnote on page 9. Generally speaking, one can attempt to construct shock-fronted limiting solutions directly as weak solutions, but there may be infinitely many; see eg. \cite{hollig}.\\

It is worth re-emphasizing the effect of the singular perturbation in \eqref{eq:master}, in terms of coordinate representations of the travelling wave. With respect to a natural Li\'enard-type representation, the unperturbed problem can be thought of as modelling the slow dynamics on the critical manifold itself. The introduction of a regularisation term is then tantamount to embedding this slow flow within a higher-dimensional space and introducing a rule (i.e. the fast layer flow) for connecting the wavefront smoothly. In \cite{li}, it is explicitly shown that viscous relaxation adds a one-dimensional layer flow, while a fourth-order nonlocal regularization term adds a two-dimensional layer flow.\\

The objective of this paper is to demonstrate the {\it spectral stability} of travelling waves arising in \eqref{eq:master}. To determine the spectral stability along a wave, we must find the spectrum $\sigma(L)$ of a corresponding linearized operator $L$. The total spectrum is decomposed into its continuous and point components as $\sigma(L) = \sigma_e(L) \cup \sigma_p(L)$. Typically, most of the work involves deducing the point spectrum $\sigma_p(L)$, by recasting the eigenvalue problem $(L-\lambda I)v = 0$  as a bifurcation problem posed on the underlying linearized subspaces along the wave.  The major point that we want to highlight is that although the existence problem remains amenable to the usual GSPT techniques, the analysis of the stability problem (in particular, the computation of the point spectrum) is quite distinct. \\

The literature on geometric and analytic techniques for linearized stability problems is substantial (see e.g. \cite{kapprom} and \cite{bjornstability}  for comprehensive surveys). In this paper we follow the geometric framework developed in a series of seminal papers by Alexander, Gardner, and Jones \cite{AGJ,GJ,jones}. The key point is that $L$ also inherits a slow-fast structure, enabling a reduction via linearized subsystems defined along either the slow manifold or along the fast fibres. The typical approach is then to hope that the spectrum of the full system, for $\eps >0$ sufficiently small, is close to the `fast' and/or `slow' spectra of these linearized subsystems, which are defined for $\eps = 0$.\\

Our principal result is Theorem \ref{thm:main}. As a consequence of the abovementioned high order of the singular perturbation term representing the regularisation, the stability problem for small $\eps > 0$ is now close to a {\it slow eigenvalue problem} defined on the critical manifold. We give an outline of the proof in Section \ref{sec:outline}. To summarize, we show in Sections \ref{sec:bundles} and \ref{sec:convergence} that the eigenvalues for the perturbed problem are the same as the eigenvalues for the reduced problem, which are calculated in Section \ref{sec:slowevans}. Together with the fact that the essential spectrum lies in the left half complex plane as shown in Section \ref{sec:spectrum},  we can conclude that there is no spectrum of solutions to \eqref{eq:master} in any fixed contour $K$ in the right half complex plane. \\

 The principal estimate for constructing these reduced eigenvalue problems is the {\it elephant trunk lemma},  which characterises an attracting set over a fast unstable subbundle in the linearized space when $\lambda$ is {\it not} a fast eigenvalue; see the construction in the paper of Gardner \& Jones \cite{GJ}. We highlight that a nondegenerate fast eigenvalue problem is a requirement to construct linked elephant trunks over the entire wave. Furthermore, the construction of a slow  eigenvalue problem, corresponding to the singular limit of some slow subbundle, is essentially auxiliary to that of the fast problem, and many of the key properties of the slow subbundle are indirectly {\it enforced} by the elephant trunk over the fast subbundle. \\

In order for us to define the promised slow eigenvalue problem, it is essential that we determine how exactly slow linear data on one branch of the slow manifold is transported across fast fibres to another branch of the slow manifold. When an elephant trunk lemma is available, it enforces certain nice properties, such as continuity of the slow subbundle across the fast layer in the singular limit. In other words, slow linear data from one branch of the slow manifolds is transported across fast fibres to the other branch {\it identically}, allowing a straightforward definition of slow eigenvalue problems across disjoint subsets of the slow manifolds.\\

We will show that our problem does not allow for elephant trunk-type estimates over the fast layer. The main obstruction is that, relative to a natural set of coordinates for the linearized system so that setting $\lambda = 0$ recovers the usual variational equations along the flow, the parameter $\lambda$ enters the equations `weakly,' i.e. only through $\bigO (\eps \lambda)$ terms. The singular limit of the fast linearized system therefore degenerates such that there is `always a fast eigenvalue,' i.e. there is a connection made between unstable and stable directions across the wavefront for all $\lambda \in \mathbb{C}$ as $\eps \to 0$.  In the present case the eigenvalue problem is low-dimensional enough that we can partly avoid this problem by concentrating on the slow problem only, but we return to this issue in the conclusion when we consider fourth-order regularizations. \\

Consequently, it must be the case that the slow linearized dynamics is responsible for the generation of eigenvalues near the tails of the wave. In this context, we think of our work as complementing Jones' early result on the stability of a travelling pulse in the FitzHugh-Nagumo system, in which the eigenvalue problem of the full system is strongly controlled by the fast eigenvalue problems along the  wave front and back \cite{jones}; as well as the result of Gardner and Jones on the stability of travelling waves in predator-prey systems, where both fast and slow eigenvalue problems are constructed \cite{GJ}. \\

While we do not have access to a fast elephant trunk over the entire wave, we nevertheless retain control of the linearized dynamics for small values of $\eps$ since $\lambda$ enters the equations only weakly. The key observation here is that the eigenvalue problem can be  thought of as a `weak' ($\mathcal{O}(\eps)$)-perturbation of the standard variational problem. 
We will show that the corresponding transport of slow data is not the identity map, but instead is replaced by a nontrivial {\it jump map} that we can write down explicitly in the limit as $\eps \to 0$. The key technical tools that we use are desingularized linearized slow flows;  $\eps$-dependent rescalings along the fast fibres; and a continuous differentiability criterion that holds across the fold. We remark that this analysis holds also for the $\lambda = 0$ case (corresponding to the regular variational dynamics carried by the wave), but to our knowledge such a statement about jump compatibility is new in the GSPT literature.\\

 The {\it exchange lemma}, another fundamental GSPT technique (see e.g. \cite{joneskopell,jonestin}),\footnote{Strictly speaking, in this paper we consider the case of an {\it inclination lemma} (see e.g. \cite{brunovsky,schecter}). We also point out that the adapation of the standard $(k+\sigma)$-exchange lemma to the case of zero unstable fast directions is straightforward.} can also be adapted to our eigenvalue problem, since it is $\mathcal{O}(\eps)$-close to the standard variational problem. Exchange lemma-type estimates allow us to track solutions of the eigenvalue problem as they leave the fast layer and enter small neighborhoods of the slow manifolds, where they enter partial elephant trunks and remain well-controlled along the tails of the wave. The key result here is that solutions entering neighborhoods of the slow manifolds in a generic way will be aligned closely to a {\it slow subbundle} defined on the slow manifold after sufficiently long times. This subbundle is defined from the eigenvalue problem, and is therefore not generally equivalent to the tangent bundle of the slow manifold (but it is nearby). \\

These new techniques cause no extra trouble for the topological arguments that allow us to calculate the point spectrum of the perturbed problem: the corresponding evaluation of first Chern numbers of certain {\it augmented unstable bundles} (defined in \cite{AGJ}) follows, as usual. Following the approach of Gardner \& Jones in \cite{GJ}, we construct a homotopy between the augmented unstable bundle $\mathcal{E}_{\eps}(K)$ of the `full' problem with $\eps > 0$ and that of a reduced problem $\mathcal{E}_0(K)$, when $\eps$ is sufficiently small. The reduced vector bundle is defined over separate hemispheres using the slow eigenvalue problem. The jump map now plays an essential role in defining the appropriate clutching function, which glues these vector bundles together along the edges of the hemispheres. We also highlight that the translational eigenvalue at $\lambda =0$ is now counted by the {\em slow} problem. \\

In the course of proving stability in our problem, we also demonstrate explicitly the geometric consequence of eigenvalue crossings in the full system for $0 < \eps \ll 1$: the projectivized solutions wind around according to the number of eigenvalues crossed, and these winds occur entirely within the slow dynamics of the critical manifold. The propagation of winding in the reduced problem, as well as the generation of new winds (eigenvalues) can be identified by locating the poles and zeroes, respectively, of a meromorphic {\it Riccati-Evans function} \cite{harley}, defined with respect to a judiciously chosen section of the projective dynamics. This also gives us another opportunity to demonstrate the computational utility of the Riccati-Evans function to count eigenvalues via the argument principle.\\

We wish to highlight some results in the literature in the context of our present problem. Slow eigenvalues appear in a doubly-diffusive Fitzhugh-Nagumo system which is considered as an early application of the augmented unstable bundle theory in  \cite{AGJ}, but here the extra fast dynamics is relatively trivial; furthermore, the `slow' problem is in fact fast-slow relative to another small parameter, with its eigenvalues arising from the corresponding fast front and back.  The occurrence of such a nontrivial slow eigenvalue problem as described in our current analysis appears to be quite atypical; we highlight the work of Bose \cite{bose} in this vein as one of the only comparable examples we have found in the literature. Bose provides similar geometrical results for pulse (meta)stability in singularly-perturbed nonlocal RDEs.  The singular perturbation occurs at the second order, so the resulting nonlocal eigenvalue problem can then be studied with traditional Sturm-Liouville techniques, and an oscillation theorem is proven. Relative to Bose's result, we have two additional technical issues to settle: (i) we must first `reduce' our singularly-perturbed third-order eigenvalue problem (see \eqref{eq:fulleigprob}) by showing closeness to a second-order eigenvalue problem defined on the critical manifold, and (ii) we must also find the corresponding compatibility condition to connect the eigenvalue problem across disjoint branches of the critical manifold.\\

 We also wish to highlight de Rijk et. al.'s recent work on analytical slow-fast factorizations of Evans functions \cite{derijk}; the relationship between the geometric obstructions discussed above and analytic conditions to produce reduced Evans functions in de Rijk et. al.'s paper bears further exploration, but we do not attempt this in our paper.\\

We find it sensible to closely follow the general structure in \cite{GJ} in organizing our paper.  In Sec. \ref{sec:waves}, we describe the construction of a one-parameter family of travelling waves limiting onto a singular shock-fronted travelling wave using geometric singular perturbation theory. In Sec. \ref{sec:geometry}, we define the relevant geometric spaces in which our objects of interest lie and write down some key facts about them. In Sec. \ref{sec:spatialeigs} we write down the spatial eigenvalue problem for the associated linearized operator of \eqref{eq:master}, and we write down some leading-order estimates for the eigenvalues and eigenvectors of the associated asymptotic systems. We also define the linearized slow and fast subsystems that we will use often to estimate the linear dynamics of the `full' system. In Sec. \ref{sec:spectrum}, we describe some general stability theory for travelling waves and show that the essential spectrum is bounded away from the imaginary axis in the left-half complex plane. In Sec. \ref{sec:bundles}, we construct the augmented unstable vector bundles over a complex contour enclosing the remaining point spectrum, and we also take the opportunity to outline the proof of the main stability theorem in the paper. In Sec. \ref{sec:convergence}, we prove the main uniform estimates for the linearized dynamics in terms of the singular subsystems. We use these estimates to construct a homotopy between some augmented unstable bundles. In Sec. \ref{sec:slowevans}, we characterize the point spectrum of the reduced subsystem, finally allowing us to deduce the stability of the wave of the full system for $\eps > 0$ sufficiently small. We conclude in Sec. \ref{sec:conclusion}, and the most technical lemmas are relegated to the Appendices.

\section{Existence of travelling waves} \label{sec:waves}

Our first task is to construct constant-speed travelling waves for the PDE \eqref{eq:master}. In this section we summarise the setup and analysis in \cite{li}. We highlight that the techniques and terminology used here are entirely standard in the context of GSPT; we point to the usual references \cite{joneslec,kuehn}  for general definitions.

\subsection{The travelling wave equations}

The PDE \eqref{eq:master} is expressed in terms of the frame $(\zeta,t') = (x-ct,t)$ (where $c$ is a constant parameterizing the wavespeed) is given by

\begin{align}\label{eq:travelingwavemaster}
\bar{U}_{t'} &= c \bar{U}_{\zeta} +  \frac{\del}{\del \zeta} \left( D(\bar{U}) \frac{\del \bar{U}}{\del \zeta} \right) + R(\bar{U}) +\eps \left( \frac{\del^3 \bar{U}}{\del \zeta^2 \del t'} - c  \frac{\del^3 \bar{U}}{\del \zeta^3} \right).
\end{align}

Li et al. (2021) provide the precise conditions necessary on $D(\bar{U})$ and $R(\bar{U})$ for the existence of a shock-fronted travelling wave solution. We refer the reader to \cite{li} for an in-depth discussion about the modelling assumptions and potential generalisations underlying these definitions. Here we continue with the definitions used by  \cite{li} in their computations for the sake of consistency: we take a quadratic nonlinear diffusion term
\begin{align} \label{eq:diffusionterm}
D(\bar{U}) &= 6\left( \bar{U} - 7/12 \right)\left( \bar{U} -3/4 \right)
\end{align}
and cubic reaction term
\begin{align} \label{eq:reactionterm}
R(\bar{U}) &= 5\bar{U}(1-\bar{U})(\bar{U}-1/5).
\end{align}
We also record the following potential function (i.e. integral) of $D(\bar{U})$, which will be used to define the vector field in the travelling wave frame:
\begin{align} \label{eq:potentialterm}
F(\bar{U}) &= 2\bar{U}^3-4\bar{U}^2+\frac{21}{8}\bar{U}.
\end{align}

Note that the diffusion term is negative within the range $7/12 < \bar{U} < 3/4$, and the reaction term is pinned at $\bar{U} = 0$ and $\bar{U} = 1$. \\

A travelling wave solution is found as a steady-state to the above equation:

\begin{align}
-R(\bar{U}) &= \left( c \bar{U} +   \frac{\del F(\bar{U})}{\del \zeta}  -\eps c\frac{\del^2 \bar{U}}{\del \zeta^2} \right)_{\zeta}.
\end{align}

Letting 
\begin{align*}
\bar{P} &= c \bar{U} +   \frac{\del F}{\del \zeta}  -\eps c \frac{\del^2 \bar{U}}{\del \zeta^2} \\
\bar{V} &= F(\bar{U}) - \eps c \bar{U}_{\zeta},
\end{align*}

we arrive at the following {\it slow travelling wave equations} (where $\dot{~} := d/d\zeta$)

\begin{align}
\eps\dot{\bar{U}} &= \frac{1}{c} \left( F(\bar{U})-\bar{V}\right) \nonumber\\
\dot{\bar{P}} &= -R(\bar{U}) \label{eq:travelodesslow}\\
\dot{\bar{V}} &= \bar{P}-c\bar{U}. \nonumber
\end{align}

We also record the equivalent {\it fast travelling wave equations} in terms of the stretched variable $\xi := \zeta/\eps$:\footnote{From now on we use lowercase letters throughout to distinguish the variables and functions defined with respect to the stretched scaling, following the convention in \cite{GJ}.} 

\begin{align}
\bar{u}' &= \frac{1}{c} \left( F(\bar{u})-\bar{v} \right) \nonumber\\
\bar{p}' &= -\eps R(\bar{u}) \label{eq:travelodesfast}\\
\bar{v}' &= \eps(\bar{p}-c\bar{u}), \nonumber
\end{align}

where ${~}' := d/d\xi$. From the definitions of $R$ and $F$, the systems  \eqref{eq:travelodesslow}--\eqref{eq:travelodesfast}  evidently admit three fixed points for $\eps > 0$. Two of them are
\begin{equation} \label{eq:fps}
\begin{aligned}
\bar{z}^- &=  (\bar{u}^-, \bar{p}^-, \bar{v}^-) = (0, 0, 0) \\
\bar{z}^+ &=  (\bar{u}^+, \bar{p}^+, \bar{v}^+)  = (1, c, 5/8),
\end{aligned}
\end{equation}

lying on $S^{a,-}_0$ and $S^{a,+}_0$, respectively, and the third is the point $\bar{z}^{m} = (0.2,0.2c,0.381)$. Travelling waves of \eqref{eq:master} linking the states $\bar{U} = 0$ and $\bar{U} = 1$ correspond to heteroclinic connections between the fixed points $\bar{z}^-$ and $\bar{z}^+$.\\

 Since the systems \eqref{eq:travelodesslow}--\eqref{eq:travelodesfast} are odd-dimensional, it is convenient to choose an orientation that minimizes the dimension of the unstable manifolds $W^u(\bar{z}^-)$ and $W^u(\bar{z}^+)$. We therefore apply the orientation reversal  $\zeta \mapsto -\zeta$ and $\xi \mapsto -\xi$ in the subsequent analysis and work henceforth with the pair of systems
\begin{align}
\eps\dot{\bar{U}} &= \frac{1}{c} \left( \bar{V}-F(\bar{U})\right) \nonumber\\
\dot{\bar{P}} &= R(\bar{U}) \label{eq:travelslow}\\
\dot{\bar{V}} &= c\bar{U}  - \bar{P}\nonumber
\end{align}
and
\begin{align}
\bar{u}' &= \frac{1}{c} \left( \bar{v} - F(\bar{u}) \right) \nonumber\\
\bar{p}' &= \eps R(\bar{u}) \label{eq:travelfast}\\
\bar{v}' &= \eps(c\bar{u}-\bar{p}). \nonumber
\end{align}

Let us now assess the linear stability of the fixed points in \eqref{eq:fps}.

\begin{lem} \label{lem:eighierarchy}
Let $c > 0$ be fixed. For each sufficiently small $\eps > 0$, the fixed points $\bar{z}^-$ and $\bar{z}^+$ of \eqref{eq:travelfast} are both saddle-type equilibria. Furthermore, the three eigenvalues of the linearisation of \eqref{eq:travelfast} evaluated at both $\bar{z}^-$ and $\bar{z}^+$ are all real, having the hierarchy
\begin{align}
\mu_{f} \ll \mu_{s,1} < 0  < \mu_{s,2}.
\end{align}

Specifically, with respect to the scaling in \eqref{eq:travelfast} we have $\mu_{f} = \bigO(1)$ and $\mu_{s,i} = \bigO(\eps)$ where $i = 1,2$. 
\end{lem}

{\it Proof:} We verify the spectral hierarchy for $\bar{z}^-$; the steps for $\bar{z}^+$ are identical. Let $f(\bar{u},\bar{p},\bar{v})$ denote the vector field of \eqref{eq:travelfast}. The associated Jacobian matrix is
\begin{align} \label{eq:jacfast}
Df(\bar{u},\bar{p},\bar{v}) &= \begin{pmatrix} -\frac{D(\bar{u})}{c} & 0 & \frac{1}{c}\\
\eps R'(\bar{u}) & 0 & 0 \\
\eps c & -\eps & 0 \end{pmatrix}.
\end{align}

At $\bar{z}^-$ we have $D(0) = 21/8$ and $R'(0) = -1$. The characteristic polynomial of $Df(\bar{z}^-)$ is
\begin{align*}
p(\mu) &= -\mu^3- \frac{21}{8c}\mu^2 + \eps \mu + \frac{1}{c^2}\eps^2.
\end{align*}

Therefore,  the eigenvalue of largest magnitude has the expansion
\begin{align}
\mu_{f}(\eps) &= \frac{-21}{8c} + \bigO(\eps).
\end{align}

Thus $\mu_f(\eps)<0$ and uniformly bounded away from 0 for $\eps$ sufficiently small when $c > 0$. To estimate the remaining two $\bigO(\eps)$ eigenvalues $\mu_{s,i}$, $i=1,2$, we consider the scaling $\mu = \eps \nu$. Then
\begin{align*}
p(\nu) &= -\eps^3 \nu^3- \frac{21}{8c}\eps^2 \nu^2 + \eps^2 \nu + \frac{1}{c^2}\eps^2.
\end{align*}

Dividing a factor of $\eps^2$, regular perturbation theory again gives the solution of $p(\nu) = 0$ in orders of $\eps$ as
\begin{align}
\nu_{\pm} &= \frac{2}{21}\left( 2 \pm \sqrt{4 + \frac{42}{c^2}} \right) + \bigO(\eps).
\end{align}

Since the roots limit to distinct real values as $\eps \to 0$, the pair must perturb to distinct real values. Furthermore $\mu_{s,1} = \eps \nu_- < 0 < \eps\nu_+ = \mu_{s,2}$ for each $\eps>0$ sufficiently small. $\Box$\\

Following standard GSPT, we now define two subproblems which characterize the slow and fast singular limits of \eqref{eq:travelslow} and \eqref{eq:travelfast}. These are used to define a {\it singular heteroclinic orbit} by concatenating solutions of these subsystems. 

\subsection{Fast and slow singular limits}

We first characterize the fast dynamics by considering the singular limit $\eps \to 0$ of \eqref{eq:travelfast}.

\begin{defi} \label{def:critman}
The {\it layer problem} of \eqref{eq:travelfast} is given by
\begin{align}
\bar{u}' &= \frac{1}{c} \left( \bar{v} - F(\bar{u}) \right) \nonumber\\
\bar{p}' &= 0 \label{eq:layerproblem}\\
\bar{v}' &= 0. \nonumber
\end{align} 
\end{defi}

In this limit, the slow variables $\bar{p},\bar{v}$ are constant, parametrising the {\it fast fiber bundle}; away from the zero set of $u'$, the layer problem then specifies the one-dimensional fast dynamics fiberwise.

\begin{figure}[t!] 
\centering
(a)\includegraphics[width=0.95\textwidth]{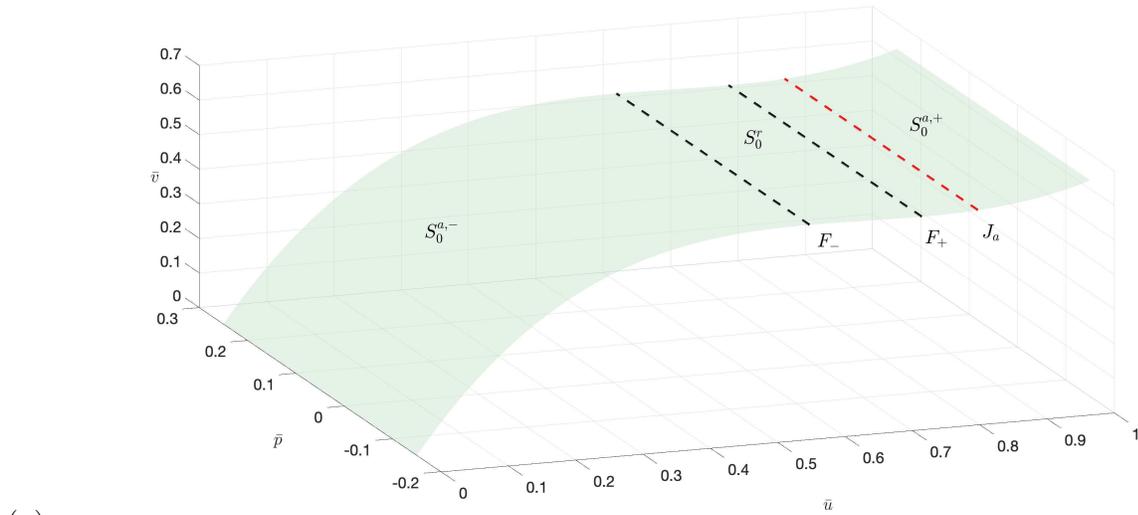}\\
(b)\includegraphics[width=0.95\textwidth]{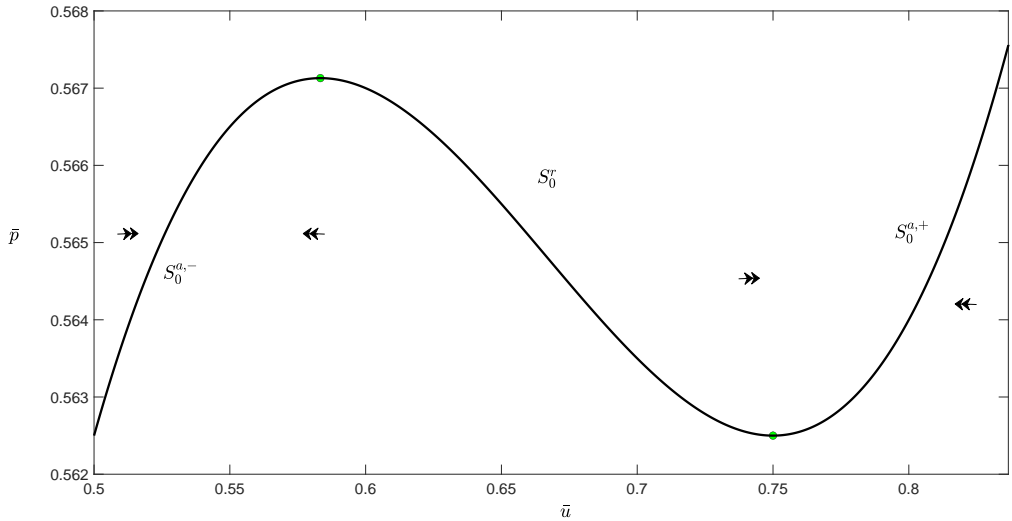}
\caption{(a) Depiction of the critical manifold $S_0$ from Def. \ref{def:critman} in $(\bar{u},\bar{p},\bar{v})$-space. (b) Dynamics of the layer problem \eqref{eq:travelfast} away from $S_0$, projected onto the $(\bar{u},\bar{p})$ plane.}
\label{fig:layerproblem}
\end{figure}

\begin{defi} \label{def:critman}
The {\it critical manifold} is $S_0 = \{(\bar{u},\bar{p},\bar{v}):\,\bar{v} = F(\bar{u})\}$. We decompose $S_0$ into the pieces
\begin{align*}
S_0 &= S^{a,-}_0 \cup F_- \cup S^r_0 \cup F_+ \cup S^{a,+}_0,
\end{align*}

with
\begin{align*}
S^{a,-}_0 &= S_0 \cap \{\bar{u} < 7/12\}\\
F_- &= S_0 \cap \{\bar{u} = \bar{u}_F := 7/12\}\\
S^{r}_0 &= S_0 \cap \{7/12 < \bar{u} < 3/4\}\\
F_+ &= S_0 \cap \{\bar{u} = 3/4\}\\
S^{a,+}_0 &= S_0 \cap \{\bar{u} > 3/4\}.
\end{align*}

Here, $F_-$ and $F_+$ denote resp. the left and right folds of the cubic manifold $S_0$, where $D(\bar{U})$ vanishes. We also define the {\it jump-on curve} $J_a \subset S^{a,+}_0$, which is the horizontal projection of the fold $F_-$ onto $S^{a,+}_0$:
\begin{align*}
J_a &= S_0 \cap \{\bar{u} = \bar{u}_J := 5/6\}.
\end{align*}
These sets are drawn in Fig. \ref{fig:layerproblem}(a).
\end{defi}
 The Jacobian matrix of \eqref{eq:layerproblem} has two trivial zero eigenvalues for each $x \in S_0$, corresponding to the locally two-dimensional tangent space at each point.  The sign of the remaining {\it nontrivial eigenvalue} determines the dynamics along the fast fibers nearby.

\begin{defi}
We say that $S_0$ is {\it normally hyperbolic} at $x \in S_0$ if its nontrivial eigenvalue $\lambda$ does not lie on the imaginary axis. Furthermore, we say that $S_0$ is {\it normally hyperbolic attracting} if $\lambda <0$ and {\it normally hyperbolic repelling} if $\lambda > 0$.
\end{defi}

The nontrivial eigenvalue of the linearization of \eqref{eq:layerproblem} is 
\begin{align*}
\lambda &= -\frac{D(\bar{u})}{c},
\end{align*}
from which it immediately follows that $S^{a,-}_0$ and $S^{a,+}_0$ are normally hyperbolic attracting and $S^r_0$ is normally hyperbolic repelling when $c > 0$;  see Fig. \ref{fig:layerproblem}(b). Well-known theorems of Fenichel \cite{fenichel} specify the existence of normally hyperbolic invariant {\it slow manifolds} $S^{a,\pm}_{\eps}$ and $S^r_{\eps}$  and fast fiber bundles near to compact normally hyperbolic subsets of their singular counterparts $S^{a,\pm}_0$ resp. $S^r_0$. \\

The critical manifold loses normal hyperbolicity via a simple zero eigenvalue crossing at the folds $F_{\pm}$.  Slow dynamics may be defined on $S_0$ away from $F_{\pm}$ by considering the singular limit of \eqref{eq:travelslow}. Differentiating both sides of $\bar{V} = F(\bar{U})$ with respect to $\zeta$ we have
\begin{align} \label{eq:UVrelation}
\dot{\bar{V}} &= D(\bar{U})\dot{\bar{U}}.
\end{align}

Solutions of the reduced problem blow up as they approach $F_{\pm}$, and so we rescale $\zeta$ by the function $D(\bar{U})$ to extend solutions of the reduced problem smoothly across the folds. Altogether, we define the following two systems. 
\begin{defi} \label{eq:reducedproblem}
The {\it reduced problem} defined on $S_0$ is given by

\begin{equation} \label{eq:reducedproblem}
\begin{aligned} 
D(\bar{U})\dot{\bar{U}} &= c\bar{U}-\bar{P}\\
\dot{\bar{P}} &= R(\bar{U}).
\end{aligned}
\end{equation}

Rescaling $\zeta$ by the diffusion $D(\bar{U})$, we obtain the {\it (desingularized) slow flow}

\begin{equation} \label{eq:slowflow}
\begin{aligned} 
\dot{\bar{U}} &= c\bar{U}-\bar{P}\\
\dot{\bar{P}} &= R(\bar{U}) D(\bar{U}).
\end{aligned}
\end{equation}
~\\
\end{defi}

The slow flow extends the solutions of the reduced problem across the lines of singularities given by $F_{\pm}$ in finite time. Note that although topological equivalence of solutions of \eqref{eq:reducedproblem} and \eqref{eq:slowflow} is assured on $S^{a,\pm}_{0}$ and $S^r_0$, there is an orientation reversal on $S^{r}_0$ owing to the change in sign of $D(\bar{U})$.\\

The slow flow \eqref{eq:slowflow} admits five fixed points corresponding to the zeroes of $R(\bar{U})$ and $D(\bar{U})$. A calculation of the eigenvalues of the Jacobian of \eqref{eq:slowflow} shows that the fixed points 
\begin{equation} \label{eq:fpsreduced}
\begin{aligned}
\bar{Z}^- &= (0,0)\\
\bar{Z}^+ &= (1,c)
\end{aligned}
\end{equation}

are saddle-type when $c > 0$. These correspond to the singular limits of the fixed points \eqref{eq:fps} of the full system as $\eps \to 0$. 

\begin{remk}
The set $F_{-}$ consists mostly of so-called {\it regular} (or generic) {\it fold points}; we refer the reader to the precise definition in \cite{kuehn}.\footnote{There are also isolated {\it canard points}, corresponding to folded singularities where trajectories beginning near an attracting branch of a slow manifold are able to cross over to the repelling branch, tracking it for  $\bigO(1)$ time. In the travelling wave literature such points are referred to as `holes in the wall,'  relating to the observation that $D(\bar{U}) = 0$ defines `walls' of singularities corresponding to the fold lines; see \cite{harley2}, \cite{li2}, and \cite{wexpet10}.} One of the seminal achievements of GSPT is to rigorously characterise the flow of typical orbits across folds, where Fenichel theory breaks down. This requires extending the slow manifolds (which a priori exist only over compact normally hyperbolic subsets of the critical manifold) across such folded singularities using geometric blow-up theory; see e.g. \cite{krupa,relaxoscr3,aproposcanards}. 
\end{remk}

  We are interested in the existence of a {\it singular heteroclinic orbit} that connects $\bar{Z}^-$ to $\bar{Z}^+$ by jumping across a regular fold point on $F_-$. This singular orbit is a hybrid object defined by concatenating solutions of \eqref{eq:layerproblem} and \eqref{eq:slowflow}. 

\begin{hyp} \label{hyp:singhet}
There exists a {\it singular heteroclinic orbit} $\Gamma_0$ for $c = c_0>0$, connecting $\bar{U} = 0$ to $\bar{U} = 1$. This orbit is formally defined as the concatenation of the following solution segments:
\begin{itemize}
\item the portion of the unstable manifold $W^u(\bar{Z}^-)$ which connects $\bar{Z}^-$ to a regular fold point $(7/12,\bar{p}_f)$ on $F_-$, corresponding to a unique trajectory $X_R(\zeta)$ of \eqref{eq:slowflow} defined on $\zeta \in (-\infty,0]$
\item the horizontal fast fiber connecting $W^u(\bar{Z}^-) \cap F_-$ to a point $(5/6,\bar{p}_f,\bar{v}_f)$ on the jump curve $J_a$, corresponding to a unique trajectory $x_R(\xi)$ of  \eqref{eq:layerproblem} defined on $\xi \in (-\infty,\infty)$, and 
\item the portion of the stable manifold $W^s(\bar{Z}^+)$ which connects $(5/6,\bar{p}_f)$ to $\bar{Z}^+$, corresponding to a unique trajectory $X_R(\zeta)$ of \eqref{eq:slowflow} defined on $\zeta \in [0,\infty)$. We also suppose that the slow flow is transverse to $J_a$ at $W^s(\bar{Z}^+)\cap J_a$.\\
 \end{itemize}
 
We furthermore suppose that the heteroclinic connection $\Gamma_0$ is transversal with respect to variation in $c$; i.e. extending the reduced systems with the equation $c' = 0$, we have $\Gamma_0 = W^u(\bar{Z}^-,c) \transv W^s(\bar{Z}^+,c)$, where $W^u(\bar{Z}^-,c) = \cup_c W^u(\bar{Z}^-(c))$ and $W^s(\bar{Z}^+,c) = \cup_c W^s(\bar{Z}^+(c))$ denote, respectively, the extended two-dimensional center-unstable (resp. three-dimensional center-stable) manifolds, continued across the fast jump in the usual way.\\
\end{hyp}

\begin{figure}[t!] 
\centering
\includegraphics[width=1.05\textwidth]{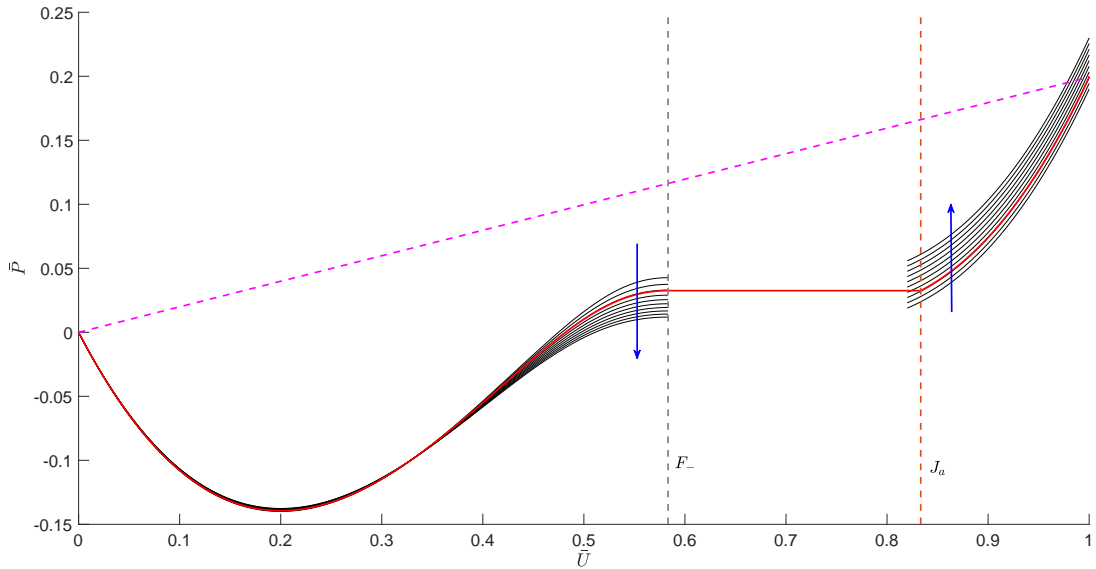}
\caption{Segments of $W^u(\bar{Z}^-)$ and $W^s(\bar{Z}^+)$ arising from the slow flow \eqref{eq:slowflow}--\eqref{eq:fpsreduced}, as $c$ is varied within the range $[0.19,0.23]$. The stable manifolds are extended slightly past $J_a$ to illustrate the transversality of the slow flow at $J_a$. The singular homoclinic orbit $\Gamma_0$ is depicted by the red curve for $c_0 \approx 0.199362$. The blue arrows indicate how the (un)stable manifolds vary as $c$ is increased. The magenta dashed line denotes the $\bar{U}$-nullcline for $c = c_0$.}
\label{fig:transversehet}
\end{figure}

Hypothesis \ref{hyp:singhet} can be readily verified numerically; this was done already in \cite{li}. Here we report a singular heteroclinic connection for $c_0 \approx 0.199362$ which is clearly transversal with respect to variation in $c$; see Fig. \ref{fig:transversehet}. Based on the relative location of $\Gamma_0$ with respect to the $\bar{U}$-nullcline, we also record the following monotonicity hypothesis.

\begin{hyp} \label{hyp:monotoneU}
The portions of $\Gamma_0$ (from Hypothesis \ref{hyp:singhet}) lying in $S^{a,\pm}_0$ are strictly monotone in $\bar{U}$; i.e.  $\dot{\bar{U}}(\zeta) > 0$ (equivalently, $\dot{\bar{V}}(\zeta) > 0$) for all $\zeta < 0$ on $S_0^{a,-}$ and for all $\zeta > 0$ on $S_0^{a,+}$.
\end{hyp}

This monotonicity condition turns out to be useful when analyzing the reduced eigenvalue problem introduced in Sec. \ref{sec:slowevans}: it is equivalent to the statement that the projectivized dynamics of the slow coordinates lives entirely on one chart of $\mathbb{CP}$ specified by $S = (P/V)$, $V \neq 0$, when the temporal eigenvalue parameter $\lambda = 0$. Here $P$ and $V$ denote respectively the linearisations of $\bar{P}$ and $\bar{V}$.\footnote{Our notation  is equivalent to the ``$\delta$'' convention   to denote linearized variables, i.e. $V = \delta\bar{V}$.}

\subsection{Travelling waves for $0 < \eps \ll 1$}

We now use hypothesis \ref{hyp:singhet}, to show the existence of a one-parameter family of travelling waves $\{\Gamma_\eps\}$ for sufficiently small values of $\eps > 0$; we write it in such a way as to make clear the uniform convergence onto the singular wave $\Gamma_0$.

\begin{lem} \label{hyp:fullhet}
For any $\delta > 0$, there exists $\bar{\eps} > 0$ such that for each $\eps \in (0,\bar{\eps})$, there exists a wavespeed $c(\eps)$ and a heteroclinic orbit $\Gamma_{\eps}(\xi)$ for \eqref{eq:travelfast} connecting $\bar{z}^-$ to $\bar{z}^+$ with $d_H(\Gamma_{\eps},\Gamma_0) < \delta$ and $|c(\eps) - c_0| < \delta$. 
\end{lem}

\begin{figure}[t!] 
\centering
\includegraphics[width=1.0\textwidth]{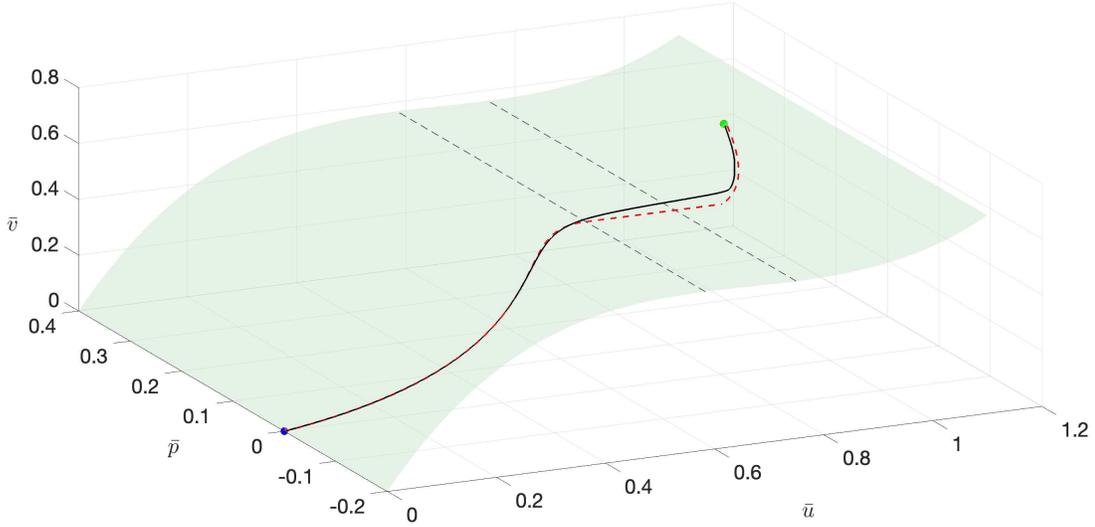}
\caption{ Numerical approximation of a heteroclinic orbit (black solid curve) connecting the saddle equilibria $(0,0,0)$ and $(1,c,F(1))$ for $(\eps,c) = (0.001, 0.20637)$, together with the singular heteroclinic orbit  $\Gamma_0$ (red dashed curve) from Hypothesis \ref{hyp:singhet}.}
\label{fig:fullhet}
\end{figure}

Shock-fronted travelling waves for small fixed values of $\eps > 0$ are numerically approximated in \cite{li}; see also Fig. \ref{fig:fullhet}. We sketch a geometrical proof showing that Hypothesis \ref{hyp:singhet} implies that Lemma \ref{hyp:fullhet} holds. Our sketch follows standard arguments for constructing heteroclinic orbits by estimating the relative orientations of slow manifolds and their fast fiber bundles for $\eps > 0$ small,  see eg. \cite{szmolyan}. Here we also address the issue that the jump occurs across a regular fold, necessitating a slightly different estimate.\\

{\it Proof:} Fix $\alpha > 0$ small enough that $\Gamma_0$ is the unique transverse singular heteroclinic orbit for $c \in B_{\alpha} := (c_0 - \alpha, c_0 + \alpha)$. We also fix a section which intersects $\Gamma_0$ in its interior in the middle of the jump, say $\Sigma = \{\bar{u} = 0.7\}$. For $\eps > 0$ sufficiently small, Fenichel theory provides the existence of slow manifolds $S^{a,-}_{\eps}$ and $S^{a,+}_{\eps}$, which are $\bigO(\eps)$-close in Hausdorff distance to compact normally hyperbolic subsets of $S^{a,-}_0$ and $S^{a,+}_0$.\\

 The slow manifold $S^{a,-}_{\eps}$ can be continuously extended across a regular fold of a two-dimensional critical manifold in $\mathbb{R}^3$ using blow-up theory; after crossing the fold, $S^{a,-}_{\eps}$ lies $\bigO(\eps^{2/3})$ close to the singular fast fiber subbundle $\mathscr{J}$ extending from the fold $F_-$ to the jump-on curve $J_a$, at the section $\Sigma$ (see Theorem 1 in \cite{relaxoscr3} for details). The slow manifold $S^{a,-}$ then flows along the fast fiber bundle, reaching $\Sigma$ in $\bigO(1)$ time with respect to the scaling \eqref{eq:travelfast} (this intersection must exist for $\eps > 0$ sufficiently small since $\Gamma_0$ intersects $\Sigma$ by construction). Thus, $S^{a,-}_{\eps}$ remains $\bigO(\eps^{2/3})$-close to $\mathscr{J}$ at the point of intersection on $\Sigma$ . \\
 
 On the other hand, consider a ball $V$ inside $W = \cup_{c \in B_{\alpha}}W^s(\bar{z}^-(c))$ which straddles $J_a$. Fenichel theory implies that the (nonlinear) fast fiber bundle over the basepoint set $V$ perturbs to an $O(\eps)$-close fast fiber bundle over the corresponding basepoints now given on $S^{a,+}_{\eps}$ \cite{fenichel}. Extend this fiber bundle backwards to $\Sigma$. Therefore at $\Sigma$, the transverse intersection of the projection of $F$ with the stable fast fiber bundle over $S^{a,+}_0$ persists over a small range in $\eps$, and for each such $\eps$ the intersection occurs at a corresponding wavespeed $c(\eps)$. Evaluated at $\Sigma$, the perturbed intersection is $\bigO(\eps^{2/3})$-close to the original intersection in the extended $(\bar{u},\bar{p},\bar{v},c)$-space, i.e. $c(\eps) \to c_0$ as $\eps \to 0$.   $\Box$\\

\begin{remk}
We highlight a high-dimensional generalization of this construction, due to Lin and Wechselberger \cite{lin}. They prove a generalization of Theorem 1 in \cite{relaxoscr3}, and they also require exchange/lambda lemma-type estimates to pick out unique waves from two-parameter families. In our case, the transversal intersection $\Gamma_{\eps}$ for each corresponding wavespeed $c(\eps)$ must be unique, which can be verified directly with dimension counting.
\end{remk}

We are now able to state the main theorem proved in this paper.

\begin{thr} \label{thm:main}
Suppose Hypotheses \ref{hyp:singhet}, and \ref{hyp:monotoneU} hold. Fix a contour $K$ in the right half of the complex plane. Then there exists $\bar{\eps} > 0$ so that for each $0 < \eps \leq \bar{\eps}$, the wave $\Gamma_{\eps}$ does not have spectrum inside of $K$.  
\end{thr}

The proof of this theorem is developed mostly in the second half of the paper; see Sec. \ref{sec:outline} for a general outline of the strategy. A characterization of the spectral stability of a travelling wave must first be developed. This is accomplished by the end of Sec. \ref{sec:bundles}. \\

We finish this section with a preliminary estimate showing that we have strong control over the dynamics of the wave  via the reduced subsystems away from the transition sets $F_-$ and $J_a$; c.f. Corollary 2.2 in \cite{GJ}. 

\begin{coro} \label{cor:waveestimate}
Given $\eps > 0$ sufficiently small, let $X(\zeta,\eps) = x(\xi,\eps)$  parametrize the travelling wave solution from Hypothesis \ref{hyp:fullhet} corresponding to the wavespeed $c = c(\eps)$. The following are true for each fixed $a > 0$:
\begin{align*}
\lim_{\eps \to 0} x(\xi,\eps) &= x_R(\xi) \,\text{ uniformly for } |\xi| \leq a\\
\lim_{\eps \to 0} x(\zeta,\eps) &= X_R(\zeta) \text{ uniformly for } |\zeta| \geq a.
\end{align*}

In the second limit, we take $X_R$ to mean the embedding of the corresponding parametrised solutions of $\Gamma_0|_{S^{a,-}_0 \cup S^{a,+}_0}$ into $\mathbb{R}^3$ via the map $(\bar{U},\bar{P}) \mapsto (\bar{U},\bar{P},F(\bar{U})) \in S_0$. 
\end{coro}

{\it Proof:} We sketch the argument for the first uniform limit.  Fix an $a > 0$. Then $|x(\xi,\eps) - x_R(\xi)|$ can be uniformly bounded through a combination of Gronwall's inequality and the distance estimate in Hypothesis \ref{hyp:fullhet}. The argument for the second uniform limit is similar. $\Box$

\section{Geometric preliminaries} \label{sec:geometry}

Ahead of our analysis in the remaining sections, we explicitly describe the geometric spaces in which our objects of interest live. We will track the dynamics of fibers of complex vector bundles,  corresponding to points in $\mathbf{Gr}(k,n)$, the {\it Grassmannian} of complex $k$-planes in $\mathbb{C}^n$. The general set-up is standard: we will study dynamics on the Grassmannian indirectly via embedded systems on {\it complex projective space}.\\

We first introduce some essential notation for complex projective space. The space $\mathbb{CP}^n$ is the quotient space $(\mathbb{C}^{n+1}-\{0\})/\sim{ }$  subject to the equivalence relation that identifies complex rays, i.e. $y_1,y_2 \in \mathbb{C}^{n+1}-\{0\}$ are identified if $y_1=\alpha y_2$ for some $\alpha \in \mathbb{C}$. We denote the projection map by $\pi: \mathbb{C}^{n+1}-\{0\} \to \mathbb{CP}^n$, and for each $y \in \mathbb{C}^n$ we use the notation $[y]$ or $\hat{y}$ for the image $\pi(y) \in \mathbb{CP}^n$ (this notation also holds for sets $S \subset \mathbb{C}^n-\{0\}$). \\

Later, in Section \ref{sec:estimate1} we will use the following elementary lemma to interpret metric neighborhoods on $\mathbb{CP}^{2}$ in terms of those on $\mathbb{C}^{3}$.

\begin{lem} \label{lem:metbound}
Let $\rho$ be any metric on $\mathbb{CP}^n$.  Then there is a constant $K>0$, depending only on $\rho$, such that the following is true for each $\delta > 0$:\\

If $y \in \mathbb{C}^{n+1}$ is a complex vector with unit modulus and if $P \subset \mathbb{C}^{n+1}$ is a complex $n$-plane passing through the origin with $\hat{y} \notin N_{\delta}(\hat{P})$, then $$|w| \geq K \delta,$$ where $w$ is the component of $y$ lying orthogonal to $P$, and $N_{\delta}(\hat{P})$ denotes the $\delta$-neighborhood of $\hat{P}$ with respect to $\rho$, i.e. $
N_{\delta}(\hat{P}) = \{\hat{x} \in \mathbb{CP}^n: \rho(\hat{x},\hat{p}) < \delta \text{ for some } \hat{p} \in \hat{P}\}. $
\end{lem}

\begin{remk}
The norm $|x| := \sqrt{\bar{x}^T x}$ in Lemma \ref{lem:metbound} is the usual Euclidean norm, and the orthogonal component is defined using the Hermitian inner product $x\cdot_{\mathbb{C}}y := \sum_{i=1}^n x_i \bar{y}_i$. Essentially the lemma says that if a vector is inclined far away from a plane, then the component of the vector lying orthogonal to the plane cannot be too small, up to a scaling factor that depends on the details of how this inclination is measured. This lemma was implicitly used in Lemma 5.4 in \cite{GJ} for the case $n = 3$, so we provide an explicit proof for arbitrary $n$.
\end{remk}

{\it Proof of Lemma \ref{lem:metbound}:} Suppose that there is no such uniform constant for some fixed choice of $\delta > 0$. Then we can find a sequence $\{y_k\}$ lying on the unit sphere in $\mathbb{C}^{n+1}$ such that for each $k$, we have $\hat{y}_k \notin N_{\delta}(\hat{P})$ and $|w_k| < \delta/k$, where $w_k$ denotes the orthogonal projection of $y_k$ with respect to $P$. Choose a convergent subsequence $y_k \to y_*$, relabeling $k$ if necessary. Then $y_*$ clearly has zero component in the direction orthogonal to $P$, implying that $y_* \in P$ and hence that $\hat{y}_* \in N_{\delta}(\hat{P})$. On the other hand, we have $\hat{y}_k \to \hat{y}_* \notin N_{\delta}(\hat{P})$ since the quotient map is continuous and the complement of $N_{\delta}(\hat{P})$ is closed, giving a contradiction. $\Box$\\

 Let $v_1,\cdots,v_k \in \mathbb{C}^n$ denote a spanning set of vectors for the $k$-plane $W$. The well-known {\it Pl\"{u}cker embedding} 
\begin{align*}
\mathbf{Gr}(k,n) &\xhookrightarrow{\psi} \mathbb{P}(\bigwedge^{k} \mathbb{C}^n)  \cong \mathbb{CP}^{\binom{n}{k}-1}\\
\psi(W) &:= [v_1 \wedge \cdots \wedge v_k]
\end{align*}

provides homogeneous coordinates for the Grassmannian, called Pl\"{u}cker coordinates. The image $\psi(W)$ can be coordinatized by writing down the maximal minors/cofactors of the $(n\times k)$ matrix  $M_W = [v_1~ \cdots~ v_k]$.\\

{\it Example.} Consider a complex 2-plane $W \in \mathbf{Gr}(2,3)$ with spanning complex vectors $v_1 = (a_1,b_1,c_1)^T$ and $v_2 = (a_2,b_2,c_2)^T$. A coordinate representation of $[v_1 \wedge v_2] \in \mathbb{C}\mathbb{P}^2$ is given by the list of $2\times 2$ cofactors of the $3\times 2$ matrix $[v_1~~v_2]$,
namely
\begin{align*}
N_W := \begin{pmatrix}  b_1 c_2 - c_1 b_2 \\ a_2 c_1 - a_1 c_2 \\ a_1 b_2 - a_2 b_1\end{pmatrix}.
\end{align*}
Observe that the Pl\"ucker coordinate representation $N_W$ can be interpreted as a parametrisation of the plane $W$ by a (nontrivial) normal vector.

\begin{remk}
The Pl\"{u}cker coordinates are generally redundant: we have $\binom{n}{k}$ of them, whereas the dimension of $\mathbf{Gr}(k,n)$ is $k(n-k)$. In this paper our analysis is restricted to the cases $\mathbf{Gr}(1,3)$ and $\mathbf{Gr}(2,3)$, since we will be occupied with tracking complex line (resp. plane) bundles. In the case $\mathbf{Gr}(1,3)$, the redundancy arises from nonunique scalings of the complex ray. In the case $\mathbf{Gr}(2,3)$, an analogous scaling nonuniqueness arises when choosing a normal vector to represent a complex plane (see the example above). 
 \end{remk}

The Grassmannian spaces are metrizable. Here we construct an explicit metric, because we will need a way to measure `angles' between lines and planes. Using a complex version of the construction given in \cite{jonestin}, we take advantage of the Pl\"{u}cker embedding to equip $\mathbf{Gr}(k,n)$ with the {\it Fubini-Study} metric. Complex projective space can be embedded into the space of Hermitian $(n+1)\times(n+1)$ matrices $\mathcal{H}$ via the following isometric immersion:
\begin{align*}
\mathbb{CP}^n &\xhookrightarrow{\varphi} \mathcal{H}\\
\varphi([x]) &:= \frac{x\bar{x}^T}{|x|^2},
\end{align*}

where $x \in \mathbb{C}^{n+1}-\{0\}$ represents the ray $\mathbb{C}x$. The space of Hermitian matrices is in turn endowed with an inner product given by
\begin{align*}
d(P,Q) &= \frac{1}{2} \text{Tr}(PQ),
\end{align*}

which allows us to define a norm $\norm{P} = d(P,P)^{1/2}$, and hence a metric induced by the norm given by $\text{dist}(P,Q) := \norm{P - Q}$. We can pull back this metric on $\mathcal{H}$ to one on $\mathbb{CP}^n$, and hence on the Grassmannian:
\begin{align*}
\text{dist}([x],[y]) &:= \norm{ \frac{x\bar{x}^T}{|x|^2} - \frac{y\bar{y}^T}{|y|^2}}\\
&= \sqrt{\frac{1}{2} \text{Tr} \left[  \left( \frac{x\bar{x}}{|x|^2} -\frac{y\bar{y}^T}{|y|^2}  \right)^2 \right]}\\
&= \sqrt{1- \frac{|(x\cdot_{\mathbb{C}} y)|^2}{|x|^2 |y|^2}}.
\end{align*}

We observe that $\text{dist}([x],[y])$ recovers the sine of the so-called {\it Hermitian angle} \cite{complexangles} between a corresponding pair of nontrivial complex vectors representing $[x],\,[y]$. For example, this metric on $\mathbf{Gr}(2,3)$ measures the Hermitian angle between a pair of normal vectors parametrising a corresponding pair of complex $2$-planes.\\

We will use the Pl\"{u}cker embedding to study (linear) dynamics which describe the evolution of $k$-planes in the Grassmannian. Concretely, let $y' = b(\xi) y$ denote a nonautonomous linear system on $\mathbb{C}^{n+1}$. This linear system induces a nonlinear flow on $\mathbb{CP}^n$, and hence on the corresponding Grassmannian, which we denote by 
\begin{align} \label{eq:projdyn}
\hat{y}' &= \hat{b}(\hat{y},\xi).
\end{align}

We mention two elementary facts about about these induced flows. The first is that $y$ is an eigenvector of $b(\xi)$ if and only if $\hat{y}$ is a critical point of $\hat{b}(\hat{y},\xi)$. The second is that if $\{\lambda_i\}$ denotes the eigenvalues of the linearization of a (constant) matrix $b$ and if $y_i$ is an eigenvector corresponding to the eigenvalue $\lambda_i$, then the eigenvalues of the linearization of $\hat{b}$ near the corresponding critical point $\hat{y}_i$ are $\lambda_j - \lambda_i$ for $j \neq i$. As a consequence of these facts, if the largest eigenvalue of a linear system is simple, then the corresponding eigendirection in the projectivized system is a {\it stable} fixed point.\\

We will find it useful to study the system \eqref{eq:projdyn} in its own right as a nonlinear flow on appropriate charts on $\mathbb{CP}^n$, but it is also convenient to treat the natural system $y' = b(\xi) y$ as the Pl\"{u}cker coordinate representation of the system on projective space. To clarify this, we introduce the notion of {\it induced variational fields}. We repeat the general treatment of \cite{jonestin}, Sec. 3.2. Let $V$ be a vector space of dimension $n$ and let $1 \leq \sigma \leq n$. Given a linear map $A: V \to V$, an induced derivation of order $\sigma$ is the multilinear map $A^{\sigma}: \bigwedge^{\sigma} V \to \bigwedge^{\sigma} V$ defined by
\begin{align} \label{def:Asigma}
v_1 \wedge \cdots \wedge v_{\sigma} \mapsto \sum_{i=1}^{\sigma} v_1 \wedge \cdots \wedge A v_i \wedge \cdots \wedge v_{\sigma}.
\end{align}
We have
\begin{lem} \label{lem:twowedgeeigs}
Let $\{\lambda_i\}_{i=1}^n$ denote the set of eigenvalues of $A^{(1)} = A$. Then the set of eigenvalues of $A^{(2)}$ is given by $\{\lambda_i + \lambda_j\}_{i<j,i=1}^{i=n}$.
\end{lem}

{\it Proof:} Choose an ordered basis for $V$ consisting of eigenvectors of $A$, and use it to construct an ordered basis for $\bigwedge^2 V$. The conclusion follows from a straightforward calculation using the definition \eqref{def:Asigma}. $\Box$\\

If $x' = Ax$ is a linear vector field, then $P' = A^{(\sigma)}P$ defines a linear vector field on the $\sigma$th exterior power space. By using the Pl\"{u}cker embedding, the latter induced vector field can be interpreted as evolving points in $\mathbf{Gr}(\sigma,n)$ (see Sec. 3.2 in  \cite{jonestin} for a complete description of this construction in the real setting; the extension to complex $k$-planes is direct). \\

\section{Spatial eigenvalue problems} \label{sec:spatialeigs}
We now assume that Hypotheses \ref{hyp:singhet},  \ref{hyp:monotoneU}, and  \ref{hyp:fullhet} hold, and we proceed to define the spatial eigenvalue problem associated with the one-parameter family of waves $(\Gamma_{\eps},c({\eps}))$.  Let us fix $\eps > 0$ sufficiently small and represent $\Gamma_{\eps}$  by the solution $\bar{X}(\zeta,\eps) = (\bar{U}(\zeta,\eps),\bar{P}(\zeta,\eps),\bar{V}(\zeta,\eps))$. We consider perturbations of the form
\begin{align}
\tilde{U}(\zeta,t) &= \bar{U}(\zeta,\eps) + \delta e^{\lambda t} U(\zeta)+\bigO(\delta^2).
\end{align}

Assuming that $\tilde{U}(\zeta,t)$ solves \eqref{eq:master} and applying the fact that $\bar{U}(\zeta,\eps)$ is a stationary solution of \eqref{eq:travelingwavemaster}, we arrive at the eigenvalue problem by collecting $\bigO(\delta)$ terms:
\begin{align} \label{eq:fulleigprob}
cU_{\zeta} + (D(\bar{U})U)_{\zeta\zeta} + R'(\bar{U})U + \eps (\lambda U_{\zeta \zeta}-cU_{\zeta\zeta\zeta}) &= \lambda U.
\end{align}

Collecting derivatives and defining
\begin{align*}
V &= D(\bar{U})U + \eps (\lambda U - cU_{\zeta})\\
P &= cU + V_{\zeta},
\end{align*}

we obtain the following closed three-dimensional system: 
\begin{align}
\eps \dot{U} &= \frac{1}{c} \left((\eps \lambda + D(\bar{U}))U  - V\right) \nonumber\\
\dot{P} &= (\lambda-R'(\bar{U}))U \label{eq:visclinslow}\\
\dot{V} &= P-cU.\nonumber
\end{align}
 
 In the stretched scaling $\xi = \zeta/\eps$, we obtain 
\begin{align}
u' &= \frac{1}{c} \left( -v+(\eps \lambda + D(\bar{u}))u \right) \nonumber\\
p' &=\eps (\lambda-R'(\bar{u}))u \label{eq:visclinfast}\\
v' &=\eps (p-cu). \nonumber
\end{align}

In view of the orientation reversal giving the pair of wave equations \eqref{eq:travelslow}--\eqref{eq:travelfast}, we will instead work with the following pair of linearized systems throughout the paper:
\begin{align}
\eps \dot{U} &= \frac{1}{c}(V-(\eps \lambda+D(\bar{U}))U) \nonumber\\
\dot{P}&= (R'(\bar{U})-\lambda) U  \label{eq:linreverseslow}\\
\dot{V} &= cU-P. \nonumber
\end{align}
and
\begin{align}
u' &= \frac{1}{c}(v-(\eps \lambda+D(\bar{u}))u) \nonumber\\
p' &= \eps (R'(\bar{u})-\lambda) u  \label{eq:linreversefast}\\ 
v' &= \eps(cu-p). \nonumber
\end{align}

\begin{remk}
Setting $\lambda = 0$ in \eqref{eq:linreverseslow}--\eqref{eq:linreversefast} gives the standard variational equations of \eqref{eq:travelslow}--\eqref{eq:travelfast} along the wave $(\Gamma_{\eps},c(\eps))$ for each $\eps > 0$. In particular, we can account for the existence of the translational eigenvalue; the variational equations are satisfied by a nontrivial uniformly bounded solution, namely the derivative of the wave (i.e. the vector field).
\end{remk}

Collecting $y = (u,p,v)$ and $Y = (U,P,V)$, we write the pair of systems \eqref{eq:linreverseslow}--\eqref{eq:linreversefast}  more compactly in terms of the matrices $A,a$:
\begin{align}
\dot{Y} &= A(\bar{U},\lambda,\eps)Y \label{eq:bigaeps}\\
y' &= a(\bar{u},\lambda,\eps)y. \label{eq:littleaeps}
\end{align}

We end this section by recording the {\it fast projectivized eigenvalue problem} $\hat{y} = \hat{a}(\bar{u},\lambda,\eps)$ corresponding to the linear eigenvalue problem \eqref{eq:linreversefast}. This system is defined on the space $\mathbb{CP}^2$; with respect to the chart $(\beta_1,\beta_2) = (p/u,v/u)$ with $u \neq 0$, we have

\begin{equation} \label{eq:projsystemu0}
\begin{aligned}
\beta_1' &= \frac{1}{c}(D(\bar{u}) \beta_1- \beta_1 \beta_2) +  \eps (R'(\bar{u}) + \lambda(\beta_1/c - 1))\\
\beta_2' &= \frac{1}{c}(D(\bar{u}) \beta_2 - \beta_2^2) + \eps\left(c-\beta_1 + \frac{\lambda}{c}\beta_2\right).
\end{aligned}
\end{equation}

{\it Slow} projectivized eigenvalue problems can be similarly defined using \eqref{eq:linreverseslow}. We record for future convenience the following system defined with respect to the chart $(\beta_1,\beta_2) =(U/V,P/V)$ with $V \neq 0$:
\begin{equation} \label{eq:projsystemu0slow}
\begin{aligned}
\eps \dot{\beta_1} &= \frac{1}{c} \left( 1-D(\bar{U}) \beta_1 \right) + \eps (\beta_1 \beta_2 - c \beta_1^2 - (\lambda/c)\beta_1)\\
 \dot{\beta_2} &= \beta_2^2 - c \beta_1 \beta_2 + (R'(\bar{U})-\lambda) \beta_1.
\end{aligned}
\end{equation}

\subsection{Asymptotic and far-field hyperbolicity}
We turn to the asymptotic systems associated with \eqref{eq:bigaeps}--\eqref{eq:littleaeps}.  Broadly speaking, these systems encode `far-field' information about the wave; for instance, they are useful in determining the boundaries of the essential spectrum (see \cite{kapprom}). The matrices $a(\xi,\lambda,\eps)$ and $A(\zeta(\xi),\lambda,\eps)$ tend to the limits $a^{\pm}(\lambda,\eps)$ and $A^{\pm}(\lambda,\eps)$ as $\xi \to \pm\infty$. Here the minus (resp. plus) superscript corresponds to an evaluation at $\bar{z}^-$ (resp. $\bar{z}^+$). We have
\begin{equation}\label{eq:asymptmat}
\begin{aligned}
a^{\pm}(\lambda,\eps) &= \begin{pmatrix} - \frac{1}{c} ( D(u^{\pm}) + \eps \lambda) & 0 &\frac{1}{c}  \\ \eps (R'(u^{\pm}) - \lambda) & 0 & 0 \\ \eps c & - \eps & 0 \end{pmatrix},\\
A^{\pm}(\lambda,\eps) &= (1/\eps) a^{\pm}(\lambda,\eps).
\end{aligned}
\end{equation}

Note from \eqref{eq:reactionterm} that $R'(1) < R'(0) < 0$. Let $\beta$ be some fixed real number with $R'(0) < \beta < 0$ and define 
\begin{align} \label{eq:omegadef}
\Omega &= \{\lambda \in \mathbb{C}: \text{Re}\,\lambda > \beta\}.
\end{align}

We will show in Sec. \ref{sec:spectrum} that the essential spectrum lies entirely in the left-half complex plane, and  bounded away from $\Omega$.  We now record some estimates for the eigenvalues and eigenvectors of the asymptotic matrices \eqref{eq:asymptmat} for any $\lambda \in \Omega$, which will be important in determining the point spectrum.

\begin{lem} \label{lem:asympt}
Let $\lambda \in \Omega$. Then for $\eps > 0$ sufficiently small, the asymptotic matrices \eqref{eq:asymptmat} have three distinct eigenvalues $\mu^{\pm}_j(\lambda,\eps)$ for $a^{\pm}(\lambda,\eps)$ (resp. $\Delta^{\pm}_j(\lambda,\eps)$ for $A^{\pm}(\lambda,\eps)$), which satisfy
\begin{equation} \label{eq:eigshierarchy}
\begin{aligned}
\text{Re}\,\mu^{\pm}_f \ll \text{Re}\,\mu^{\pm}_{s,1} < &\,0 < \text{Re}\,\mu^{\pm}_{s,2}\\
\Delta^{\pm}_j &= (1/\eps) \mu^{\pm}_j.
\end{aligned}
\end{equation}

In particular, we have $\mu^{\pm}_f = \bigO(1)$ and $\mu^{\pm}_{s,j}= \bigO(\eps)$, $ j= 1,2$. Furthermore, let $e_f^{\pm}(\lambda,\eps)$ and $e_{s,j}^{\pm}(\lambda,\eps)$ denote the eigenvectors associated with $\mu^{\pm}_f$ and $\mu^{\pm}_{s,j}$, respectively, where $j = 1,2$. Then each $e^{\pm}(\lambda,\eps)$ limits to a `reduced' eigenvector $r^{\pm}(\lambda)$ as $\eps \to 0$. These reduced eigenvectors are given explicitly by
\begin{equation}
\label{eq:redvecs}
\begin{aligned}
r_f^{\pm} &= (1,0,0)^{\top}\\
r_{s,1}^{\pm} &= (1/D(\bar{u}^{\pm} ), \nu_{p,\pm} , 1 )^{\top}\\
r_{s,2}^{\pm} &= (1/D(\bar{u}^{\pm} ), \nu_{m,\pm} , 1)^{\top},
\end{aligned}
\end{equation} 
where
\begin{equation}
\label{eq:redvecs2}
\begin{aligned}
\nu_{p,\pm} &= \frac{c_0+\sqrt{c_0^2 + 4D(\bar{u}^{\pm})(\lambda-R'(\bar{u}^{\pm}))}}{2D(\bar{u}^{\pm})}\\
\nu_{m,\pm} &= \frac{c_0-\sqrt{c_0^2 + 4D(\bar{u}^{\pm})(\lambda-R'(\bar{u}^{\pm}))}}{2D(\bar{u}^{\pm})}.
\end{aligned}
\end{equation} 
\end{lem}

\begin{remk}
The smallness of $\eps$ required for Lemma \ref{lem:asympt} is dependent on $\lambda$, but a uniform bound $\eps \leq \bar{\eps}$ over closed, bounded contours in $\mathbb{C}$ can be extracted using compactness.
\end{remk}

{\it Proof of Lemma \ref{lem:asympt}:} We outline the proof as it is similar to the proof of Lemma \ref{lem:eighierarchy}. The characteristic polynomial of $a^{\pm}(\lambda,\eps)$ is 
\begin{align}\label{eq:charfast}
p^{\pm}(\mu) &= - \mu^3 - \mu^2 \left( \frac{D(u^{\pm})+\eps \lambda}{c} \right) + \mu\eps  + \eps^2 \frac{\lambda - R'(u^{\pm})}{c},
\end{align}

from which we obtain the expansion $\mu^{\pm}_f = -D(\bar{u}^{\pm})/c_0 + \bigO(\eps)$. Leading-order expansions of the $\bigO(\eps)$ eigenvalues are obtained by using the scaling $\mu = \eps \nu$ as before. Their leading order terms are given by $\nu_{p,\pm}$ and $\nu_{m,\pm}$; this can be checked by direct calculation. The remaining small eigenvalues therefore fall into the hierarchy stated by the lemma at both $u^{-}$ and $u^{+}$ provided that $\text{Re}\,\lambda >R'(0)$. The expressions for the singular limits of the eigenvectors can also be checked with algebra. $\Box$\\

For fixed $a >0$, there exists $\bar{\eps} > 0$ sufficiently small so that the wave $x(\xi,\eps)$ lies near $S^{a,-}_0 \cup S^{a,+}_0$  for $|\xi| \geq a/\eps$ for all $0 < \eps \leq \bar{\eps}$. The corresponding coefficient matrix $a(\xi,\lambda,\eps)$ has three eigenvalues $\mu_i$ satisfying the hierarchy $\mu_f \ll \mu_{s,1} < 0 < \mu_{s,2}$, inherited from the asymptotic hyperbolicity described in Lemma \ref{lem:asympt}. Let $f_1(\xi,\lambda,\eps),\,f_2(\xi,\lambda,\eps),\,f_f(\xi,\lambda,\eps)$ denote a choice of eigenvectors corresponding to the eigenvalues $\mu_{s,1},\,\mu_{s,2},\,\mu_f$ (note: if $\mu_{s,1}$ and $\mu_{s,2}$ coalesce, we let $f_1,\,f_2$ denote the generalized eigenvectors instead). We furthermore suppose that the $f_i$ are normalized so that $|f_i|_{\infty} = 1$ for all $(\xi,\lambda,\eps)$.\\

\begin{defi} \label{def:slowsubbundle}
The {\it slow subbundle} $\sigma_s$ of the linearized system is 
\begin{align} \label{eq:slowsub}
\sigma_s(\xi,\lambda,\eps) &= \bigcup_{|\xi|\geq a/\eps}\text{span}\{f_1(\xi,\lambda,\eps),\,f_2(\xi,\lambda,\eps)\}
\end{align}
where the base space is $|\xi|\geq a/\eps$ and $\lambda \in \Omega$. 
We denote the slow subbundle with respect to the timescale $\zeta = \xi/\eps$ by $\Sigma_s$, i.e.
\begin{align} \label{eq:slowsub}
\Sigma_s(\zeta,\lambda,\eps) &= \sigma_s(\xi/\eps,\lambda,\eps).
\end{align}
\end{defi}

\subsection{Complex 2-plane dynamics induced by the eigenvalue problem}

We now write down a concrete representation of the evolution equations of complex $2$-planes that are induced by the eigenvalue problem \eqref{eq:linreversefast}, and recording some essential properties about this induced system using the formulation described in Sec. \ref{sec:geometry}. We have
\begin{align} \label{eq:2planeseig}
\begin{pmatrix} u \wedge p \\ p \wedge v\\ v \wedge u\end{pmatrix}' &= \begin{pmatrix} -\frac{1}{c} \left( \eps \lambda + D(\bar{u}) \right) & -\frac{1}{c} & 0 \\ - \eps c & 0 & -\eps (R'(\bar{u})-\lambda) \\ \eps & 0 & -\frac{1}{c} (\eps \lambda + D(\bar{u}))\end{pmatrix} \begin{pmatrix} u \wedge p \\ p \wedge v \\ v \wedge u \end{pmatrix}.
\end{align}

The coupled system \eqref{eq:travelfast}--\eqref{eq:2planeseig} acts on the complexified tangent bundle $ T^{\mathbb{C}}U \cong \cup_{p\in U} T_{p}(U) \otimes_{\mathbb{R}} \mathbb{C}$, where $U \subset \mathbb{R}^3$. In view of Corollary \eqref{cor:waveestimate} and Lemmas \ref{lem:asympt} and \ref{lem:twowedgeeigs}, there exists $\bar{\eps} > 0$ and $a > 0$ so that for each $\eps \in (0,\bar{\eps})$ and $|\xi| > a/\eps$, the eigenvalues $\mu^{(2)}_1(\xi),\,\mu^{(2)}_2(\xi), \mu^{(2)}_3(\xi)$ of the linearization of \eqref{eq:2planeseig} can be ordered so that 
\begin{align*}
\text{Re}(\mu^{(2)}_1) &= \mathcal{O}(1) < 0\\
\text{Re}(\mu^{(2)}_2) &= \mathcal{O}(1) < 0\\
\text{Re}(\mu^{(2)}_3) &= \mathcal{O}(\eps).
\end{align*}

Let us provide some geometric intuition to clarify these estimates. When the wave lies near the slow manifolds, the eigenvalue problem \eqref{eq:linreversefast} provides an eigenvector lying near the direction of strong contraction onto the slow manifolds, and a pair of  eigenvectors lying near to the tangent space of the slow manifold. These eigenvectors can be chosen pairwise to construct an eigenbasis for the linearisation of system \eqref{eq:2planeseig}, consisting of three complex $2$-planes. Two of the constructed 2-planes have an axis lying near the fast fibers of the slow manifolds, so they strongly contract in forward time under the flow; these are the eigenvectors corresponding to the two eigenvalues of $\mathcal{O}(1)$ negative real part. The remaining 2-plane lies near the tangent space of the slow manifold. Depending on the relative magnitudes of the singular eigenvalues along the wave, this eigenvector provides a weakly (un)stable direction.\\

We can further clarify the local dynamics near the tails of the wave, in view of the elementary facts about projectivized systems discussed just after \eqref{eq:projdyn}.  The projectivized system induced by \eqref{eq:2planeseig} on $\mathbb{CP}^2$, written in terms of the coordinate representation $(\beta_1,\beta_2) = \left( \frac{p\wedge v}{u\wedge p}, \frac{v\wedge u}{u\wedge p} \right)$, $u\wedge p \neq 0$, is given by
\begin{equation} \label{eq:proj2planes}
\begin{aligned} 
\beta_1 ' &= -\eps c + \frac{\left( \eps \lambda + D(\bar{u}) \right)}{c} \beta_1 - \eps (R'(\bar{u})-\lambda) \beta_2+ \frac{1}{c}\beta_1^2 \\
\beta_2' &= \eps + \frac{1}{c}\beta_1 \beta_2,
\end{aligned}
\end{equation}

where we remind the reader of the nonautonomous nature of the flow due to the inclusion of the phase space variable $\bar{u} = \bar{u}(\xi)$. associated with the system \eqref{eq:proj2planes} is a family of {\it frozen systems}, defined for each fixed $\eps$ by replacing the nonautonomous system by a one-parameter family of autonomous ODEs where the variable $\xi$ is formally replaced by an independent parameter $\gamma$ on the right hand side (see Sec. \ref{sec:elephanttrunks} for the explicit construction). A straightforward perturbation argument shows that for $\eps > 0$ sufficiently small with $a >0$ defined as above and for each fixed $\gamma$ with $|\gamma|  > a/\eps$, there is an (isolated) attracting fixed point of the frozen family for each such $\gamma$, denoted $\beta_0(\gamma,\eps)$, which is $\mathcal{O}(\eps)$-close to $(\beta_1,\beta_2) = (-D(\bar{u}(\gamma)),0)$.  For fixed $\eps > 0$, the function $\beta_0(\gamma,\eps)$ defines curves of attracting critical points on the left and right branches of the slow manifold as $\gamma$ is varied, corresponding precisely to the family of eigenvector $2$-planes near the slow manifolds whose construction was just described. This is exactly the slow subbundle in Def. \ref{def:slowsubbundle}.\\
 
These curves of attracting critical points serve as the organizing centers for so-called {\it relatively invariant sets} (see \cite{GJ}), which we construct over the wave tails in Sec. \ref{sec:elephanttrunks}.

\subsection{Reduced eigenvalue problems}

The dynamics of the linearized systems \eqref{eq:visclinslow} and \eqref{eq:visclinfast} are equivalent when $\eps > 0$, but they limit to distinguished problems as $\eps \to 0$. In analogy to the definitions of the layer and reduced problems for the wave, we now define two linearized subsystems using \eqref{eq:linreverseslow}--\eqref{eq:linreversefast}. The desingularized slow eigenvalue problem will be used in Sec. \ref{sec:sloweigs} to define a {\it slow eigenvalue}. This will require the derivation of a jump condition for two disjoint solution segments, defined separately on $S^{a,-}_0$ and $S^{a,+}_0$.  \\

\begin{defi} \label{def:linlayer}
The {\it fast eigenvalue problem}  is the linear subsystem defined by the singular limit of \eqref{eq:linreversefast}, along the fast fiber bundle:

\begin{align}
u' &= -\frac{1}{c}D(\bar{u})u \label{eq:linlayer}
\end{align}
with $p = v = 0$. 
\end{defi}

\begin{remk}
The limiting system \eqref{eq:linlayer} is the variational equation of \eqref{eq:travelfast}, and thus is degenerate as an eigenvalue problem. There is {\it always} a nontrivial, uniformly bounded solution to \eqref{eq:linlayer}. 
\end{remk}

\begin{defi} \label{def:linreduced}
The {\it slow eigenvalue problem} is the linear subsystem on $S^{a,-}_{\varepsilon} \cup S^{a,+}_{\varepsilon}$ defined by the singular limit of \eqref{eq:linreverseslow}, subject to the constraint $V = D(\bar{U})U$:

\begin{equation} \label{eq:linreduced}
\begin{aligned}
\dot{P} &= (R'(\bar{U})-\lambda)\frac{V}{D(\bar{U})}\\
\dot{V} &= \frac{cV}{D(\bar{U})} - P,
\end{aligned}
\end{equation}
or more compactly as
\begin{align}
\dot{W} &= A_0(\bar{U},\lambda, c) W. \label{eq:linlargeA}
\end{align}

It is convenient to record here the {\it projectivized slow eigenvalue problem} on the chart $S = P/V$ with $V \neq 0$:
\begin{align} \label{eq:projlinreduced}
\dot{S} &= \frac{1}{D(\bar{U})} \left( R'(\bar{U}) - \lambda - cS + D(\bar{U})S^2 \right).
\end{align}

We also record here a {\it desingularized slow eigenvalue problem}, which is defined by first appending \eqref{eq:linreduced} to \eqref{eq:reducedproblem} and then rescaling the frame variable of the resulting autonomous system by $D(\bar{U})$:
\begin{equation} \label{eq:linslow}
\begin{aligned}
\dot{P} &= (R'(\bar{U})-\lambda)V\\
\dot{V} &= cV - D(\bar{U})P.
\end{aligned}
\end{equation}
\end{defi}
Projectivizations of the desingularised problem \eqref{eq:linslow} will be written down and used in Sec. \ref{sec:slowevans} to identify the slow eigenvalues.

\section{Essential \& Absolute Spectrum} \label{sec:spectrum}

Our goal is to determine {\it spectral stability} of each member of the one-parameter family of waves $\{\Gamma_{\eps}\}_{\eps \in (0,\bar{\eps}]}$ for some sufficiently small value of $\bar{\eps}>0$. Following the standard approach, we will compute the essential and point spectra, defined as in \cite{kapprom}, and compute them separately. It is more convenient to write the eigenvalue problem \eqref{eq:fulleigprob} as a first order system. Following \cite{bjornstability}, we consider the (equivalent) family $\cL$ of linear operators from \eqref{eq:linreverseslow} and \eqref{eq:linreversefast} 
\begin{align}\label{eq:toperator}
\cL(\lambda) & := \frac{d}{d\zeta} - A(\bar{U},\lambda,\eps) =  \frac{d}{d\zeta} -
 \begin{pmatrix}
- \frac{1}{c\ve}(D(\bar{U}) + \ve \lambda) & 0 & \frac{1}{c \ve} \\ 
(R'(\bar{U}) - \lambda) & 0 & 0 \\ 
c& -1 & 0 
\end{pmatrix}
\end{align}
where $\cL(\lambda): H^{1}(\R,\C^{3}) \to L^{2}(\R,\C^{3})$. \\

The {\em essential spectrum} $\sigma_{e}(\cL)$ of the family of linear operators $\cL(\lambda)$ defined in \eqref{eq:toperator} is the set of $\lambda \in \C$ such that either $\cL(\lambda)$ is not Fredholm, or $ \cL(\lambda)$ is Fredholm but the (Fredholm) index is not 0. The {\em point spectrum} $\sigma_{p}(\cL)$ is the set of values $\lambda \in \C$ where $\cL(\lambda)$ is not invertible, but does have index 0. We will denote by $\sigma(\Gamma)$ or $\sigma(\cL)$, the union $\sigma_{p} \cup \sigma_{e}$ of the spectrum of the family $\cL$ associated with the travelling wave $\Gamma$. 

\begin{defi} (Def. 4.1.7 in \cite{kapprom}) \label{def:specstab}
A wave $\Gamma$ corresponding to a stationary solution of equation \eqref{eq:travelingwavemaster} is called {\it spectrally stable} if $\sigma(\Gamma) \cap \{\lambda \in \mathbb{C}: \text{Re}(\lambda) \geq 0\} = \emptyset$, except possibly at $\lambda = 0$.   
\end{defi}

\begin{figure}
\begin{center}
\includegraphics[scale=0.22]{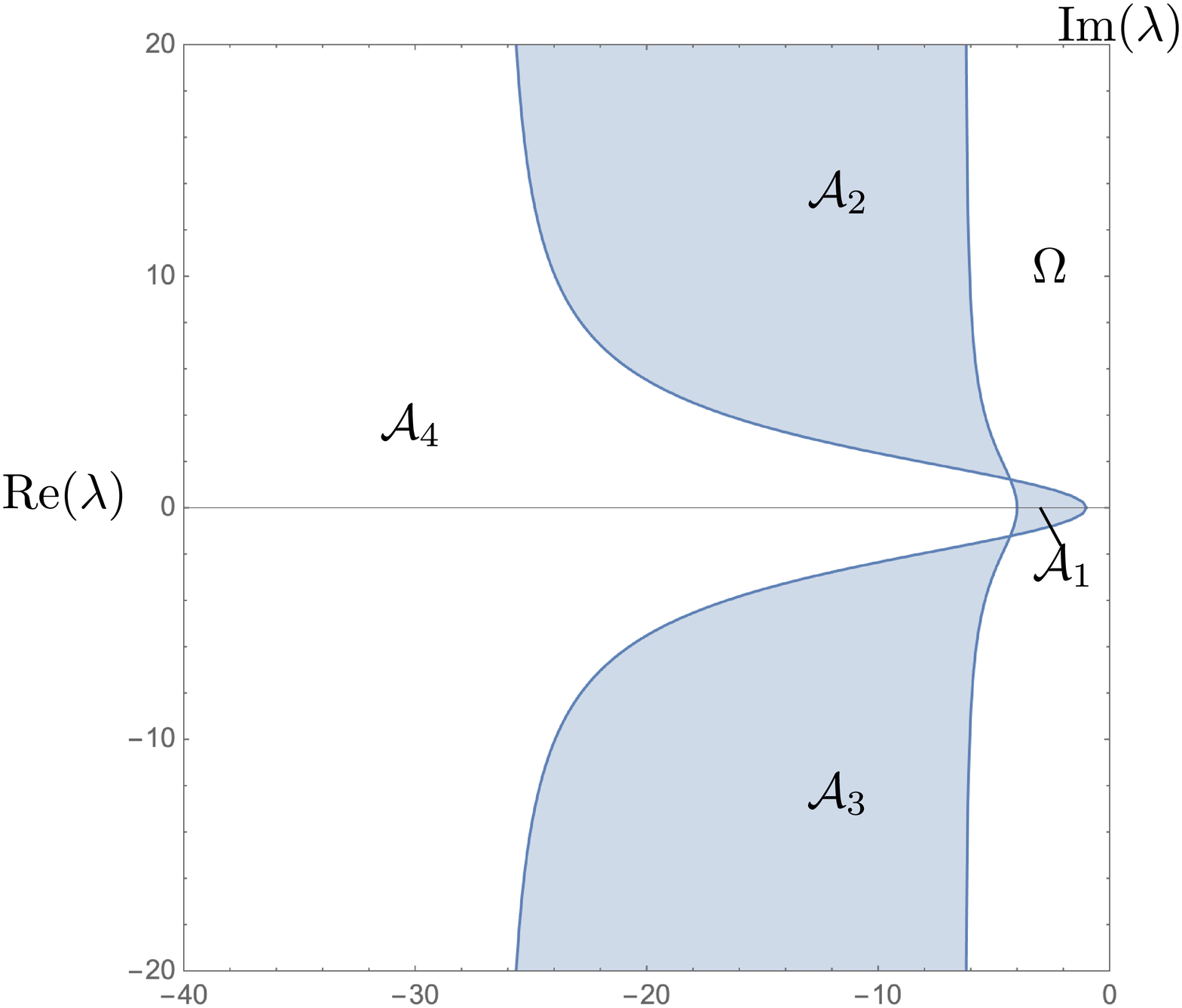}
\includegraphics[scale=0.22]{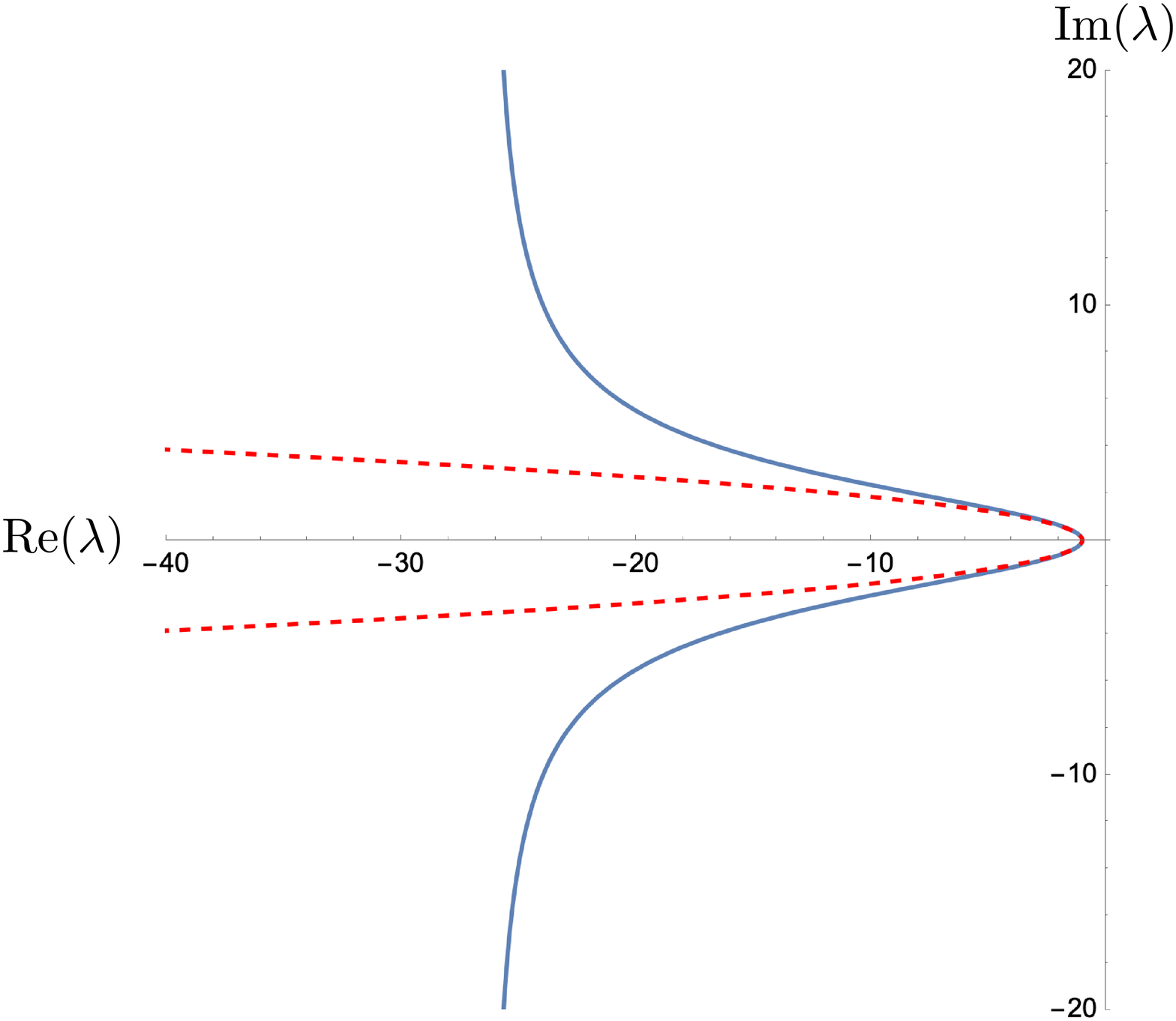}
\caption{Left: a plot of the essential spectrum (shaded regions) of the family $\cL$ illustrating its qualitative features, and noting the regions $\Omega$ and the $\cA_{i}$ partitioning the complex plane via the Fredholm borders (continuous spectrum). The asymptotics of the essential spectrum for large spatial modes means that the family of operators $\cL$, for non-zero $\ve$ is not sectorial. Right: a plot comparing the Fredholm border of $A^{-}(\lambda, \ve)$ (solid line - blue online) to the Fredholm border of the reduced problem linearised about the constant solution $\bar{U}^{-} =0$ with $\ve =0$ (dashed line - red online). As $\ve \to 0$, the solid line follows the dashed line for a longer region in the left half plane, before eventually `flaring' up. For the figures, $\ve =0.1$ (where it is present) and $c =1$, with $D$ and $R$ as in \eqref{eq:diffusionterm} and \eqref{eq:reactionterm}.} \label{fig:essential}
\end{center}
\end{figure}

For nonzero $\ve$, because $\bar{U}$ approaches its end states exponentially in $\zeta$ (equivalently $\xi$), determining the essential spectrum amounts to determining when the matrices $A^{\pm}(\lambda,\ve)$ from \eqref{eq:asymptmat} have different signatures (this is the content of Weyl's essential spectrum theorem, Theorem 2.2.6 from \cite{kapprom} as well as and Lemma 3.1.10, also from \cite{kapprom}). The edges of the essential spectrum are called {\em dispersion relations} or the {\em Fredholm borders} as well as the {\em continuous spectrum}. To find the dispersion relations we look for values where $A^{\pm}$ has a purely imaginary (spatial) eigenvalue $i k$, and solve for the temporal eigenvalue parameter $\lambda$. This gives a pair of curves $\lambda_{\pm}(k;\ve)$ in the complex plane, parametrised by the spatial modes, where the asymptotic matrices fail to be hyperbolic. Working with \eqref{eq:linreverseslow}, and denoting $\lim_{\zeta \to \pm \infty} \bar{U}$ by $\bar{U}^{\pm}$, we have that the dispersion relations are:
\begin{align}
\lambda_{\pm}(k;\ve) & = - i c k  + \frac{R'(\bupm)- D(\bupm)k^{2}}{1+k^{2}\ve}
\end{align}
We have a pair of curves partitioning $\C$ into five disjoint sections: the set to the right of both curves (where we're going to look for point spectrum) which we denote $\Omega$, together with four more regions which we denote by $\cA_{j}$ for $j=1,2,3,4$. (See Fig. \ref{fig:essential}). 

Since the characteristic polynomials of $A^{\pm}$ are both cubics: 
\begin{align}\label{eq:charslow}
P_{A^{\pm}}(\mu) = - \mu^{3} - \frac{1}{c}\left( \lambda + \frac{D(\bupm)}{\ve} \right) \mu^{2} + \frac{\mu}{\ve} + \frac{1}{c\ve}\left( \lambda - R'(\bupm)\right),
\end{align}
everything can be checked explicitly. We are working in the slow variables, but everything that follows in this section can be computed for the (equivalent) fast systems, since
$$
P_{A^{\pm}}(\mu) = \frac{1}{\ve^{3}}p^{\pm}(\ve \mu)
$$
from \eqref{eq:charfast}.
Defining the transformation
$$
\eta := \mu - \frac{1}{3}\left( \lambda + \frac{D(\bupm)}{c \ve} \right),  
$$
we have 
$$
P_{A^{\pm}}(\eta) = -\eta^{3} + H_{\pm} \eta + K_{\pm}, 
$$
where 
$$
H_{\pm}:= \frac{\lambda^{2}}{3c^{2}} + \lambda \frac{2 D(\bupm)}{3 c^{2}\ve} + \frac{D(\bupm)^{2}+3c^{2}\ve}{3c^{2} \ve^{2}},
$$
and 
\begin{align*}
K_{\pm} & := - \frac{2}{27 c^{3}} \lambda^{3}- \frac{2 D(\bupm)}{9 c^{3}\ve}\lambda^{2} + \frac{18c^{2}\ve^{2}- 6 \ve D(\bupm)^{2}}{27c^{3}\ve^{3}}\lambda \\ 
& \qquad -\frac{9c^{2}\ve D(\bupm) + D(\bupm)^{3}+ 27 c^{2}\ve^{2}R'(\bupm)}{27c^{3}\ve^{3}}.
\end{align*}

We note that the discriminants $\Delta_{H_{\pm}}$ of $H_{\pm}$ are polynomials in $\lambda$ that are independent of $\bupm$:
$$
\Delta_{H_{\pm}} = -\frac{4}{3c^{2}\ve} < 0
$$
which means that for $\lambda \in \R$, we have $H_{\pm}>0$. Hence for each $\lambda \in \R$, 
$$
\frac{H_{\pm}^{3}}{27}+ \frac{K_{\pm}^{2}}{4}  \geq 0.
$$
We thus infer that for $\lambda \in \R$, $P_{A^{\pm}}(\mu)$ has three real roots (counting multiplicity in edge cases). We can then use Descartes' rule of signs to see that for $\lambda = 0$ (and subsequently for all $\lambda \in \Omega$), we have that the signature of both $A^{+}$ and $A^{-}$ is $(-,-,+)$. That is $A^{\pm}$ has two eigenvalues with negative real part, and one with positive real part. Checking the signatures in the rest of the regions $\cA_{j}$ similarly produces the following table: 
\begin{center}
\begin{tabular}{c|| c | c }
Region & $\sgn(A^{-})$ & $\sgn(A^{+})$ \\ 
\hline
$\Omega$ & $(-, -, +)$ & $(-,-,+)$ \\ 
$\cA_{1}$ & $(-, +, +)$ & $(-,-,+)$ \\
$\cA_{2}$ & $(-, -, +)$ & $(-,+,+)$ \\
$\cA_{3}$ & $(-, -, +)$ & $(-,+,+)$ \\
$\cA_{4}$ & $(-, +, +)$ & $(-,+,+)$ \\
\end{tabular}\label{tab:sigs}
\end{center}
Thus we find that the essential spectrum is the union of regions $\cA_{2,3,4}$ together with their boundaries. This set in the complex plane with a bounded real part and an unbounded imaginary component (see Fig. \ref{fig:essential}). \\

The dispersion relations are symmetric about $k$, and they cross the real axis at $R'(\bar{U}^{\pm})$ when $k =0$. The dispersion relations are vertical lines when $\ve = \ve^{\pm}_{*} := -\frac{D(\bupm)}{R'(\bupm)}$, and are parabolic in $k$ near $k = 0$, opening to the left when $\ve < \ve_{*}^{\pm}$ and to the right when $\ve> \ve_{*}^{\pm}$ before `flaring' up to asymptote towards vertical lines. \\

It is worthwhile to investigate whether the family $\cL$ is {\it sectorial} when $\eps > 0$, since the property of sectoriality allows us to strengthen spectral stability to linear stability. (see eg. Henry \cite{henry}). When $\ve =0$ the dispersion relations define a pair of parabolas opening leftward, and as such the linearized operator associated with the unperturbed problem is sectorial.  On the other hand, for any fixed $\eps > 0$, the following asymptotics characterize the behavior of the Fredholm bolders as $k \to \pm \infty$:
\begin{align}
 \lambda_{\pm} \sim -\frac{D(\bupm)}{\ve} + -i ck. \quad 
\end{align}

Hence, the essential spectrum is contained in the left half plane, but the family $\cL$ is {\em not} sectorial since the borders are vertically asymptotic; see Fig. \ref{fig:essential}.\\

We note that as $\ve \to 0$, the continuous spectrum (Fredholm borders) of the perturbed problem will converge to the continuous spectrum (Fredholm borders) of the linearised unperturbed problem, linearised about the steady-state solutions $\bupm$. That is we set $\ve =0$ in \eqref{eq:travelingwavemaster}, and consider the dispersion relations of the steady state solutions $\bupm$. We see that these are given by
$$\lambda_{\pm}(k;0) = - {k^{2}}{D(\bupm)} + R'(\bupm) +ick, $$
which corresponds to setting $\ve =0$.  \\

 We highlight the subtlety underlying the asymptotic convergence of the essential spectrum of the linearized operator as $\eps \to 0$: for any large bounded set intersecting the essential spectrum of the operator for $\eps > 0$ and for any fixed $\delta > 0$, there exists a sufficiently small uniform bound $\bar{\eps} > 0$ so that for each $\eps \in (0,\bar{\eps})$, the Fredholm borders for the operator when $\eps > 0$ can be made $\delta$-close to the (parabolic) borders of the reduced operator within this bounded set; on the other hand, the uniform bound $\bar{\eps}$ must be made smaller and smaller as we take larger and larger bounded sets and demand $\delta$-closeness. See Fig. \ref{fig:essential}.

\begin{remk}
In view of the nonsectoriality of the family $\cL$, spectral stability does not immediately imply {\it linear exponential stability}, i.e. it is unknown whether nearby translates of the spectrally stable wave decay exponentially onto the wave in forward time. 
\end{remk}

\begin{remk}
Lastly, we note that as the essential spectrum is entirely contained in the left half plane, any {\em absolute spectrum} will also be contained in the left half plane, and so will play no role in destabilising the underlying wave. 
\end{remk}

\section{Augmented unstable bundles} \label{sec:bundles}

Having shown that the essential spectrum lies to the left of $\Omega$ in the left-half complex plane, we now focus on finding the point spectrum $\sigma_p(\cL)$ consisting of isolated eigenvalues of finite multiplicity. To that end, we fix a a simple, closed curve $K \subset \Omega$, such that $\sigma_{p} \cap K$ is empty. In view of the hyperbolicity of the asymptotic systems $a^{\pm}(\lambda,\eps)$, we henceforth refer to the values in the point spectrum as {\it eigenvalues}. Specifically, we rely on the following geometric characterization of the eigenvalues of \eqref{eq:fulleigprob}.

\begin{defi} \label{def:eigenvalues}
The number $\lambda \in \Omega$ is an {\it eigenvalue} of $\cL(\lambda)$ if there exists a nontrivial uniformly bounded solution $y(\xi)$ of \eqref{eq:linreverseslow} (equivalently, if the bounded solution $y(\xi)$ decays to zero as $\xi \to \pm \infty$; see \cite{AGJ}).
\end{defi}

In this section we construct two {\it augmented unstable bundles} $\mathcal{E}_{\eps}(K)$ and $\mathcal{E}_0(K)$, over a cylindrical base space formed by extending $K$ along a compactified time variable.  As we will outline in Sec. \ref{sec:outline}, our goal is to use a particular correspondence theorem that relates a certain topological invariant of these complex vector bundles to the eigenvalue count of $\cL(\lambda)$ inside $K$. 

\subsection{Construction of the augmented unstable bundle $\mathcal{E}_\eps(K)$ for $\eps > 0$}  \label{sec:augeps}
Following \cite{GJ}, we now use general facts in \cite{AGJ} together with Lemma \ref{lem:asympt} to construct (un)stable complex vector bundles of solutions along the wave for each sufficiently small $\eps > 0$.

\begin{defi} \label{def:bundles}
Let $\lambda \in \Omega$. For $\eps > 0$ sufficiently small, the {\it unstable bundle} and {\it stable bundle} $\varphi^{-}(\xi,\lambda,\eps)$ resp. $\varphi^{+}(\xi,\lambda,\eps)$  are the linear subspaces of solutions along the wave $\Gamma_{\eps}$ defined by the conditions
\begin{equation} \label{eq:bundles}
\begin{aligned}
\varphi^{-}(\xi,\lambda,\eps) &\to \text{span}\{e_{s,2}^{-}(\lambda,\eps)\} &&\text{as~~} \xi\to -\infty \hspace{1cm}\\
\varphi^{+}(\xi,\lambda,\eps) &\to \text{span}\{e_{s,1}^{+}(\lambda,\eps), e_f^+(\lambda,\eps)\} &&\text{as~~} \xi\to +\infty.\hspace{1cm}
\end{aligned}
\end{equation}

The asymptotic convergence to the corresponding (un)stable subspaces of $a^{\pm}(\lambda,\eps)$ is defined with respect to the topologies of the appropriate Grassmannians $G_{1,3}$ and $G_{2,3}$.
\end{defi}

The linear subspace $\varphi^-$ (resp. $\varphi^+$) can be viewed as a complex line bundle (resp. 2-plane bundle) over the base space $(\xi,\lambda) \in \mathbb{R}\times \Omega$. In view of Definitions \ref{def:eigenvalues} and \ref{def:bundles}, $\lambda \in \Omega$ is an eigenvalue if and only if $\varphi^{-}(\xi,\lambda,\eps)$ and $\varphi^+(\xi,\lambda,\eps)$ intersect nontrivially for some (and hence all) $\xi \in \mathbb{R}$.\\
 
 We now describe a general construction for obtaining an {\it augmented unstable bundle} $\mathcal{E}_{\eps}(K)$ over the real 2-sphere $S^2$, using $\varphi^{-}$; we refer to \cite{AGJ} for details. The approach is to first compactify $\xi$ using the change of variable $\xi \to \tau$, where
 \begin{align*}
\tau' &= \eps \kappa (1-\tau^2),\\
\tau(0)  &= 0,
\end{align*}
 
for $\kappa > 0$ chosen such that $\tau(\xi)\to \pm 1$ as $\xi \to \pm \infty$. With respect to this compactified `timescale', $\varphi^{-}(\tau,\lambda,\eps)$ then specifies a line bundle over $(-1,1) \times \Omega$; a bundle can then be defined over the restricted base space $ (-1,1) \times K$. This base space is homeomorphic to an infinite cylinder (equivalently, a finite cylinder with its caps removed). Our first goal is to {\it continuously extend} the line bundle over the closure of the cylinder (i.e. over the base space $[-1,1] \times K$). After this we fill in the caps, altogether producing a complex vector bundle over the closed cylinder, i.e. over $S^2$. \\

Let us explain this continuous extension. First recall that every complex vector bundle can be realized as a pullback of the universal bundle over a Grassmannian, as follows. Given a (complex) vector bundle $\mathcal{E}: E\subset B \times \mathbb{C}^n \to B$ with projection map $\pi$, there is a natural map $\hat{e}:B \to G_{k,n}$, which assigns to each $b\in B$ the complex $k$-dimensional linear fiber $\pi^{-1}(b)$, considered as a $k$-plane in  $G_{k,n}$. On the other hand, the universal bundle $\Gamma_k(\mathbb{C}^n) \to G_{k,n}$ is a trivial (tautological) construction: the fiber above each element in $G_{k,n}$ is defined to be the corresponding $k$-plane. Then $\mathcal{E}$ can be realized as $\hat{e}^*\Gamma_k(\mathbb{C}^n)$:
$$
\begin{tikzcd}
E \arrow[r,"  "'] \arrow[d,"\pi "']  &  \Gamma_k(\mathbb{C}^n) \arrow[d,"  "] \\
   B   \arrow[r," \hat{e} "']  & G_{n,k}
 \end{tikzcd}
$$

Thus, we say that $\varphi^{-}|_{(-1,1)\times K}$ continuously extends to a bundle over $[-1,1]\times K$ if the natural map $\hat{e}:(-1,1) \times K \to G_{1,3}$ extends continuously to $\tilde{e}:[-1,1] \times K \to G_{1,3}$. The extension to $\tau = -1$ is clear: for each $\lambda \in K$, $\varphi^-$ continuously extends to $e^-_{s,2}(\lambda,\eps)$ as $\tau \to -1$ by definition. At the other end, we need to specify that $\lambda \in K$ is {\it not} an eigenvalue of \eqref{eq:linreverseslow}. General considerations about $\omega$-limit sets relative to hyperbolic fixed points (see \cite{AGJ}) then imply that $\varphi^-(\tau,\lambda,\eps)$ continuously extends to the corresponding unstable eigenvector $e^{+}_{s,2}(\lambda,\eps)$ as $\tau \to +1$. We note that this limit can be made uniform in $\lambda$ over the closed contour $K$.\\

It remains to fill in the interiors of the caps $C_- = \{\tau = -1\} \times K^0$ and $C_+ = \{\tau = +1\} \times K^0$ at either end of the cylinder. For each $\lambda \in K^0$, this is achieved for both caps by gluing the corresponding unstable eigenvector $e^{\pm}_{s,2}(\lambda,\eps)$. 

\begin{defi}
The {\it augmented unstable bundle} $\mathcal{E}_\eps(K)$ is the complex line bundle over the cylindrical base space $C_- \cup [-1,1] \times K \cup C_+$, constructed from $\varphi^{-}$ as described above.
\end{defi} 
This base space is homeomorphic to $S^2$. Our goal will be to relate the eigenvalue count to a particular characteristic class defined for complex vector bundles.

\subsection{An outline of the proof of the stability theorem} \label{sec:outline}
We make use of the main theorem proved in \cite{AGJ}: the total number of eigenvalues (counted with algebraic multiplicity) inside $K$ is equal to the first Chern number, $c_1(\mathcal{E}_{\eps}(K))$, of the complex vector bundle constructed above. The main advantage of this topological characterisation lies in the homotopy invariance of $c_1$.\\

 We can now outline the steps of the proof of Theorem \ref{thm:main}:
\begin{enumerate}
\item Construct a {\it reduced} augmented unstable vector bundle $\mathcal{E}_{0}(K)$ using a complex two-dimensional reduced (or {\it slow}) eigenvalue problem defined on the critical manifold $S_0$. This is accomplished in Sec. \ref{sec:aug0}.
\item Show that for sufficiently small $\eps > 0$, we can continue $\mathcal{E}_{\eps}(K)$ to $\mathcal{E}_{0}(K)$; i.e. we can construct a homotopy between them. By homotopy invariance, $c_1(\mathcal{E}_{\eps}(K))= c_1(\mathcal{E}_{0}(K))$. This homotopy is constructed at the end of Sec. \ref{sec:convergence} only after several key estimates are proven; see Theorem \ref{thm:approx}.
\item Use the correspondence theorem in $\cite{AGJ}$ in the `other direction,' i.e. compute $c_1(\mathcal{E}_{0}(K))$ by finding the slow eigenvalues using Sturm-Liouville techniques.  We will verify that the reduced problem admits only a simple translational eigenvalue $\lambda_0 = 0$, and another simple eigenvalue $\lambda_1$ with $R'(0) < \text{Re}\,\lambda_1 < 0$. The slow eigenvalues are computed in Sec. \ref{sec:slowevans}.
\end{enumerate}

Steps 1 and 3 turn out to be relatively straightforward; most of the technical issues lie in proving the uniform estimates required to construct the homotopy in step 2. To summarize, we will show that for any fixed contour $K \subset \mathbb{C}$ to the right of the essential spectrum, there exists $\bar{\eps} > 0$ so that for each $\eps \in (0,\bar{\eps}]$, the following statements hold:
\begin{align*}
\{\# \text{ of eigenvalues of }  \cL \text{ inside } K\} &= c_1(\mathcal{E}_{\eps}(K)) = c_1(\mathcal{E}_{0}(K)) \leq 2,
\end{align*}

where the first equality follows from the general theory in \cite{AGJ};  the second follows from the Step 2 above; and the third follows from Steps 1 \& 3. By choosing smaller contours $K_0$ and $K_1$ surrounding only $\lambda_0$ resp. $\lambda_1$ (and shrinking $\bar{\eps}$ finitely many times if necessary so that the corresponding homotopies can be constructed for all three contours $K,\,K_0$, and $K_1$ simultaneously), we obtain linearized stability of the one-parameter family of waves $\Gamma_{\eps}$ for $0 < \eps < \bar{\eps}$.

\subsection{Slow eigenvalues and the jump map} \label{sec:sloweigs}

We now proceed with step 1 of the stability proof. In order to construct the reduced bundle, we must first define the appropriate geometric criteria for a {\it slow eigenvalue} $\lambda \in \mathbb{C}$ for \eqref{eq:visclinslow}. Let $(\bar{U}(\zeta),\bar{P}(\zeta))$ parametrize the segments of the singular heteroclinic orbit  $\Gamma_0$ from Hypothesis \ref{hyp:singhet} on $S^{a,-}_0 \cup S^{a,+}_0$, so that the fast jump occurs at $\zeta = 0$.

\begin{defi} \label{def:sloweigs}
 The complex number $\lambda \in \Omega$ is called a {\it slow eigenvalue} of \eqref{eq:visclinslow} if there exists a pair of nontrivial, uniformly bounded solutions on $S^{a,-}_0 \cup F_-$ resp. $S^{a,+}_0$, denoted 
\begin{equation}
\label{eq:slowsols}
\begin{aligned}
Y_L(\lambda,\zeta) &\text{ for }  \zeta \leq 0 \text{ and}\\
Y_R(\lambda,\zeta) &\text{ for } \zeta \geq 0\\
\end{aligned}
\end{equation} 

of the desingularized slow eigenvalue problem \eqref{eq:linslow}, such that
\begin{align} \label{eq:jump}
Y_R(\lambda,0) &= J_{\lambda}\left(Y_L(\lambda,0)\right),
\end{align}

where $J_{\lambda}:\mathbb{C}^2 \to \mathbb{C}^2$ is the (linear) {\it jump map}  defined by
\begin{align} \label{eq:jumpmatrix}
J_{\lambda}(P,V) &= \begin{pmatrix} 1 & \frac{R(\bar{U}_J) - R(\bar{U}_F) - \lambda(\bar{U}_J-\bar{U}_F)}{c\bar{U}_F - \bar{P}_F} \\ 0 & \frac{c\bar{U}_J - \bar{P}_F}{c\bar{U}_F - \bar{P}_F} \end{pmatrix} \begin{pmatrix} P \\ V  \end{pmatrix}.
\end{align}
\end{defi}

The jump map will be derived rigorously as a singular limit in Sec. \ref{sec:jumpmapconv}.

\begin{remk} \label{rem:slow0}
 Suppose that $\lambda$ is a slow eigenvalue. Owing to the asymptotic hyperbolicity of \eqref{eq:linslow} as $\zeta \to \pm \infty$, the pair of nontrivial, uniformly bounded solutions \eqref{eq:slowsols} must decay to zero in forward (resp. reverse) time. Given two nontrivial, uniformly bounded `half-solutions' $Y_L(\lambda,\zeta)$ and $Y_R(\lambda,\zeta)$ as given in \eqref{eq:slowsols}, any nontrivial complex scalar multiple of these solutions $\tilde{Y}_L(\lambda,\zeta) = K Y_L(\lambda,\zeta)$ and $\tilde{Y}_R(\lambda,\zeta) = K Y_R(\lambda,\zeta)$ also satisfy the jump condition \eqref{eq:jump} owing to the linearity of the jump map \eqref{eq:jumpmatrix}. \\

It is instructive to consider the case $\lambda = 0$, corresponding to a reduction of the translational eigenvalue. The reduced vector field provides (complex scalar multiples of) nontrivial, uniformly bounded linear solutions on $S^{a,-}_0\cup F_-$ and $S^{a,+}_0$ as required in \eqref{eq:slowsols}. Indeed, any such solution pair $Y_L,\,Y_R$, separately defined along the segments of the wave lying in $S^{a,-}_0\cup F_-$ and $S^{a,+}_0$, must satisfy
\begin{align*}
Y_L(\lambda,0) &= (P_F,V_F) = K_1 (R(\bar{U}_F),c\bar{U}_F-\bar{P}_F)\\
Y_R(\lambda,0) &=(P_J,V_J) = \,K_2 (R(\bar{U}_J),c\bar{U}_J-\bar{P}_F),
\end{align*}
 where $K_1,\,K_2 \neq 0$ denote complex scalars. On the other hand, 
\begin{align*}
J_0(P_F,V_F) &= K _1   \begin{pmatrix}   R(\bar{U}_F) + (R(\bar{U}_J) - R(\bar{U}_F)) \\
  (c\bar{U}_J - \bar{P}_F)   \end{pmatrix}\\
  &= K_1 \begin{pmatrix} R(\bar{U}_J) \\ c\bar{U}_J - \bar{P}_F  \end{pmatrix}\\
  &= (K_1/K_2) (P_{J},V_{J}).
\end{align*}

The jump condition \ref{eq:jump} is therefore satisfied by rescaling one of the solutions if necessary, and hence $0$ is a slow eigenvalue. Observe that the reduced vector field provides a natural scaling $K_1 = K_2 = 1$ such that the jump condition is automatically satisfied, but this is a special constraint that arises from the relationship between the reduced vector field along the singular limit of the wave and  the variational equations, i.e. the linearised problem when $\lambda = 0$. We do not generally expect (or need) to find such a natural scaling when $\lambda \neq 0$.\\

We end this remark by comparing our definition of a slow eigenvalue to that in \cite{GJ}. Motivated by the behavior of their slow line subbundle in the fast inner layer as $\eps \to 0$, their slow eigenvalue problem is defined by asking for {\it continuous} solutions across the jump. Hence, 0 is necessarily {\it not} a slow eigenvalue in their construction, since the tangent vector of a continuous solution at the end of the jump is transverse to the reduced vector field at the jump-on point on $S^{a,+}_0$. The translational eigenvalue is instead counted by their fast eigenvalue problem. 
\end{remk}

\subsection{The construction of $\mathcal{E}_0(K)$} \label{sec:aug0}

We work with the projectivization of \eqref{eq:linslow} on $\mathbb{CP}$:
\begin{align} \label{eq:projlinslow}
\dot{\hat{Z}} &= \hat{W}_0(\hat{Z},\zeta,\lambda).
\end{align}

Let $\hat{Z}_L(\zeta,\lambda) $ denote the (unique) solution of \eqref{eq:projlinslow} on $S^{a,-}_0$ which tends to the image of the unstable eigenvector of the (projective) asymptotic reduced system $\hat{W}^{-}_{0}(\hat{Z},\lambda)$ as $\zeta \to -\infty$. We write $Z_L(\zeta,\lambda) = (P_L(\zeta,\lambda),V_L(\zeta,\lambda)) \in \pi^{-1}(\hat{Z}_L)$. Let us denote by $Z_R(\zeta,\lambda) = (P_R(\zeta,\lambda),V_R(\zeta,\lambda))$ the corresponding solution on $S^{a,+}_0$ for $\zeta \geq 0$ which satisfies the jump condition \eqref{eq:jump} (note that this gives a concrete initial condition for the solution on $S^{a,+}_{0}$).  We find the corresponding embedded solutions in $\mathbb{R}^3$ using 
\begin{align*}
E_L(\zeta,\lambda) &= \iota_0(\zeta,\lambda) Z_{L}(\zeta,\lambda)\\
E_R(\zeta,\lambda) &= \iota_0(\zeta,\lambda) Z_{R}(\zeta,\lambda),
\end{align*}
 where $\iota_0$ denotes the inclusion map into $TS_0$.\\
 
Finally, we consider the projectivizations $\hat{E}_{L,R}(\zeta,\lambda) \in \mathbb{CP}^2$. In analogy to the construction of $\mathcal{E}_{\eps}(K)$, we compactify the time variable and work instead with $T(\zeta)$, where
\begin{align} \label{eq:timescalingreduced}
T' &= \kappa(1-\tau^2)\\
T(0) &= 0,
\end{align}
with $-1 \leq T \leq 1$. When there are no slow eigenvalues within $\lambda\in K$, we can then continuously extend each of $\hat{E}_{L,R}$ to the appropriate cap of the compactified cylinder by taking
\begin{align*}
\hat{E}_L(-1,\lambda) &= \hat{e}_{2,0}^{-}(\lambda)\\
\hat{E}_R(+1,\lambda) &= \hat{e}_{2,0}^{+}(\lambda)
\end{align*}
for each $\lambda \in K\cup K^0$. \\

We construct $\mathcal{E}_0(K)$ by gluing together vector bundles over two hemispheres, defined separately over the singular wave on $S^{a,\pm}_0$. In particular, let
\begin{align}
i^{\pm}_0(\lambda) &= \lim_{\zeta \to 0^{\pm} }i_0 (\zeta,\lambda)
\end{align}
and take
\begin{align*}
B_- &= B \cap \{\tau \leq 0\}\\
B_+ &= B \cap \{\tau \geq 0\},
\end{align*}
denoting the two sets splitting the base space over left and right hemispheres. By pulling back $\hat{E}_L$ and $\hat{E}_R$ with respect to these restricted base sets, we can obtain two bundles $\mathcal{E}^{\pm}_R$. It remains to glue these bundles together over the hemisphere $B_- \cap B_+$. 

The fibers over $B_- \cap B_+$ of  each bundle $\mathcal{E}^{\pm}_0$ is specified by
\begin{align*}
\mathcal{E}^{-}_0|(B_- \cap B_+) &=\text{span} \{ i_0^-(\lambda)Z_{L,0}(\lambda)\}\\
\mathcal{E}^{+}_0|(B_- \cap B_+) &=\text{span} \{ i_0^+(\lambda)Z_{R,0}\} =  \text{span}\{i_0^+(\lambda) J_\lambda(Z_{L,0})(\lambda)\},
\end{align*}

 where $Z_{L,0} \in \pi^{-1}\hat{Z}_L(0,\lambda)$ and $Z_{R,0} \in \pi^{-1}\hat{Z}_R(0,\lambda)$. Then there is a bundle isomorphism
 \begin{align} \label{eq:clutch}
\varphi^H:\mathcal{E}^{-}_0|\{0\} \times K \to \mathcal{E}^{+}_0|\{0\} \times K
\end{align}
defined by (the linear extension of)
\begin{align}
\varphi^H(\lambda)i_0^{-}(\lambda)Z_{L,0} &= i_R^{+}(\lambda)Z_{R,0} =  i_R^{+}(\lambda) J_{\lambda}Z_{L,0}.
\end{align}
The basic properties of bundle isomorphism can be checked by using the linearity of the embeddings and the jump map.\\

\begin{defi}
The {\it reduced augmented unstable bundle} $\mathcal{E}_0(K)$ is defined by the clutching operation
\begin{align}
\mathcal{E}_0 &= \mathcal{E}^-_0 \cup_{\varphi^H} \mathcal{E}^+_0,
\end{align}
i.e. it is the bundle over the closed cylinder ($\cong S^2$) obtained by gluing together the hemispheric bundles along their common boundary via the equivalence relation
\begin{align*}
y_- &\sim \varphi^H(\lambda)y_-.
\end{align*}
using the clutching map \eqref{eq:clutch}.
\end{defi}

\begin{remk} \label{rem:projhomotopy}
We refer the reader to the standard treatments in differential topology for a general overview of clutching operations (see e.g. \cite{atiyah,hatcher}). We will find it useful to work instead with an equivalent bundle $\tilde{\mathcal{E}}_0(K) \cong \mathcal{E}_0(K)$ obtained by projecting the fiber at each point over the base set onto the corresponding two components. It is straightforward to show that these two bundles are indeed homotopy equivalent by constructing the explicit homotopy that continuously bends each fiber as a homotopy parameter goes from $0$ to $1$; see the remark below Lemma 6.2 of \cite{GJ} for details.
\end{remk}

\section{Comparison of the bundles $\mathcal{E}_{\eps}(K)$ and $\mathcal{E}_{0}(K)$} \label{sec:convergence}
With the reduced objects from Sec. \ref{sec:aug0} in hand, let $\hat{E}(\zeta,\lambda,\eps)$ denote the (unique) projective solution of the $\eps$-dependent family of solutions tending to the unstable subspace as $\zeta \to -\infty$. By general arguments in \cite{AGJ}, we have that $\lambda \notin \sigma_p(\cL)$ if and only if $\hat{E}(\zeta,\lambda,\eps) \to \hat{e}_f^+(\lambda,\eps)$ (defined in Lemma \ref{lem:asympt}) as $\zeta \to +\infty$. We want to show that particular solutions of the projectivization of the slow eigenvalue problem \eqref{eq:linslow} uniformly approximate $\hat{E}(\zeta,\lambda,\eps)$ when $\eps > 0$ is sufficiently small .\\

Let $\Pi: \mathbb{C}^3 \to \mathbb{C}^2$ denote projection onto the slow coordinates, i.e.
\begin{align}
\Pi(u,p,v) = (p,v).
\end{align}

We can also define the inclusion maps $\iota,\iota_R$ in the obvious way. If $Y,Y_R$ are solutions to  the full system (in the slow timescale, with $\eps \neq 0$) and the linearized slow system (i.e. with $\eps = 0$), respectively,  then we use the variables $Z, Z_0 \in \mathbb{C}^2$ to denote the projections:
\begin{align}
Z &= \Pi(Y)\\
Z_0 &= \Pi(Y_0).
\end{align}

Owing to the linearity of the eigenvalue problem, we can also induce a dynamical system for $\hat{Z}$ on projective space: 
\begin{align}
\dot{\hat{Z}} &= \hat{B}_0 (\hat{Z}, \zeta,\lambda).
\end{align}

The goal now is to compare  two solutions, denoted $\hat{Z}_*$ and $\hat{Z}_0$, defined via a singular limit of a solution of the full system and a solution of the reduced linear system, respectively (strictly speaking, we compare their embeddings in $\mathbb{CP}^3$). The solution $\hat{Z}_0$ is the unique solution of the projectivized linearized slow system \eqref{eq:projlinslow}. On the other hand, the solution $\hat{Z}_*$ is defined as the projection of the full solution $E(\zeta,\lambda,\eps)$ onto the slow manifolds in the singular limit $\eps \to 0$; i.e. take $Z = \Pi E$ and choose a subsequence $\eps_n \to 0$ so that the limits
\begin{align*}
\hat{Z}(\pm 1,\lambda,\eps_n) &= \hat{Z}^{\pm}
\end{align*}
 both exist, and then define two solutions $\hat{Z}^{\pm}(\zeta,\lambda)$ of the reduced problem  \eqref{eq:projlinslow} on $\zeta < 0$ and $\zeta > 0$ with initial conditions $\hat{Z}^{\pm}$, setting
 \begin{align*}
\hat{\zeta}_*(\zeta,\lambda) &=     \begin{cases}
      \hat{\zeta}^-(\zeta,\lambda) & \text{if} \zeta < 0 \\
       \hat{\zeta}^+(\zeta,\lambda) & \text{if} \zeta > 0.
    \end{cases}
\end{align*}

  Finally, define the following two inclusions:
\begin{equation}
\begin{aligned}
E_0(\zeta,\lambda) &=  \iota_0(\zeta,\lambda) Z_0(\zeta,\lambda)\\
E_*(\zeta,\lambda) &=  \iota_0(\zeta,\lambda) Z_*(\zeta,\lambda).
\end{aligned}
\end{equation}

We re-emphasize that the starred data is obtained from approximation of the full solution, whereas data denoted by `$0$' is obtained only from solutions of the reduced problem on $S^{a,-}_0\cup F$ and $S^{a,+}_0$, concatenated by an algebraic jump condition.\\

 An approximation theorem can now be stated.

\begin{thr} \label{thm:approx}
There exists $\bar{\eps} > 0$ such that for $0 < \eps < \bar{\eps} $ the following are true:

\begin{itemize}
\item[(a)] $\lim_{\eps \to 0}\hat{E}(\zeta,\lambda,\eps) = \hat{E}_*(\zeta,\lambda)$ uniformly in $a \leq |\zeta| \leq A$, for each $0 < a < 1 < A$.
\item[(b)] Let  $Z^{0-}_*(\lambda) := \lim_{\zeta \to 0^-} Z_*(\zeta,\lambda)$. Then
\begin{align}
J_{\lambda}(Z^{0-}_*(\lambda)) = Z_*(0,\lambda) = \lim_{\zeta\to 0^+} Z_*(\zeta,\lambda).
\end{align}

\item[(c)] Suppose that $\lambda$ is {\bf not} a slow eigenvalue. Then 
\begin{align*}
\hat{Z}_*(\zeta,\lambda) &= \hat{Z}_0(\zeta,\lambda),
\end{align*}
where 
\begin{align}
\hat{Z}_0(\zeta,\lambda) &= \begin{cases}
        \hat{Z}_L(\zeta,\lambda) \text{ for } \zeta < 0
        \\
        \hat{Z}_R(\zeta,\lambda) \text{ for } \zeta \geq 0 
        \end{cases}
\end{align}

is defined by projectivizations of a pair of solutions defined on $S^{a,-}_0\cup F$ resp. $S^{a,+}_0$ that tends to the projectivization of the unstable eigenvector in either direction. The pair also satisfies the projectivized jump condition from \eqref{eq:projectivejump}:
\begin{align*}
\hat{Z}_R(0,\lambda) &= s(\bar{u}_J,\lambda).
\end{align*}
Furthermore, $\hat{E}(\zeta,\lambda,\eps) \to \hat{e}_f^+(\zeta,\lambda,\eps)$ as $\zeta \to +\infty$.
\end{itemize}

\end{thr}

\begin{remk}
This theorem should be directly compared to the analogous convergence theorem 5.3 in \cite{GJ}. Observe that the continuity condition in part (b) is replaced by a jump condition. In general we no longer expect the projected limiting solution $Z_*$ to be continuous at $\zeta = 0$ (however, it is still right-continuous).
\end{remk}

The subsections \ref{sec:estimate1}, \ref{sec:jumpmapconv}, and \ref{sec:estimate2} are devoted to proving  Theorem \ref{thm:approx} (a),(b), and (c), respectively. In the final subsection \ref{sec:homotopy}, we use this theorem to construct a homotopy between the bundles $\mathcal{E}_{\eps}(K)$ and $\mathcal{E}_0(K)$.

 ~\\

\subsection{Preliminary estimates} \label{sec:preliminaries}

Before proving Theorem \ref{thm:approx}, we need some preliminary estimates that will allow us to control the dynamics of the eigenvalue problem within the slow manifolds, and as they connect from the fast layers to the slow manifolds. The key points are that although we do not have access to an elephant trunk lemma, (i) we may still define partial relatively invariant sets over the attracting branches $S^{a,\pm}_{\eps}$ of the critical manifold $S_{\eps}$ for sufficiently small values of $\eps > 0$, and (ii) the way in which the eigenvalue parameter $\lambda$ enters the equations suggests that we can retain control of the dynamics using {\it exchange lemma}-type estimates. Let us highlight that the exchange lemma estimates are  apparently new relative to the techniques introduced in \cite{GJ}. Furthermore, the use of relatively invariant set theory applied to the case of two-plane dynamics (i.e. the construction of a so-called {\it slow elephant trunk} over a slow subbundle) also appears to be new. \\

\subsubsection{Fast and slow elephant trunks over $S^{a,\pm}_{\eps}$} \label{sec:elephanttrunks}
 
In this section, we construct relatively invariant attracting sets for \eqref{eq:linreversefast} defined within compact neighborhoods of $S^{a,\pm}_{\eps}$. Here we closely follow the treatment of Sec. IV in \cite{GJ}. Let $$\Omega = \cup_{\xi \in I} \Omega (\xi) \times \{\xi\}$$ denote a subset of $\mathbb{C}^2 \times \mathbb{R}$ such that $I$ is an open interval in $\mathbb{R}$ and $\Omega(\xi)$ is a neighborhood in $\mathbb{C}^2$ for each $\xi$, such that $\partial \Omega \cap \mathbb{C}^2 \times I$ is a smooth manifold (with $\Omega(\xi)$ varying smoothly in $\xi$). Furthermore, consider a sufficiently smooth nonautonomous system on $\mathbb{C}^2$ of the form 
\begin{align} \label{eq:arbnonaut}
\beta' &= G(\beta,\xi,\eps).
\end{align}

\begin{defi} \label{def:relinvset}
The set $\Omega$ is {\it (positively) invariant relative to $I$} if for any solution $\beta(\xi)$ of \eqref{eq:arbnonaut} with $\beta(\xi_0) \in \Omega(\xi_0)$ for some $\xi_0 \in I$, we have $\beta(\xi) \in \Omega(\zeta)$ for all $\zeta \geq \zeta_0$ for which $\zeta \in I$.
\end{defi}

 We will define two distinct collections of relatively invariant sets over the disjoint branches $S^{a,\pm}_{\eps}$, referring to one type as  {\it fast elephant trunks}, denoted $\Omega_{\pm}^f$, and to the other as {\it slow elephant trunks}, denoted $\Omega_{\pm}^s$. \\
 
\begin{figure}[t!] 
\centering
(a)\includegraphics[width=0.46\textwidth]{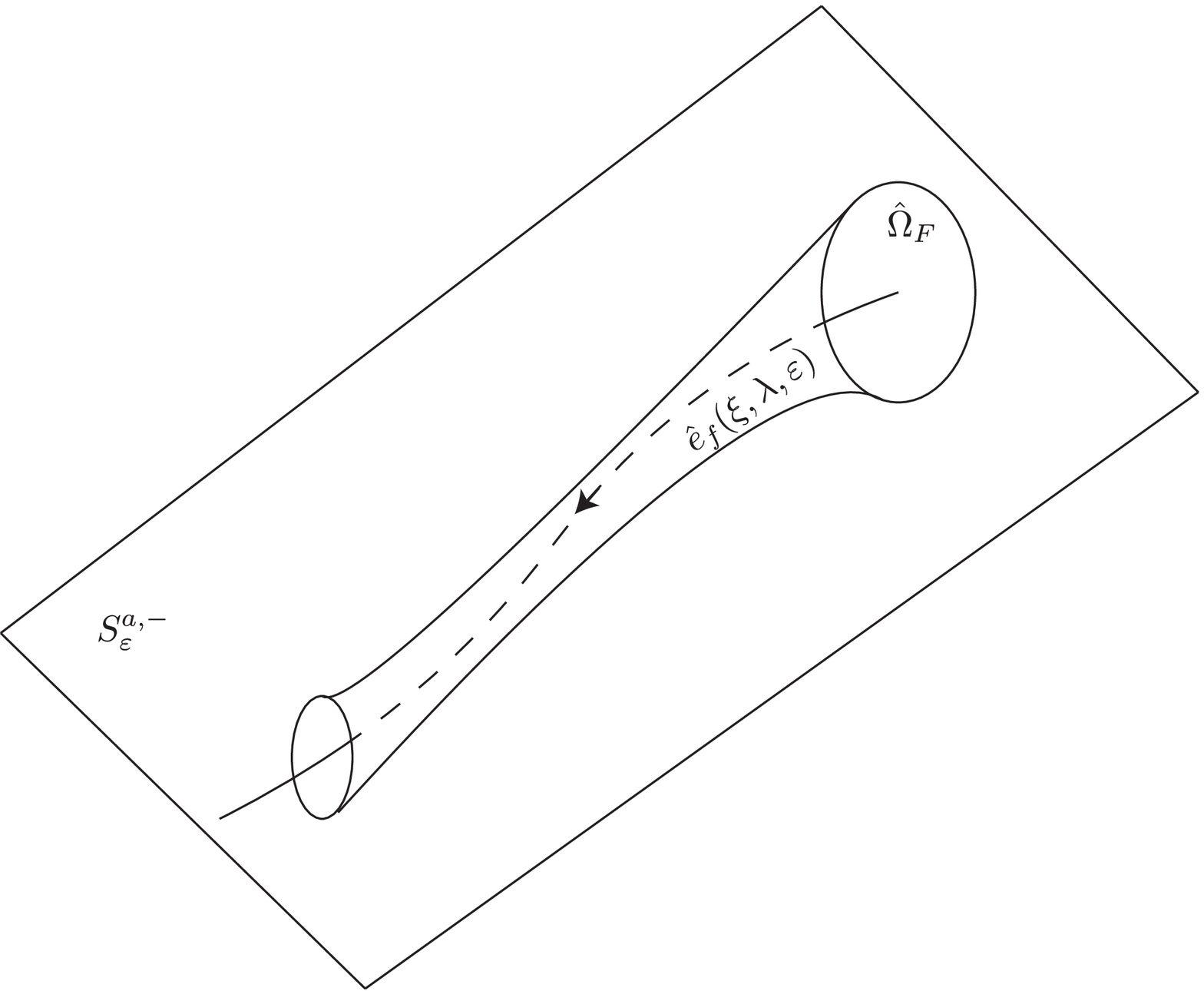}
(b)\includegraphics[width=0.46\textwidth]{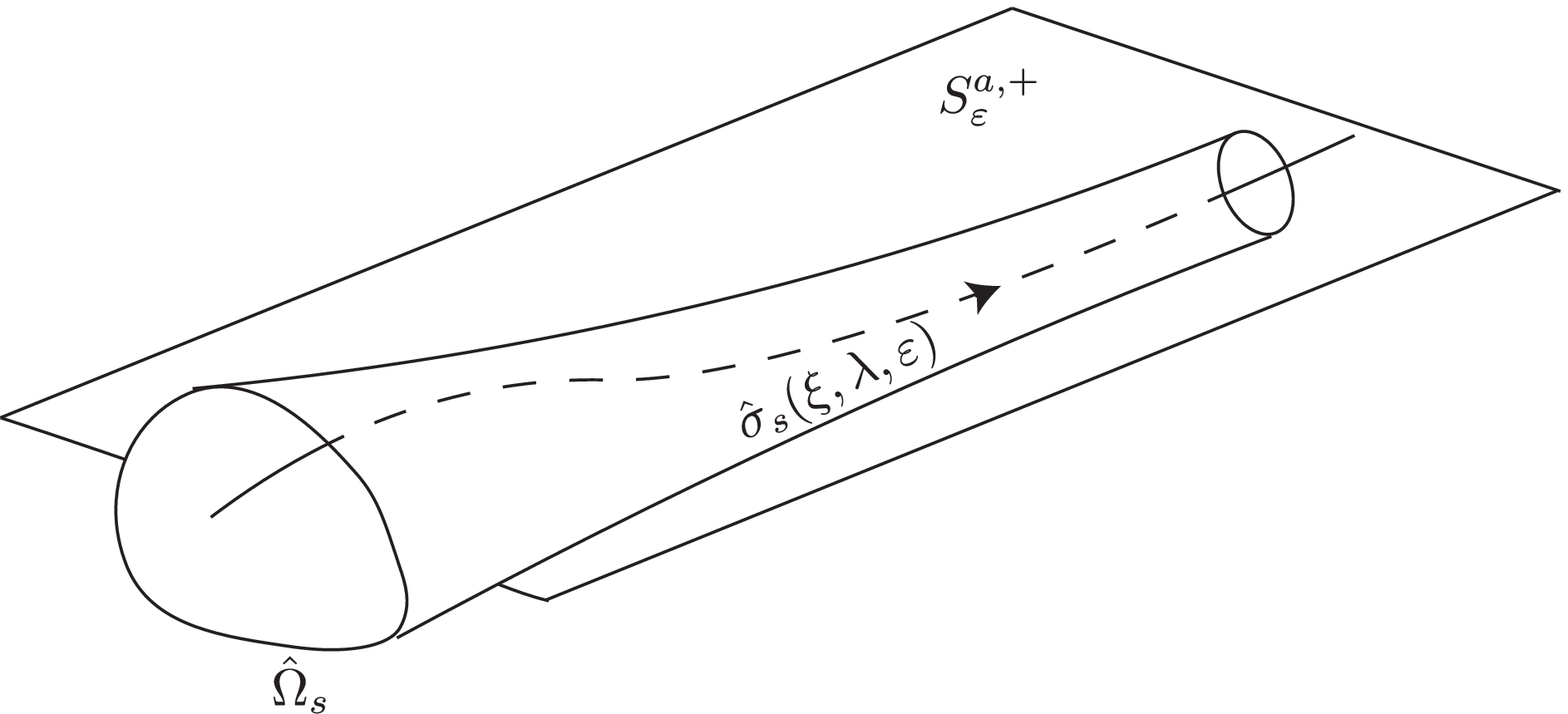}
\caption{Sketches of the (a) fast resp. (b) slow elephant trunks constructed in Lemmas \ref{lem:fasttrunks}--\ref{lem:slowtrunks}.}
\label{fig:elephanttrunks}
\end{figure}

  The starting point is to consider the auxiliary {\it autonomous} family of frozen systems associated with \eqref{eq:arbnonaut}:
  \begin{align} \label{eq:arbnonautfreeze}
\frac{d\beta}{d\xi} &= G(\beta,\gamma,\eps)
\end{align}

and to assume that \eqref{eq:arbnonautfreeze} admits  a smooth curve of critical points $\beta_0(\gamma,\eps)$ (i.e. for each $\eps >0$, $\beta_0$ depends smoothly on $\gamma$). Thinking of \eqref{eq:arbnonautfreeze} as the frozen family corresponding to a projectivization of a linearized system over a reference chart, such curves arise naturally by continuing the eigenvectors of the linear system along the travelling wave, which is parametrized by $\gamma$. We will consider subsets of the following form as candidates for elephant trunks:	
 \begin{align}
\Omega(d,\eta,\eps) &= \{(\beta,\eta): |\beta - \beta_0(\gamma,\eps)|_{\gamma} < \eta, \gamma \in I(d)\},
\end{align}

where the real parameter $d$ parametrizes a nested family of subintervals $I(d)$ (i.e. $I(d_1) \subset I(d_2)$ whenever $d_1 < d_2$), and the metric $|\cdot|_{\gamma}$  is defined according to the construction in Sec. IV. B. in \cite{GJ}. As noted there, the parameter $d$ can be thought of as characterizing the slowly varying character of \eqref{eq:arbnonaut} (with respect to $\xi$ when $\eps > 0$ is sufficiently small), while $\eta$ characterizes the width of the tube. See Fig. \ref{fig:elephanttrunks} for a sketch of the elephant trunks we construct.

 We begin with the construction of the fast elephant trunks which are defined for \eqref{eq:visclinfast}, i.e. \eqref{eq:linreversefast} with the frame variable having reversed orientation $\xi \mapsto -\xi$. Strictly speaking, we work with the (orientation-reversed)  projectivized system \eqref{eq:projsystemu0}. \\
 
 \begin{lem} \label{lem:fasttrunks}
For each $\lambda \in \Omega$ (see \eqref{eq:omegadef}), there exist $\eps_0 > 0$ and $d_0 > 0$ such that for each $\eps \in (0,\eps_0]$ and $d \in (0,d_0]$, there is a pair of fast elephant trunks $\Omega_{-}^f := \Omega_{-}^f(d,\eta,\lambda,\eps)$ (over $S^{a,-}_{\eps}$) and $\Omega_{+}^f :=\Omega_{+}^f(d,\eta,\lambda,\eps)$ (over $S^{a,+}_{\eps}$); namely, $\Omega_{-}^f$ (resp. $\Omega_{+}^f$) is positively invariant relative to the nested families of sub-intervals $I_-(d)$ (resp. $I_+(d)$), where
$$I_-(d) =\{\xi:\xi \leq \xi_-(\delta(d))<0\}$$  and 
$$I_+(d) =\{ \xi: \xi \geq \xi_+(\delta(d))>0\},$$ 
where $\delta:[0,d_0) \to [0,\infty)$ varies smoothly in $d$ with $\delta(0) = 0$,  and $\xi_{\pm}(\delta(d))$ are defined so that:
\begin{itemize}
\item  the travelling wave $\bar{x}(\xi,\eps)$ lies $\delta$-close to the singular limit of the travelling wave for all $\xi\leq \xi_-(\delta)\leq -\bar{\xi}$ and for all $\xi \geq \xi_+(\delta) \geq \bar{\xi}$, where $\bar{\xi}$ is defined so that:
 \item the matrix $A(\xi,\lambda,\eps) := G_{\beta}(\beta_0(\xi,\lambda,\eps),\xi,\lambda,\eps)$ has eigenvalues $\mu_2 - \mu_1$ and $\mu_3 - \mu_1$ uniformly negative and bounded away from 0 for all $|\xi| \geq \bar{\xi}$, where $G$ denotes the orientation reversal of the projectivized system \eqref{eq:projsystemu0}.
 \end{itemize}
 \end{lem}

{\it Proof:} See Appendix \ref{app:trunks}.\\

Lemma \ref{lem:fasttrunks} provides elephant trunks about the fast unstable directions near the slow manifolds in reverse time. In fact we will only require the fast elephant trunk near the fixed point at $\bar{u} = 0$ (i.e. over the branch $S^{a,-}_{\eps}$), in order to verify one part of a  uniform closeness estimate of the unique solution $E(\xi,\lambda,\eps)$ of \eqref{eq:linreversefast} with $|E(0,\lambda,\eps)| = 1$ and $E(\xi,\lambda,\eps) \to u^{-}$ as $\xi \to -\infty$. As we track $E(\xi,\lambda,\eps)$ for $\xi$ sufficiently large, we will also require a slow elephant trunk over $S^{a,+}_{\eps}$ which guarantees uniform closeness of $E$ near the slow subbundle. We work with the  projectivized system \eqref{eq:proj2planes}.

 \begin{lem} \label{lem:slowtrunks}
For each $\lambda \in \Omega$ (see \eqref{eq:omegadef}), there exist $\eps_0 > 0$ and $d_0 > 0$ such that for each $\eps \in (0,\eps_0]$ and $d \in (0,d_0]$, there is a pair of slow elephant trunks $\Omega_{-}^s := \Omega_{-}^s(d,\eta,\lambda,\eps)$ (over $S^{a,-}_{\eps}$) and $\Omega_{+}^s :=\Omega_{+}^s(d,\eta,\lambda,\eps)$ (over $S^{a,+}_{\eps}$) for system \eqref{eq:proj2planes}; namely, $\Omega_{-}^s$ (resp. $\Omega_{+}^s$) is positively invariant relative to the nested families of sub-intervals $I_-(d)$ (resp. $I_+(d)$), where
$$I_-(d) =\{\xi:\xi \leq \xi_-(\delta(d))<0\}$$  and 
$$I_+(d) =\{ \xi: \xi \geq \xi_+(\delta(d))>0\},$$ 
where $\delta:[0,d_0) \to [0,\infty)$ varies smoothly in $d$ with $\delta(0) = 0$,  and $\xi_{\pm}(\delta(d))$ are defined so that:
\begin{itemize}
\item  the travelling wave $\bar{x}(\xi,\eps)$ lies $\delta$-close to the singular limit of the travelling wave for all $\xi\leq \xi_-(\delta)\leq -\bar{\xi}$ and for all $\xi \geq \xi_+(\delta) \geq \bar{\xi}$, where $\bar{\xi}$ is defined so that:
 \item the matrix $A(\xi,\lambda,\eps) := G_{\beta}(\beta_0(\xi,\lambda,\eps),\xi,\lambda,\eps)$ has eigenvalues uniformly negative and bounded away from 0 for all $|\xi| \geq \bar{\xi}$.
 \end{itemize}
 \end{lem}
 
 {\it Proof:} The steps are essentially identical to those shown in Appendix \ref{app:trunks} for Lemma \eqref{lem:fasttrunks}. We point out that the corresponding family of frozen systems for \eqref{eq:proj2planes} admits a curve of attracting critical points, as described in the paragraph below \eqref{eq:proj2planes}. The remaining conditions are checked with direct calculation.  $\Box$

\begin{remk}
The slow elephant trunks in Lemma \ref{lem:slowtrunks} define attracting invariant neighborhoods relative to a metric defined in $\mathbf{Gr}(2,3)$. Switching instead to the Fubini-Study metric (which gives equivalent estimates up to a constant factor depending only on the metrics), it is a rather lengthy calculation to show that these planar neighborhoods can be in turn expressed as neighborhoods of complex lines relative to the slow subbundle. 
\end{remk}

\subsubsection{Exchange lemma-type estimates}
Per Fenichel's original setup \cite{fenichel}, we consider a singularly perturbed system of differential equations in $\mathbb{R}^n$ which admits a family of $k$-dimensional normally hyperbolic invariant manifolds $S_{\eps}$, with $0 < k < n$, for $\eps \in (0,\bar{\eps}]$, arising from a critical manifold $S_0$ when $\eps =0 $. We remind the reader that $n = 3$ and $k = 2$ in our model, and we restrict our interest to the case of {\it attracting} slow invariant manifolds. We use freely a {\it Fenichel normal form} defined over a common neighborhood $U$ of $S_{\eps}$, appended with its variational equations:
\begin{equation} \label{eq:vareqns}
\begin{aligned} 
b' &= \Gamma(b,y,\eps) b \\
y' &= \eps f(y,\eps) \\
db' &= \Gamma db + D_z \Gamma (dz \otimes b)\\
dy' &=\eps D_y f\,dy ,
\end{aligned}
\end{equation}

where $z := (b,y)$ and $dz:=(db,dy)$. In terms of the geometric objects introduced in Sec. \ref{sec:geometry}, the system \eqref{eq:vareqns} may be regarded as an induced derivation on the tangent bundle, with the dynamics on the tangent vectors coordinatized according to the Pl\"{u}cker embedding.  We highlight a few properties of Fenichel theory and the Fenichel normal form. The slow manifolds are given by $S_{\eps} = \{b = 0\}$, so that $y \in \mathbb{R}^k$ may be viewed as the slow variables. The stable manifold $W^s_{\eps}$ is foliated by invariant fibers, with associated projection map $\pi^-:W^s_{\eps} \to S_{\eps}$. \\

 \begin{figure}[t!] 
\centering
\includegraphics[width=0.95\textwidth]{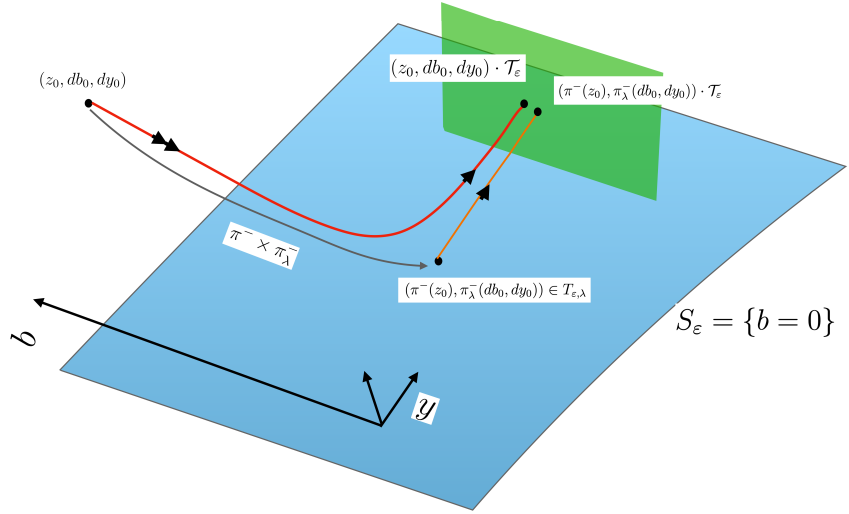}
\caption{Setup of the exchange lemma-type estimate using system \eqref{eq:vareqnseig}. Here the projection $\pi^- \times \pi^-_{\lambda}$ refers to the projection onto the invariant manifold $\{b = 0\}$, extended by the corresponding invariant manifold given by the eigenvalue problem \eqref{eq:vareqnseig}.}
\label{fig:exchange}
\end{figure}

With respect to the straightened dynamics of the normal form, this stable fibration is given by $(b,y) \mapsto y$.  Each function in \eqref{eq:vareqns} is regarded as `sufficiently' smooth for every succeeding statement to hold, and we further assume that the spectra of the matrix functions $\Gamma$ and $\eps D_y f$ satisfy the necessary conditions for Fenichel's Lyapunov type numbers. Specifically, for each $y \in U$ we list the spectrum
\begin{align*}
\text{spec}\, \Gamma(0,y,0)&:= \{\gamma_1(y) , \cdots, \gamma_{n-k}(y)\}
\end{align*}

in ascending order of their real parts, i.e. 
\begin{align*}
\text{Re}\,(\gamma_1(y)) \leq \cdots \leq \text{Re}\,( \gamma_{n-k}(y)),
\end{align*}
and we assume a spectral gap $\text{Re}\,\gamma_{n-k}(y) < 0 $ uniformly for $y \in \bar{U}$. Define
\begin{align} \label{eq:supgamma}
\gamma_0 &= \text{sup}_{y \in U} \gamma_{n-k}(y).
\end{align}

Let us summarise the geometric ideas underlying the exchange lemma. We consider families of manifolds that enter neighborhoods of normally hyperbolic slow invariant manifolds in a generic way. The primary goal is to track both the position and orientation (i.e. the corresponding evolution of tangent spaces according to the variational equations) as trajectories enter near the stable foliation and spend sufficiently long times near the slow manifold. If the entry manifolds intersect the stable foliation transversally, this generic entry is `exchanged' for exponential closeness, in both position and orientation, to the unstable foliation (or in our case, to the tangent space of the slow manifold) at exit.\\

 This geometric result can be expressed explicitly in terms of the rectified variables in  \eqref{eq:vareqns}. We follow the treatment by Jones and Tin in \cite{jonestin}.  Fix a parameter $\Delta > 0$ and  consider a one-parameter family of entry manifolds $\{M_{\eps}\}$ intersecting the stable manifold transversally at some section $\{|b| = \Delta\}$. We focus on initial conditions within the entry manifolds that remain within the box specified by 
 \begin{align} \label{eq:boxdef}
B_{\Delta} &:=  \{|y| \leq \Delta, |b| \leq \Delta\},
\end{align}

 leaving $B_{\Delta}$ through $\{|y| = \Delta\}$ only after a sufficiently long time, say $\cT_\eps = 1/\eps^s$ for some $0 < s < 1$. Let $z(t)$ parametrize an incoming trajectory $\gamma$ corresponding to such an initial condition lying on $M_{\eps}$.  The exchange lemma assures us that for $\kappa \in (0,|\gamma_0|)$ (see \eqref{eq:supgamma}), $z(t)$ is $C^1$-$\mathcal{O}(e^{-\kappa \mathcal{T}}$)-close to $S_{\eps}$ at $t = \cT_\eps$, i.e. 
 \begin{align*}
\text{dist}(z(\cT_\eps),S_{\eps}) &= \mathcal{O}(e^{-\kappa \cT_\eps})\\
\text{dist}(T_{z(\cT_\eps)}\gamma ,T_{\pi^-(z(\cT_\eps))}S_{\eps}) &= \mathcal{O}(e^{-\kappa \cT_\eps}),
\end{align*}
where the dist$(\cdot,\cdot)$ map is defined in the usual way, after specifying some appropriate metrics. In view of the rectified coordinate system for \eqref{eq:vareqns},  $T_x S_{\eps} = \{db = 0\}$ (the tangent space of $S_{\eps}$ at $x$) for each $x \in S_{\eps}$. Complete details about the formulation of the problem and a statement of the Exchange Lemma are given in \cite{jonestin}. This general formulation is referred to as the $(k+\sigma)$-{\it Exchange Lemma}, where $k$ denotes the dimension of the local unstable manifolds of the critical manifold. In the present case $k =0$, the Exchange Lemma is referred to as an {\it inclination lemma} \cite{brunovsky}.   \\

The first step is to write down the `straightened' eigenvalue problem. Although the coordinate change---i.e. the composition of straightening diffeomorphisms---which defines the Fenichel normal form is usually highly nonlinear and difficult to write down, we highlight that the {\it linearisation} of this composition acts linearly on the dynamical system defined on the tangent spaces.\footnote{Indeed, if the eigenvalue problem in the original coordinates is written as $dx' = A(x,\eps) dx + \eps \lambda M \cdot dx$, where $A(x,\eps,\lambda)dx$ denotes the original variational equations and $M$ is a constant square matrix, and $z = \varphi(x,\eps)$ is the diffeomorphism giving rise to the Fenichel normal form, then $\widetilde{dz}' = B(z,\eps)\widetilde{dz} + \eps \lambda \widetilde{M}\cdot \widetilde{dz}$, where $B\,\widetilde{dz}$ is the transformed variational equation and $\widetilde{M} = D\varphi \circ M \circ D\varphi^{-1}$.}  Hence, we obtain
\begin{equation} \label{eq:vareqnseig}
\begin{aligned} 
b' &= \Gamma(b,y,\eps) b \\
y' &= \eps (f(y,\eps) +H(b,y,\eps)b)\\
\widetilde{db}' &= \Gamma \widetilde{db} + D_z \Gamma (\widetilde{dz} \otimes b) + \eps \lambda G_1(b,y,\eps)\cdot \widetilde{dz}\\
\widetilde{dy}' &=\eps( D_y f\,\widetilde{dy} + D_z H (\widetilde{dz}\otimes b) + H \widetilde{db}) + \eps \lambda G_2(b,y,\eps)\cdot \widetilde{dz}.
\end{aligned}
\end{equation}
Here, $\widetilde{dz} :=(\widetilde{db},\widetilde{dy})$ and $G_1,\,G_2$ are smooth functions that are defined from the $\eps$-family of diffeomorphisms used to derive the Fenichel coordinates, and so are independent of $\lambda$.  Let us write the transformed eigenvalue problem in \eqref{eq:vareqnseig} in the more compact form
\begin{align}
\widetilde{dz}' &= A(z,\eps)\widetilde{dz} + \eps \lambda \widetilde{M}(z,\eps) \widetilde{dz},
\end{align}

with $\widetilde{M}_{\Delta,\eps_0} := \max_{z\in B_{\Delta},0 \leq \eps \leq \eps_0}\widetilde{M}(z,\eps)$. We drop the tilde notation for the remainder of the section for ease of reading. For $T> 0$, let $Q_T$ denote the set of initial conditions within the box of width of $\Delta>0$ so that for each $q := z(0) = (b_0,y_0) \in Q_T$, we have that $z(t)$ remains within the $\Delta$-box for each $t \in [0,T]$. 

\begin{lem} \label{lem:exchange}
Assume the hypotheses of the $(k+\sigma)$ Exchange Lemma in the singularly-perturbed case (Theorem 6.7 in \cite{jonestin}) with $k = 0$ (i.e. the normally hyperbolic critical manifold does not admit fast unstable directions). Furthermore, let $(\gamma(t,\eps),dz(t,\eps))$ be a one-parameter family of trajectories of \eqref{eq:vareqnseig} (where the components of the tangent vector are given by $dz(t,\eps) = (db(t,\eps),dy(t,\eps)$), so that $\gamma(t,\eps)$ satisfies the hypotheses for the invariant manifolds $M_{\eps}$ in \cite{jonestin}. Furthermore, let $dz(0,\eps) =: (db_0(\eps),dy_0(\eps))$ satisfy the estimate
\begin{align*}
|dy_0| \geq M \eps
\end{align*}
for some $M > 0$ independent of $\eps$, when $\eps > 0$ is sufficiently small, and fix $\beta$ with $0 < \beta < 1$. Then there exists $\bar{\eps} > 0$ so that for $0 < \eps \leq \bar{\eps}$,
\begin{align*}
\text{dist}(dz(\cT_\eps), T_{\gamma(\cT_\eps)}S_{\eps}) \leq K_2 \eps^{\beta}
\end{align*}

$K_2$ depending only on $\bar{\eps}$ and $\Delta$.
\end{lem}

{\it Proof:} See Appendix \ref{app:exchange}.

\begin{remk} \label{remk:exchange}
Let us highlight a few points about Lemma \ref{lem:exchange} relative to the $(k+\sigma)$ Exchange Lemma in the singularly-perturbed case. The first obvious difference is that the estimate is weaker, being only linear in order of $\eps$. This is because we compare a generic incoming trajectory to the tangent space of the slow manifold (i.e. the set $\{b=0,\,db=0\}$) instead of the slow subbundle, where the computation becomes more difficult. For our present purposes this weaker estimate is sufficient, since it is immediate that the slow subbundle is $\mathcal{O}(\eps)$-close to the tangent space of the slow manifold relative to the Fubini-Study metric (note: it is crucial here that the eigenvalue problem depends `weakly' on $\lambda$, i.e. through terms of $\mathcal{O}(\eps\lambda)$ only). \\

However, we highlight that the system \eqref{eq:vareqnseig} is singularly-perturbed, with the identical normally hyperbolic critical manifold  $\{b = 0,\,db=0\}$ as that of the variational equations. By standard Fenichel theory, there exists a one-parameter family of (real) six-dimensional locally invariant attracting slow manifolds $\mathcal{M}_{\eps}$ for sufficiently small values of $\eps>0$, which lie $\mathcal{O}(\eps)$-close in Hausdorff distance to the critical manifold; two of the dimensions come from the slow directions in the phase space and the remaining four come from the slow complex directions of the eigenvalue problem. We therefore expect that the tangent spaces of incoming trajectories under the flow of the eigenvalue problem align exponentially closely to that of $\mathcal{M}_{\eps}$, i.e. there is an analogue to the standard Exchange Lemma estimate for the eigenvalue problem. Furthermore, the slow subbundle should provide the $\mathcal{O}(\eps)$ correction of the tangent bundle of the slow manifold (i.e. it is $\mathcal{O}(\eps^2)$-close with respect to the Fubini-Study metric), in analogy to the standard computations using Fenichel theory. We illustrate this in Appendix \ref{app:toy} with a toy problem.
\end{remk}

~\\
\subsection{Estimates near the slow subbundle} \label{sec:estimate1}

In this subsection we prove Theorem \ref{thm:approx}\,(a).  Our aim is to prove a similar result to Corollary 5.6 in \cite{GJ}, namely that the projectivization $\hat{E}(\zeta,\lambda,\eps)$ of any nontrivial solution $E(\zeta,\lambda,\eps)$ lies uniformly close to the slow subbundle when the wave is near to the slow manifolds. A new technical issue here is that we do not have an elephant trunk estimate over the fast layer. We will instead combine our existing `partial' fast and slow elephant trunks with our exchange lemma-type estimate to achieve this uniform closeness.  \\

In the coming analysis we will work with projectivisations of the slow subbundle, denoted $\hat{\sigma}_s$ and $\hat{\Sigma}_s$ as usual. Now fix any metric $\rho$ on $\mathbb{CP}^3$.
\begin{defi} \label{def:bundlenbd}
For any set $\hat{S}\subset \mathbb{CP}^3$ and $\delta > 0$, a {\it $\delta$-neighborhood} of $\hat{S}$ is the set
\begin{align}
N_{\delta}(\hat{S}) &= \{\hat{y} \in \mathbb{CP}^3: \rho(\hat{s},\hat{y}) < \delta \text{ for some }\hat{s} \in \hat{S}\}.
\end{align}
\end{defi}

~\\
We begin by estimating the closeness of $E(\zeta,\lambda,\eps)$ to $\Sigma_s(\zeta,\lambda,\eps)$ near $S^{a,\pm}_0$. To do so, we consider the following family of autonomous frozen systems (and their projectivizations) corresponding to the linearized problem: 

\begin{equation}
\begin{aligned}  \label{eq:frozen}
y' &= a(\gamma,\lambda,\eps) y \\
\hat{y}' &= \hat{a}(\hat{y},\gamma,\lambda,\eps). 
\end{aligned}
\end{equation}

We refer the reader to a series of technical lemmas in Appendix \ref{app:estimates}, culminating in the uniform estimate Corollary \ref{coro:slowsubbundleapprox} which is used in the following proof.\\

 {\it Proof of Theorem \ref{thm:approx}(a).} Fix two constants $a$ and $A$ with $ 0 < a <  1 < A$ and suppose $\zeta \in [a,A]$.   We select a representative solution $E \in \pi^{-1} \hat{E}$ scaled so that slow components $|Z| \leq 1$ on $a \leq \zeta \leq A$ for each $\eps \in (0,\bar{\eps})$. It quickly follows that for $\zeta = 1$ fixed,  we have $Z \to Z^+$ as $\eps \to 0$. Our primary task is to compare a solution $Z$ representing $\hat{Z}$ to a solution $Z_*$ representing the corresponding projectivisation of the reduced problem, i.e. we show that $|Z-Z_*| \to 0$ uniformly on $a \leq \zeta \leq A$ as $\eps \to 0$. \\
 
 Away from the jump, the projection onto the slow components $Z = (P,V)$ and $Z_R = (P_R,V_R)$ of the full and reduced systems, respectively, may be written as
 
 \begin{align*}
\dot{Z} &= B(\zeta,\lambda,\eps) Z + G\Gamma\\
\dot{Z}_R &= B_R(\zeta,\lambda) Z_R,
 \end{align*}

where (suppressing the dependence of the phase space coordinate $\bar{U}$ on $\zeta$ and $\eps \geq 0$):
\begin{align*}
B(\zeta,\lambda,\eps) &= \begin{pmatrix} 0 & \frac{R'(\bar{U})-\lambda}{D(\bar{U})} \\ -1 & \frac{c}{D(\bar{U})}\end{pmatrix} \\
G(\zeta,\lambda,\eps) &= R'(\bar{U})-\lambda\\
\Gamma(\zeta,\lambda,\eps) &= U - \frac{V}{D(\bar{U})}
\end{align*}
and
\begin{align*}
B_R(\zeta,\lambda) &= \begin{pmatrix} 0 & \frac{R'(\bar{U})-\lambda}{D(\bar{U})} \\ -1 & \frac{c}{D(\bar{U})} \end{pmatrix}.
\end{align*}

Here the term $G\Gamma$ can be thought of as a forcing term that is `turned on' when $\eps > 0$; specifically, $\Gamma$ measures how far (the slow projection of) the solution $Y$ is from solving the reduced eigenvalue problem. Evidently $B(\zeta,\lambda,\eps) \to B_R(\zeta,\lambda)$ as $\eps \to 0$, uniformly in the interval $a \leq \zeta\leq A$. It remains to estimate the forcing term $G\Gamma$ uniformly within this $\zeta$-interval. The term $G$ has a uniform singular limit for $\eps  \to 0$, so we concern ourselves with the singular limit for $\Gamma$. For $\delta > 0$ arbitrarily chosen, Corollary \ref{coro:slowsubbundleapprox} provides a sufficiently small $\bar{\eps} > 0$ such that
\begin{align*}
\hat{E}(\zeta,\lambda,\eps) \in N_{\delta} (\Sigma_s (\zeta, \lambda, \eps))
\end{align*}

for each $\eps \in (0,\bar{\eps}]$ and $\zeta \geq a$. A fixed representative $E(\zeta,\lambda,\eps)$ of the projectivised solution of the full linearized system can thus be made arbitrarily close to a linear combination of the basis vectors $F_1,\,F_2$ by varying $\delta$ (see the definitions of $f_1,\,f_2$ in \eqref{eq:slowsub}; we use capital letters to denote the appropriate timescale). Each of these basis vectors in turn has a singular limit. Indeed, we have
\begin{align} \label{eq:convergencetoreducedbundle}
E(\zeta,\lambda,\eps) \to \alpha_1(\zeta,\lambda) R_{s,1}(\zeta,\lambda) + \alpha_2(\zeta,\lambda) R_{s,2}(\zeta,\lambda)
\end{align}

as $\eps \to 0$, where $R_{s,1}(\zeta,\lambda)$ and $R_{s,2}(\zeta,\lambda)$ are the eigenvectors spanning the reduced slow subbundle. These are the required singular limits of the basis vectors $F_1,\,F_2$, and they admit explicit formulas:

\begin{equation} \label{eq:reducedeigvecszeta}
\begin{aligned}
F_1 \to R_{s,1} &= (1/D(\bar{U} ), V_{p} , 1 )^{\top}\\
F_2 \to R_{s,2} &= (1/D(\bar{U} ), V_{m} , 1)^{\top},
\end{aligned}
\end{equation} 

 and $\alpha_1(\zeta,\lambda)$ and $\alpha_2(\zeta,\lambda)$ are uniformly bounded coefficients. Note that the basis vectors $R_{s,i}$, $i = 1,2$ are derived similarly to the asymptotic case in \eqref{eq:redvecs}, and the auxiliary quantities $V_{p,m}$ are defined analogously to $\nu_{j,\pm}$ in \eqref{eq:redvecs2}. Finally, we write $\Gamma$ in a convenient form by including and rearranging terms:
\begin{align} \label{eq:gammaconvergence}
\Gamma(\zeta,\lambda,\eps) &= U(\zeta,\eps)- \frac{V}{D(\bar{U}(\zeta,0))} + \left( \frac{1}{D(\bar{U}(\zeta,\eps))} -  \frac{1}{D(\bar{U}(\zeta,0))}\right)V.
\end{align}
The term in parentheses converges uniformly to 0 and $|V| \leq 1$ by construction, and the first to terms also converge uniformly to zero by \ref{eq:reducedeigvecszeta}; hence $|\Gamma| \to 0$ uniformly on the required $\zeta$-interval, and thus $Z \to Z_*$ uniformly on this interval by Gronwall's inequality. Finally, the $U$-components of $E$ also converge uniformly to those of $E_*$ by \eqref{eq:convergencetoreducedbundle}, \eqref{eq:gammaconvergence}, and Corollary \ref{coro:slowsubbundleapprox}.  Altogether, the uniform convergence of $E$ to $E_*$ on $a \leq \zeta \leq A$ follows from Gronwall's inequality. The case $\zeta < 0$ proceeds identically, so we omit the proof. $\Box$\\

\subsection{The jump map as a singular limit} \label{sec:jumpmapconv}

We now construct the jump map used to define the slow eigenvalue problem in Def. \ref{def:sloweigs}. Our goal is to determine the fate of the slow components $y_s(\xi,\lambda,\eps) := (p(\xi,\lambda,\eps),v(\xi,\lambda,\eps))$ of $y(\xi,\lambda,\eps)$ across the fast layer. The difficulty in tracking the slow data in this inner layer is that $y_s(\xi,\lambda,\eps) $ remains $\bigO(\eps)$ while $u(\xi,\lambda,\eps)$ grows to $\bigO(1)$ for a time interval that is $\bigO(1/\eps)$; in other words, the linearised solution can be made to align arbitrarily closely to the fast fibres after crossing the fold, and they remain close throughout most of the jump. The directional information carried by the slow variables is not annihilated as $\eps \to 0$, however. The natural approach from the point of view of GSPT is to perform an $\eps$-dependent rescaling of the fast linearized equations: for $\eps > 0$, let 

\begin{equation}
\label{eq:epsrescaling}
\begin{aligned}
\eps \beta_1 &= p\\
\eps \beta_2 &= v.
\end{aligned}
\end{equation} 

Then we have

\begin{equation}
\label{eq:betaeps}
\begin{aligned}
u' &= \frac{1}{c} (\eps \beta_2 - (\eps \lambda + D(\bar{u})) u)\\
\beta_1' &= (R'(\bar{u})-\lambda)u\\
\beta_2' &= cu - \eps \beta_1,
\end{aligned}
\end{equation} 

where a factor of $\eps$ has been cancelled from the latter two equations. The resulting equations limit to the following linear system as $\eps \to 0$:

\begin{equation}
\label{eq:beta0}
\begin{aligned}
u' &= -\frac{1}{c}  D(\bar{u})u\\
\beta_1' &= (R'(\bar{u})-\lambda)u\\
\beta_2' &= cu.
\end{aligned}
\end{equation} 

The first equation in \eqref{eq:beta0} is the variational equation of the linearized layer problem \eqref{eq:linlayer}, and thus has a family of nontrivial bounded solutions $u(\xi) = K(\bar{v}_0 - F(\bar{u}(\xi))) = K (d\bar{u}/d\xi)$, where $\bar{v}_0$ is the $\bar{v}$-component of $F_-$. It is then possible to calculate $\beta_1(\xi)$ and $\beta_2(\xi)$ explicitly. We find it more convenient to perform the intermediate calculations in the projective space $\mathbb{CP}^2$. This allows us to fix a free parameter by applying a smoothness condition across the fold.\\

We consider the projection of the linear system \eqref{eq:beta0} on a copy of $\mathbb{CP}^2$, choosing the chart $(s,g) = (\beta_1/\beta_2, u/\beta_2)$, $\beta_2 \neq 0$:

\begin{equation}
\label{eq:beta0cp2}
\begin{aligned}
s' &= (R'(\bar{u})-\lambda)g - c sg\\
g' &= -\frac{D(\bar{u})}{c}g - (R'(\bar{u}-\lambda))g^2.
\end{aligned}
\end{equation} 

The $g$-equation decouples, so we solve it directly to find
\begin{align*}
g(\xi) &= \frac{\bar{v}_F - F(\bar{u}(\xi))}{c\bar{u}(\xi) + C},
\end{align*}

where $C \in \mathbb{C}$ is a constant. The $s$ equation is then given by
\begin{align}
s' &= ((R'(\bar{u})-\lambda) - c s)\frac{\bar{v}_F - F(\bar{u}(\xi))}{c\bar{u}(\xi) + C}. \label{eq:fiberslow}
\end{align}

\begin{figure}[t!] 
\centering
\includegraphics[width=1.0\textwidth]{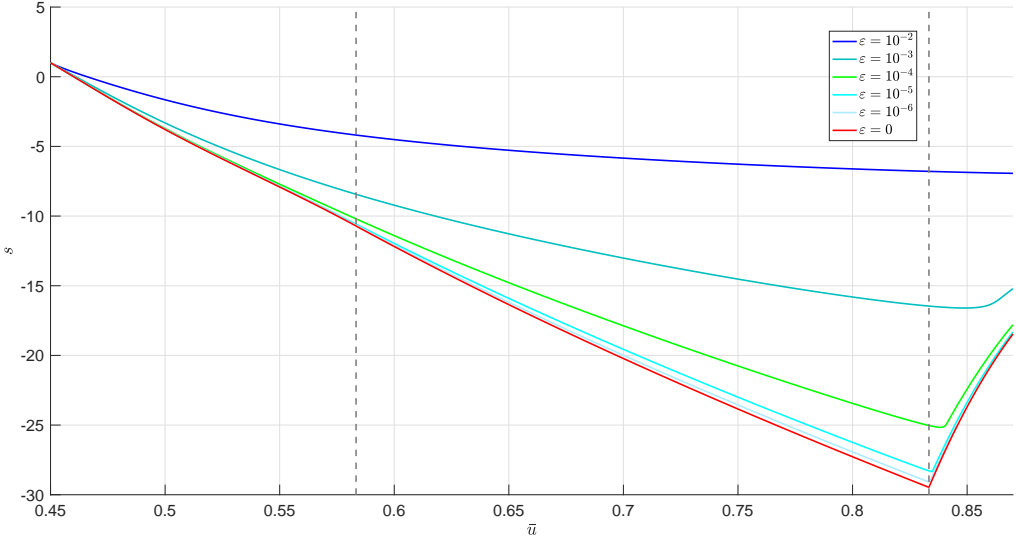}
\caption{Demonstration of the convergence of the projectivised dynamics (blueish curves) onto the reduced dynamics (red curve) as $\eps \to 0$ for $\lambda = 15$. Solutions of reduced dynamical system are defined by concatenating the solution of the projectivized flow of $dS/d\zeta$ (see \eqref{eq:compare}), defined for $\bar{U} \in (0,\bar{U}_F) \cup (\bar{U}_J,1)$,  with the graph of the map \eqref{eq:projectivejump} defined for $\bar{u} \in [\bar{u}_F,\bar{u}_J]$. The left dotted line denotes the fold $\bar{u} = \bar{u}_F$ and the right dotted line denotes the jump curve $\bar{u} = \bar{u}_J$ (and hence the fast dynamics takes place in between these dotted lines). The right dotted line is the jump curve, where the jump condition is used to concatenate the appropriate solutions. The wave speed was held at $c= c_{0} \approx 0.199362$ for each value of $\eps$.}
\label{fig:lambda15}
\end{figure}
Equation \eqref{eq:fiberslow} is now a one-dimensional nonautonomous problem and can be solved explicitly, but we must first fix the constant $C$ to define an unambiguous jump condition. Such an ambiguity arises even at the projective level. For each $\eps > 0$, the scale of the linearized solution along the wave is set by one free parameter, whereas in the singular limit, we split the linearized dynamics into subsystems defined along the singular heteroclinic orbit. We are able to choose freely the scale of each corresponding linearized solution segment. It follows that these free parameters of the reduced systems may be constrained so that they are compatible with the scalings of the `full' linearized solutions as $\eps \to 0$. \\

A suitable compatibility condition is the $C^1$-differentiability of the desingularized linearized flow across the fold. We compare \eqref{eq:fiberslow} with the projectivization of the desingularized slow eigenvalue problem \eqref{eq:linslow}, which we write over the chart $S = P/V$, with $V \neq 0$, as
\begin{align} \label{eq:projlinsloweqs}
\frac{dS}{d\zeta} &=   (R'(\bar{U})-\lambda-cS+D(\bar{U})S^2).
\end{align}

 Using the chain rule, we write 
\begin{equation}
\label{eq:compare}
\begin{aligned}
\frac{\del S}{\del \bar{P}}\frac{d\bar{P}}{d\zeta} + \frac{\del S}{\del \bar{U}} \frac{dS}{d\zeta}= \frac{dS}{d\zeta} &=   (R'(\bar{U})-\lambda-cS+D(\bar{U})S^2)\\
\frac{ds}{d\bar{u}} &= \frac{1}{c\bar{u}-C} (R'(\bar{U})-\lambda-cS).\\
\end{aligned}
\end{equation}  

At the point $\Gamma_0 \cap F_- = \{\bar{X}_F = (\bar{u}_F,\bar{p}_F,\bar{v}_F)\}$ where the singular orbit intersects the fold, we have $D(\bar{U}_F) = 0$ and thus $(dS/d\zeta)|_{\bar{X}_F} = (\del S/\del \bar{U})(c\bar{u}_F-\bar{p}_F)$. By matching $ (\del S/\del \bar{U})$ and $ds/d\bar{u}$ at $\bar{X}_F$ we find
\begin{align*}
C &= \bar{p}_F.
\end{align*}

At the projective level, the jump map $u_J \mapsto s(\bar{u_J})$ is defined from the solution $s(\bar{u})$ of the (complex) one-dimensional initial-value problem
\begin{align*}
\frac{ds}{d\bar{u}} &= \frac{1}{c\bar{u}-\bar{p}_f} (R'(\bar{u})-\lambda-cs)\\
s(\bar{u}_F) &= s_0.
\end{align*}

 Note that $c\bar{u} - \bar{p}_F > 0$ is bounded away from 0 for $\bar{u} \in [\bar{u}_F,\bar{u}_J]$ since $c\bar{u}_F - \bar{p}_F > 0$ by Hypothesis \ref{hyp:monotoneU} and $\bar{u} \geq \bar{u}_f$ across the jump. Explicitly, we have
\begin{align} \label{eq:projectivejump}
s(\bar{u},\lambda)&= \frac{R(\bar{u}) - R(\bar{u}_F)-\lambda(\bar{u}-\bar{u}_F)+ s_0(c\bar{u}_F - \bar{p}_F)}{c\bar{u}-\bar{p}_F}
\end{align}

for $\bar{u}_F \leq \bar{u} \leq \bar{u}_J$. See Fig. \ref{fig:lambda15} for a demonstration of the approximation of the `full' linearized flow to that of the hybrid reduced problem for a nonzero value of $\lambda$.\\

The linear jump map \eqref{eq:jumpmatrix} can be extracted from \eqref{eq:projectivejump} by using the chart map $s = p/v$. For $a,b \in \mathbb{C}^n$, we define the equivalence relation $a\sim b$ if $a = \gamma b$ for some complex number $\gamma \neq 0$. We then have 

\begin{equation} \label{eq:jumpmapu}
\begin{aligned}
(p_0,v_0) \mapsto (p,v) &= ((p/v)v,v) = (sv,v)\\
&\sim (s,1) = \left(   \frac{R(\bar{u}) - R(\bar{u}_F)-\lambda(\bar{u}-\bar{u}_F)+ s_0(c\bar{u}_F - \bar{p}_F)}{c\bar{u}-\bar{p}_F}, 1\right)\\
&\sim  \left(\frac{v_0(R(\bar{u}) - R(\bar{u}_F)-\lambda(\bar{u}-\bar{u}_F))+ p_0(c\bar{u}_F - \bar{p}_F)}{c\bar{u}-\bar{p}_F}, v_0\right)\\
&:= J(p_0,v_0,\lambda,\bar{u}).
\end{aligned}
\end{equation}
We can now define $J_{\lambda}(P,V) :=J(P,V,\lambda,\bar{U}_J)$. Our jump map is defined up to linear equivalence and is therefore clearly nonunique; however, any two linear choices induce the same projective jump map $\hat{J}_\lambda := s(\bar{u}_J,\lambda)$ on $\mathbb{CP}$. This essentially completes the proof of Theorem \ref{thm:approx}(b).

\subsection{Control from the reduced dynamics} \label{sec:estimate2}

{\it Proof of Theorem \ref{thm:approx}(c).} By Corollary \ref{coro:slowsubbundleapprox}, the unique solution $\hat{E}(\zeta,\lambda,\eps)$ which tends to the unstable eigenvector at $\hat{u} = 0$ as $\zeta \to -\infty$  remains uniformly close to the (projectivized) slow subbundle $\hat{\Sigma}_s$. \\

Fix an interval $a \leq |\zeta| \leq A$, where $0< a<A$ with $A \leq \infty$. For fixed $\lambda,\eps$, let $(\beta_1(\xi),\beta_2(\xi))$ denote the local coordinates of the solution $\hat{E}(\zeta,\lambda,\eps)$ which coincide with the coordinates used to define the slow projectivized system \eqref{eq:projsystemu0slow}. Then by Lemma \ref{lem:asympt}, for any $\delta > 0$ we can find some sufficiently small $\bar{\eps}= \bar{\eps}(\delta)$ so that $|\beta_1 - 1/D(\bar{U})|\leq \delta$ for all $a \leq |\zeta| \leq A$ and for all $0 < \eps \leq \bar{\eps}$. This estimate can be read off from the expressions for the reduced eigenvectors in \eqref{eq:redvecs}.\\

Now consider the solution $\hat{Z}_*(\zeta,\lambda)$ written with the coordinate $S_*(\zeta) = P(\zeta)/V(\zeta)$ with $V \neq 0$. The dynamics of $S(\zeta)$ is determined by the system \eqref{eq:projlinreduced}. At either limit $\zeta \to -\infty$ and $\zeta \to +\infty$, the corresponding asymptotic systems admit a pair of hyperbolic (saddle) fixed points, interpreted in the projectivized system as an attractor and a repeller. Let us consider the case $\zeta \to -\infty$ and denote the projectivized attractor (corresponding to the unstable eigenvector of the saddle) by $\hat{u}^-(\lambda)$ and the projectivized repeller (stable eigenvector) by $\hat{s}^-(\lambda)$. We seek to verify that the only possibility is that $\hat{Z}_*(\zeta,\lambda) \to \hat{u}^-(\lambda)$ as $\zeta \to -\infty$ (i.e. that $\hat{E}_*(\zeta,\lambda) \to \hat{r}_{s,1}^-$, with $r_{s,1}^-$ as defined in \eqref{eq:redvecs}).\\

So suppose that $\hat{E}_*(\zeta,\lambda) \to \hat{r}_{s,2}^-$ (the weak {\it stable} eigenvector) instead. Since the corresponding critical point $\hat{s}^-$ is a repeller, let $N$ denote a repelling neighborhood of radius $\eta > 0$. Then by hypothesis, there exists some $\zeta_1 = \zeta_1(\eta)$ so that $S_*(\zeta,\lambda,\eps) \in N$ for all $\zeta \leq \zeta_1$. As long as $V \neq 0$ remains bounded away from zero, for each sufficiently small $\eps > 0$ we may choose a section $E \in \pi^{-1}\hat{E}$ so that 
\begin{align*}
E(\zeta,\lambda,\eps) &= (1/D(\hat{u}),a c_0/D(\hat{u}) + (b-a) \nu_{m,-},1) + \mathcal{O}(\delta)
\end{align*}
where $a + b = 1$. Thus, $|\beta_1|$ and $|\beta_2|$ (the norms of the coordinates of $\hat{E}$ with respect to the chart  specified by \eqref{eq:projsystemu0slow}) remain uniformly bounded while $S \in N$. Now, observe that the $\beta_2$ equation in \eqref{eq:projsystemu0slow} may be written as follows:
\begin{align*}
\dot{\beta_2} &= (R'(\bar{U}^-)-\lambda) - c\beta + D(\bar{U}^-)S^2 + d(\zeta,\lambda,\eps),
\end{align*}

where $d(\zeta,\lambda,\eps)$ consists of terms characterizing the perturbation from the asymptotic system together with terms of the form $(\beta_1 - 1/D(\bar{U}))$, which can be uniformly bounded as stated in the beginning of the proof. Hence, $N$ remains a repelling neighborhood of $\hat{s}^-$ for the system above. But this implies that $\hat{E}(\zeta,\lambda,\eps)$ can be chosen to remain in any small neighborhood of $\hat{r}^-_{s,1}(\lambda)$ as $\zeta \to -\infty$, which contradicts the fact that $\hat{E}(\zeta,\lambda,\eps)\to \hat{r}^-_{s,2}$ as $\zeta \to -\infty$. Hence, $\hat{E}_*$ and $\hat{E}_0 = \iota_0 \hat{Z}_0$ coincide for $\zeta < 0$.\\

For the matching on the right hand-side, we have by Theorem \ref{thm:approx}(b) that the jump condition is identical and uniquely determined, so in fact $\hat{E}_*$ and $\hat{E}_0$ coincide for $\zeta > 0$ as well. It remains to check that $\hat{E}$ tends uniformly to $\hat{r}_{s,1}^+$. This step proceeds identically to the preceding argument, so we omit it. $\Box$

\subsection{Continuing $\mathcal{E}_{\eps}(K)$ to $\mathcal{E}_0(K)$} \label{sec:homotopy}

We now prove the following corollary of Theorem \ref{thm:approx}:

\begin{coro} \label{thm:bundlecontinuation}
There exists $\eps_0 > 0$ such that for each $\eps \in (0,\eps_0]$, $\mathcal{E}_{\eps}(K) \cong \mathcal{E}_0(K)$. In particular, $c_1(\mathcal{E}_{\eps}(K))  = c_1( \mathcal{E}_0(K))$. 
\end{coro}

  We follow the approach of Sec. VI-D in \cite{GJ}. The method is quite natural. With respect to the compactified time $T$-scaling (see \eqref{eq:timescalingreduced}), the fast dynamics over the jump occurs over an $\mathcal{O}(\eps)$-interval straddling the midpoint of the interval $[-1,1]$. The approach is to first define two hemispheric bundles over the complement of this small interval in $[-1,1]$, using the flow itself to construct the gluing map. We then construct a homotopy which closes this gap continuously as the homotopy parameter is varied. Theorem \ref{thm:approx} is then used to show that the gluing map and bundles have nice singular limits: they are precisely the ones used to construct the reduced augmented unstable bundle $\mathcal{E}_0(K)$. \\
  
{\it Proof of Corollary \ref{thm:bundlecontinuation}:} Fix the pair of real parameters $S^{a,-}_{\varepsilon},\,S^{a,+}_{\varepsilon}$ with $-1 < S^{a,-}_{\varepsilon} < 0 < S^{a,+}_{\varepsilon} < 1$. We work with the parameter pair $S = (S^{a,-}_{\varepsilon},\sigma)$, where $\sigma \in [S^{a,-}_{\varepsilon},S^{a,+}_{\varepsilon}]$ is the homotopy parameter that we will continuously vary.  \\

For $\eps > 0$ fixed, define the bundle $\mathcal{E}(S,\eps)$ as follows. The base space $B_S$ of $\mathcal{E}(S,\eps)$ is a sphere obtained by gluing the following two hemispheres along their boundaries:
\begin{align*}
b_-(S^{a,-}_{\varepsilon}) &= B\cap \{T \leq S^{a,-}_{\varepsilon}\}\\
b_+(\sigma) &= B\cap \{T \geq \sigma\}.
\end{align*}

We define two hemispheric bundles over these base spaces as follows:
\begin{align*}
\mathcal{E}^-(S^{a,-}_{\varepsilon},\eps) &= \mathcal{E}_{\eps}(K)|_{b_-(S^{a,-}_{\varepsilon})}\\
\mathcal{E}^+(\sigma,\eps) &= \mathcal{E}_{\eps}(K)|_{b_+(\sigma)}.
\end{align*}

We complete the construction by defining a gluing map $\varphi_{S,\eps}: \mathcal{E}^-(S^{a,-}_{\varepsilon},\eps)|_{b_- \cap \,b_+}\to \mathcal{E}^+(\sigma,\eps)$. Here ``$b_- \cap b_+$'' refers to the intersection curve of the hemisphere boundaries, corresponding to a copy of $K$ in $B_S$. Let $\zeta_L$ and $\zeta_{\sigma}$ be the values of $\zeta$ corresponding to $T = S^{a,-}_{\varepsilon}$ and $T = \sigma$, respectively, defined using the initial value problem \eqref{eq:timescalingreduced}. Let $E(\zeta,\lambda,\eps)$ be the solution, which has data $E_L$ at $\zeta = \zeta_L$ (by a time translation, if necessary), where $E_L$ lies in the fiber of $\mathcal{E}^-(S^{a,-}_{\varepsilon},\eps)$ over $(S^{a,-}_{\varepsilon},\lambda)$. A suitable gluing map is then the obvious one which is induced by the flow:
\begin{align}
\varphi_{S,\eps} E_L &= E(\zeta_{\sigma},\lambda,\eps).
\end{align}
 The bundle $\mathcal{E}(S,\eps)$ is then defined as
 \begin{align}
\mathcal{E}(S,\eps) &:= \mathcal{E}^-(S^{a,-}_{\varepsilon},\eps) \cup_{\varphi_{S,\eps}} \mathcal{E}^+(\sigma,\eps).
\end{align}

Note that at $\sigma = S^{a,-}_{\varepsilon}$, the map $\varphi_{S,\eps}$ is the identity (since $\zeta_L = \zeta_\sigma$). Furthermore, the flow map is continuous for each $\eps > 0$, implying that $\varphi_{S,\eps}$ forms a homotopy of isomorphisms. These facts imply (see eg. Atiyah \cite{atiyah}) that 
\begin{align} \label{eq:familyequivalence}
\mathcal{E}_{\eps}(K) &\cong \mathcal{E}(S,\eps)).
\end{align}

Note that we have fixed $\eps > 0$ up to now. Let us now fix $S^{a,-}_{\varepsilon}$ and $S^{a,+}_{\varepsilon}$  and take $\eps>0$ sufficiently small that Theorem \ref{thm:approx} holds. By Theorem \ref{thm:approx} (a) and (c), 
\begin{align*}
\lim_{\eps \to 0} \hat{E}(T,\lambda,\eps) &= \hat{E}_0(T,\lambda)
\end{align*}
uniformly for $(T,\lambda) \in b_-(S^{a,-}_{\varepsilon}) \cup b_+(S^{a,+}_{\varepsilon})$. Observe here that the information about the jump map is already encoded in the solution $\hat{E}_0(T,\lambda)$, as described in Theorem \ref{thm:approx}(c). Let us now determine the singular limit of the gluing map $\varphi_{S,\eps}$. \\

Let $\zeta_{L}$ and $\zeta_{R}$ correspond (as described earlier) to $T = S^{a,-}_{\varepsilon}$ and $T=S^{a,+}_{\varepsilon}$, respectively, and choose a solution $Z_0(\zeta,\lambda)$ of the reduced problem so that 
\begin{align*}
\hat{E}_{0} &= \reallywidehat{\iota_0(\zeta,\lambda) Z_0(\zeta,\lambda)}.
\end{align*}

It follows from Theorem \ref{thm:approx}(a) that 
\begin{align}
\varphi_{S,0}[\iota_0 (S^{a,-}_{\varepsilon},\lambda) Z_0(S^{a,-}_{\varepsilon},\lambda)] &= \iota_0(\zeta_R,\lambda) Z_0(\zeta_R,\lambda). 
\end{align}
We now define the `singular' limit of the glued hemispheric bundles $\mathcal{E}(S,\eps)$ in the obvious way:
\begin{align}
\mathcal{E}(S,0) &:= \mathcal{E}_0|b_-(S^{a,-}_{\varepsilon}) \cup_{\varphi_{S,0}} \mathcal{E}_0|b_+(S^{a,+}_{\varepsilon}).
\end{align}

We emphasize here that this is a topological limit, and so $\mathcal{E}(S,\eps) \cong \mathcal{E}(S,0)$ directly---this is the advantage of using the augmented unstable bundle. It thus follows from \eqref{eq:familyequivalence} that $\mathcal{E}_{\eps}(K) \cong \mathcal{E}(S,0)$ for each $\eps > 0$ sufficiently small. \\

We have so far held $S^{a,-}_{\varepsilon}$ and $S^{a,+}_{\varepsilon}$ fixed. Noting that the required smallness of $\eps$ from the above argument is independent of these parameters, we are free to send $S^{a,-}_{\varepsilon}$ to 0 from below and $S^{a,+}_{\varepsilon}$ to 0 from above. By Theorem \ref{thm:approx}(b), $\varphi_{S,0}$ approaches $\varphi^H$ as defined in \eqref{eq:clutch}. $\Box$

\section{Counting the slow eigenvalues} \label{sec:slowevans}

We now compute the slow eigenvalues $\lambda \in \Omega$; see Def. \ref{def:sloweigs}. We will adapt the technique presented in Section E of \cite{GJ}; we split the projectivized equation \eqref{eq:projlinsloweqs} for \eqref{eq:linslow} into its real and imaginary parts, with the spatial eigenvalue parameter written as $\lambda = \mu + i\omega$, and we consider separately the cases $\omega \neq 0$ and $\omega = 0$. Afterwards we may restrict our analysis to the real line; using comparison of solutions together with a Gronwall estimate, we demonstrate that there are no positive real eigenvalues. Here we need the monotonicity hypothesis \ref{hyp:monotoneU} to ensure that we can work with the chart defined by $V \neq 0$ for $\lambda \geq 0$.  Generally speaking, the jump condition does not introduce any new complications with respect to these arguments.\\

\begin{thr} \label{thm:sloweigscount}
There are exactly two slow eigenvalues in $\Omega$, given by $\lambda_0 = 0$ and $\lambda_1$, with $\text{Im}(\lambda_1) = 0$ and $R'(0) < \lambda_1 < 0$. Both $\lambda_0$ and $\lambda_1$ are simple. 
\end{thr}

\begin{remk}
As far as spectral stability of the wave for $\eps > 0$ is concerned, it suffices to verify that the eigenvalue with the largest real part in the point spectrum is the simple translational eigenvalue $\lambda_0$. We resort to a numerical calculation using a {\it Riccati-Evans function} to determine  the existence of the secondary eigenvalue $\lambda_1$, but the simplicity of both eigenvalues are verified rigorously.
\end{remk}

{\it Proof that $\lambda_0 = 0$ is the eigenvalue of largest real part:}

We will prove that all slow eigenvalues are real, and furthermore there is no slow eigenvalue $\lambda$ with $\text{Re}(\lambda) > 0$. The asymptotic systems associated with the projectivized linearized slow flow \eqref{eq:projlinslow} admit the fixed points $\hat{u}^-,\,\hat{s}^-$ at $\bar{U} = 0$ and $\hat{u}^+,\hat{s}^+$ at $\bar{U} = 1$, defined as follows:
\begin{equation} \label{eq:projafps}
\begin{aligned}
\hat{u}^-(\lambda)&= \frac{c- \sqrt{c^2 + D(0)(R'(0) - \lambda)}}{2D(0)}\\
\hat{s}^-(\lambda) &= \frac{c+ \sqrt{c^2 + D(0)(R'(0) - \lambda)}}{2D(0)}\\
\hat{u}^+(\lambda) &= \frac{c- \sqrt{c^2 + D(1)(R'(1) - \lambda)}}{2D(1)}\\
\hat{s}^+(\lambda) &= \frac{c+ \sqrt{c^2 + D(1)(R'(1) - \lambda)}}{2D(1)}.
\end{aligned}
\end{equation}

Linear analysis verifies that $\hat{u}^{\pm}$ are attractors and $\hat{s}^{\pm}$ are repellers with respect to the projectivized asymptotic dynamics, corresponding to the asymptotic unstable resp. stable eigendirections. We remind the reader of the geometric characterisation of eigenvalues of \eqref{eq:linslow} in terms of their asymptotic behavior; if $\lambda$ is not a slow eigenvalue, it suffices to show that `the' unique solution $S_0(\zeta,\lambda)$ which tends to $\hat{u}^-(\lambda)$ as $\zeta \to -\infty$ does not tend to $\hat{s}^{+}(\lambda)$ as $\zeta \to +\infty$.\footnote{Strictly speaking, $S_0$ refers to a pair of solutions $S^{a,-}_{\varepsilon}$ and $S^{a,+}_{\varepsilon}$, defined for $\zeta \leq 0$ and $\zeta \geq 0$, which are uniquely defined by two constraints: the aforementioned asymptotic constraint at $\zeta = -\infty$, and the jump condition defined using \eqref{eq:projectivejump}.} \\

We now write \eqref{eq:projlinslow} in terms of its real and imaginary parts. Specifically, writing $S = X + iY$ and $\lambda = \mu+ i \omega$, we have
\begin{equation} \label{eq:XYproj}
\begin{aligned}
\dot{X} &= R'(\bar{U}) - \mu - cX + (X^2-Y^2)D(\bar{U})\\
\dot{Y} &= -\omega - cY + 2XY D(\bar{U}).
\end{aligned}
\end{equation}
We consider the two subcases $\omega \neq 0$ and $\omega = 0$.  \\

{\it The subcase $\omega \neq 0$.} Let us focus on the case $\omega > 0$; the case $\omega < 0$ is similar. We have $\text{Im}\,(\hat{u}^{\pm}(\lambda)) <0$ and $\text{Im}\,(\hat{s}^{\pm}(\lambda))>0$ for each $\lambda \in \Omega$, by applying the inequality $R'(1) < R'(0)<\text{Re}\,(\lambda)$ on this set to the expressions \eqref{eq:projafps} for the asymptotic fixed points. On the other hand, the half-plane $\{Y \leq 0\}$ is forward invariant since $\dot{Y} =-\omega < 0$ along $\{Y = 0\}$. Furthermore, the jump condition is monotone decreasing in $Y$, i.e. $S^{a,+}_{\varepsilon}(0,\lambda)_Y = (J_{\lambda}(S^{a,-}_{\varepsilon}(0,\lambda)))_Y \leq S^{a,-}_{\varepsilon}(0,\lambda)_Y$. Thus, if $S(\zeta,\lambda) = X(\zeta,\lambda) + iY(\zeta,\lambda)$ remains on the chart specified by $V \neq 0$, then it is impossible that $S(\zeta,\lambda) \to s^{+}(\lambda)$ as $\zeta \to +\infty$.\\

It could happen that the solution $S_0(\zeta,\lambda)$ leaves the chart by blowing up in the $Y\to -\infty$ direction, emerging `on the other side' of the half-plane $\{Y >  0\}$ from $Y = +\infty$ to make a connection to $s^{+}(\lambda)$. Let us show that any such blow-up leads to a contradiction. Following the style of Lemma 6.6 in \cite{GJ}, we observe that if the blow-up happens at $\zeta_0 = +\infty$, then $S_0(\zeta,\lambda)$ remains inside the closure of the image of $\{Y \leq 0\}$ on the Riemann sphere, and $s^{+}(\lambda)$ is bounded away from this closed set; hence, we suppose that $\zeta_0 \in \mathbb{R}$, and that $\zeta_0$ is the smallest real number so that $S_0(\zeta,\lambda)$ becomes unbounded as $\zeta \to \zeta_0$ but for which $S_0(\zeta,\lambda)$ remains finite within the lower half-plane for all $\zeta < \zeta_0$. Then $T_0(\zeta,\lambda) :=S_0(\zeta,\lambda)^{-1}$ remains well-defined for all $\zeta$ sufficiently close to $\zeta_0$. This motivates the corresponding change in chart
\begin{align} \label{eq:chartchange}
s= \frac{X}{X^2 + Y^2}, &\hspace{1cm} t = \frac{-Y}{X^2+Y^2}.
\end{align}

The corresponding dynamical system \eqref{eq:XYproj} expressed in the new chart is
\begin{equation} \label{eq:stproj}
\begin{aligned}
\dot{s} &= - D(\bar{U}) + cs - 2\omega st - (R'(\bar{U})-\mu)(s^2-t^2)\\
\dot{t} &= ct - 2st(R'(\bar{U})-\mu) + \omega (s^2-t^2),
\end{aligned}
\end{equation}

with the solution written as $T_0(\zeta) = s(\zeta)+it(\zeta)$. If $Y\to -\infty$, then both $s$ and $t$ tend to 0 as $\zeta \to \zeta_0$; furthermore, $t(\zeta) >0$ for values of $\zeta < \zeta_0$ sufficiently close to the blow-up value according to \eqref{eq:chartchange}. The contradiction will arise by considering the behavior of $t(\zeta)$ for $\zeta$ near $\zeta_0$ using Taylor expansions. In the following argument, we suppose that $\zeta_0 \neq 0$, i.e. we are not exactly at the jump and all two-sided limits of the relevant functions exist; we discuss the case $\zeta_0 = 0$ later. By Taylor expanding $s(\zeta)$ around $\zeta = \zeta_0$ and using the first equation in \eqref{eq:stproj}, we find that
\begin{align}
s(\zeta) &= -D(\bar{U}(\zeta_0))(\zeta-\zeta_0) + \bigO(\zeta-\zeta_0)^2.
\end{align}
Using the second equation in \eqref{eq:stproj}, we note that $\dot{t} = \ddot{t} = 0$ at $\zeta = \zeta_0$; however, we have $\dddot{t}(\zeta_0) = 2\omega D(\bar{U}(\zeta_0))^2$, and so the Taylor expansion of $t(\zeta)$ gives
\begin{align}
t(\zeta) &= \frac{\omega D(\bar{U}(\zeta_0))^2}{3}(\zeta-\zeta_0)^3 + \bigO(\zeta-\zeta_0)^4.
\end{align}
But this implies that $t(\zeta) <0$ for $\zeta < \zeta_0$ and $\zeta$ sufficiently close to $\zeta_0$, which produces a contradiction. \\

If the blow-up occurs exactly at the point of discontinuity $\zeta_0 = 0$, the argument above survives by constructing continuous extensions of $D(\bar{U}(\zeta))$ and $R'(\bar{U}(\zeta))$, since their right limits still exist. Altogether, we have shown that there are no slow eigenvalues with nonzero imaginary part. \\

{\it The subcase $\omega = 0$.} The (un)stable eigenvectors associated with the asymptotic subspaces are real, and by \eqref{eq:XYproj} the subset $\{Y=0\}$ is invariant when $\omega = 0$; hence, we are able to focus on the real one-dimensional problem given by the first equation in \eqref{eq:XYproj}
\begin{align} \label{eq:XYproj1}
\dot{X} &= f(X,\bar{U};\mu):= R'(\bar{U}) - \mu - cX + X^2 D(\bar{U}).
\end{align}
By Remark \ref{rem:slow0}, $\mu = 0$ is a slow eigenvalue.  Now suppose $ \mu > 0$. Noting that $\del f/\del \mu < 0$ and that $\del \hat{u}^-/\del \mu < 0$ from  \eqref{eq:projafps}, any solution $X(\zeta,\lambda)$ which tends to $\hat{u}^-(\lambda)$ as $\zeta \to -\infty$ satisfies the comparison
\begin{align} \label{eq:comparison}
X(\zeta,\mu') \leq X(\zeta,\mu) 
\end{align}
for all $\zeta \in \mathbb{R}$ whenever $0 < \mu \leq \mu'$ (note that this comparison is preserved across the jump). Each such solution for any given $\mu > 0$ has a common upper bound, given by the solution for $\lambda = \mu = 0$, which does not blow up by Hypothesis \ref{hyp:monotoneU}. We also have $\del \hat{s}^+/\del \mu > 0$  for $\mu > 0$ from \eqref{eq:projafps}; thus, $\hat{s}^+(\mu)$ lies on the other side of the connection for $\mu = 0$ (relative to the solution $X(\zeta,\mu)$) . As long as we can show that the solution for $\mu > 0$ does not blow up (potentially allowing $X(\zeta,\mu)$ to connect to $s^+$ by `looping around' the Riemann sphere), it is impossible for the solution $X(\zeta,\lambda)$ to approach $\hat{s}^+(\mu)$ as $\zeta \to \infty$, and we will be done.\\

 According to the comparison to the common upper bound stated above, the solution $X(\zeta,\lambda)$ can only possibly blow up in the direction $X \to -\infty$. This does not happen since for each $\mu > 0$, a specific lower bound is given by Gronwall's inequality and a comparison to the system $\dot{X} = \min_{0 \leq \bar{U} \leq 1}R'(\bar{U}) - \mu - cX$, whose solutions do not blow up. \hfill $\Box$\\

{\it The existence of $\lambda_1$.} We now verify the remaining statements in Theorem \ref{thm:sloweigscount}. It is a straightforward task to numerically integrate the one-dimensional projection of the system \eqref{eq:XYproj} along the invariant subspace of real solutions $\{Y = 0\}$. The solution can blow up by leaving the chart given by $S = P/V$, $V \neq 0$. If this happens, the system \eqref{eq:stproj} can be numerically integrated without issue on the chart $T = V/P$, $P \neq 0$, around a small interval surrounding the pole, and then we can return to the original chart.\\

 \begin{figure}[t!] 
\centering
(a)\includegraphics[width=0.37\textwidth]{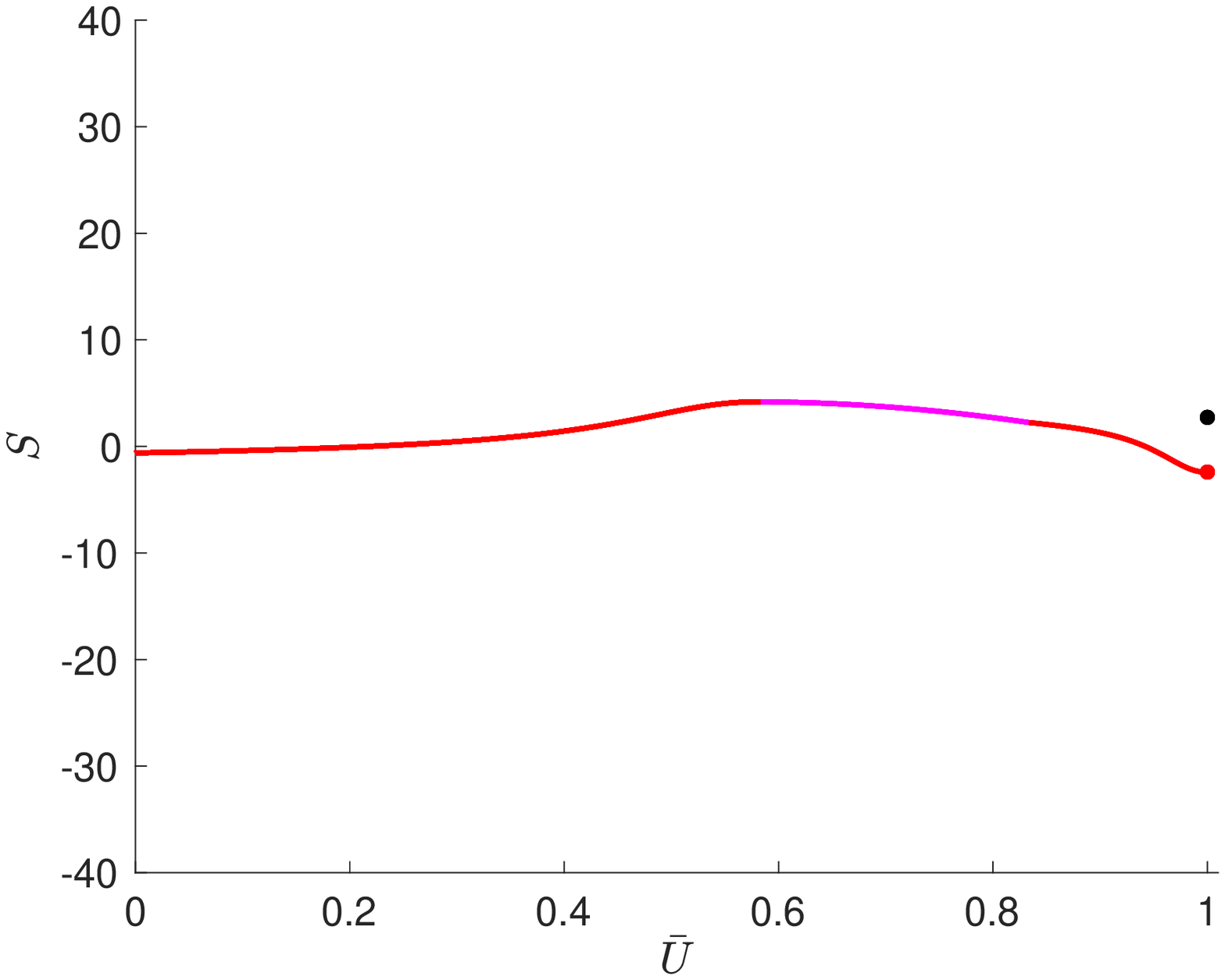}
(b)\includegraphics[width=0.37\textwidth]{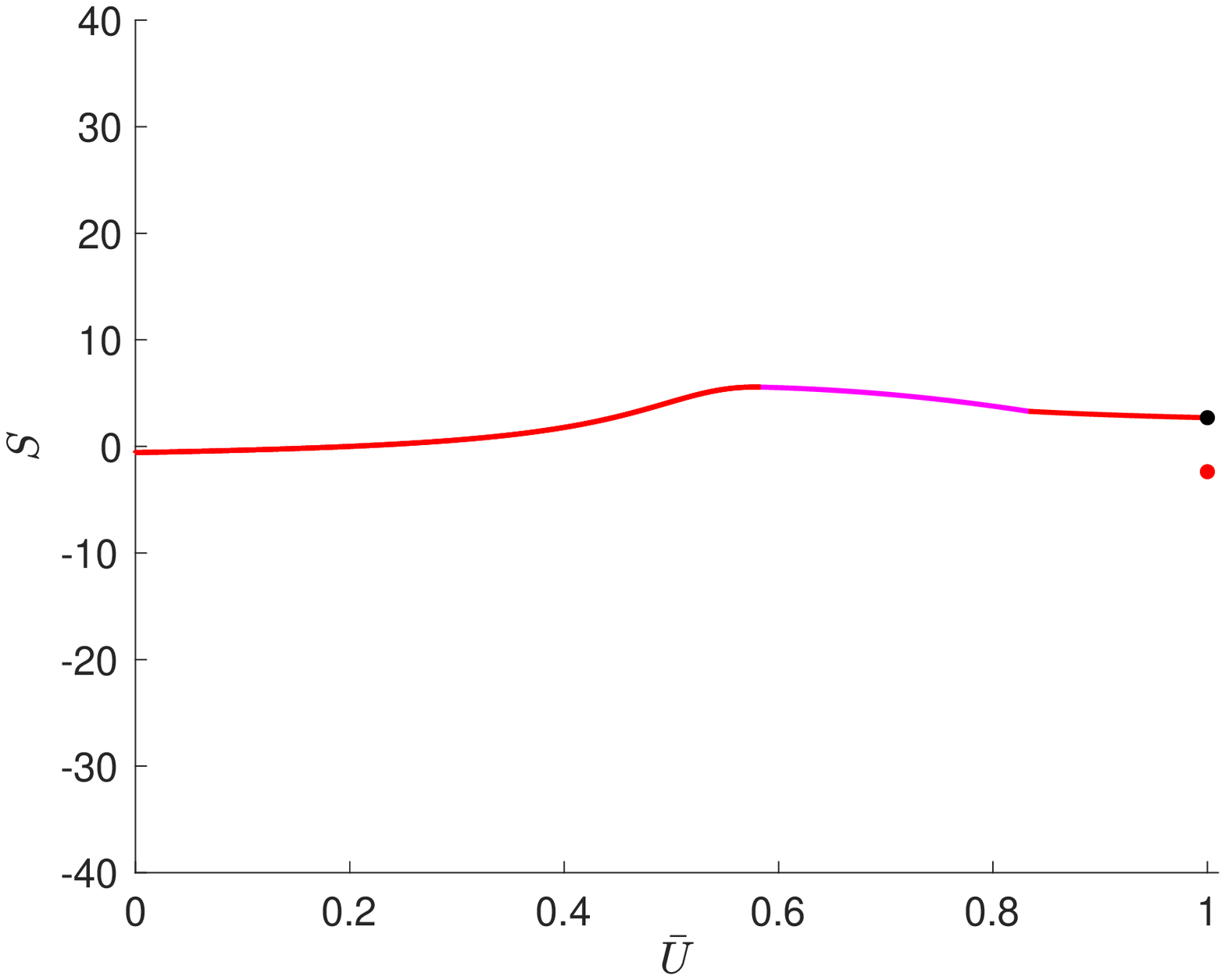}\\
(c)\includegraphics[width=0.37\textwidth]{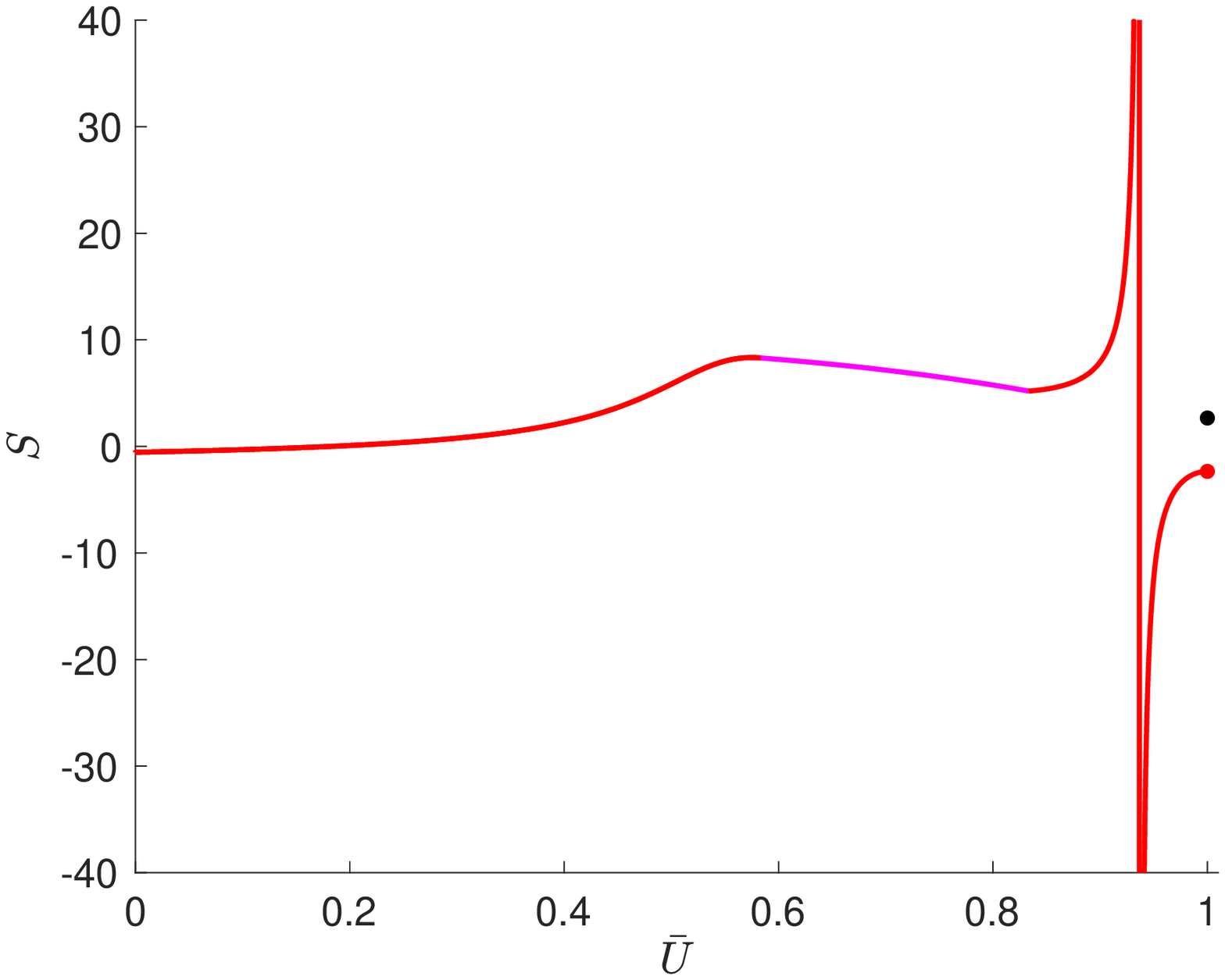}
(d)\includegraphics[width=0.37\textwidth]{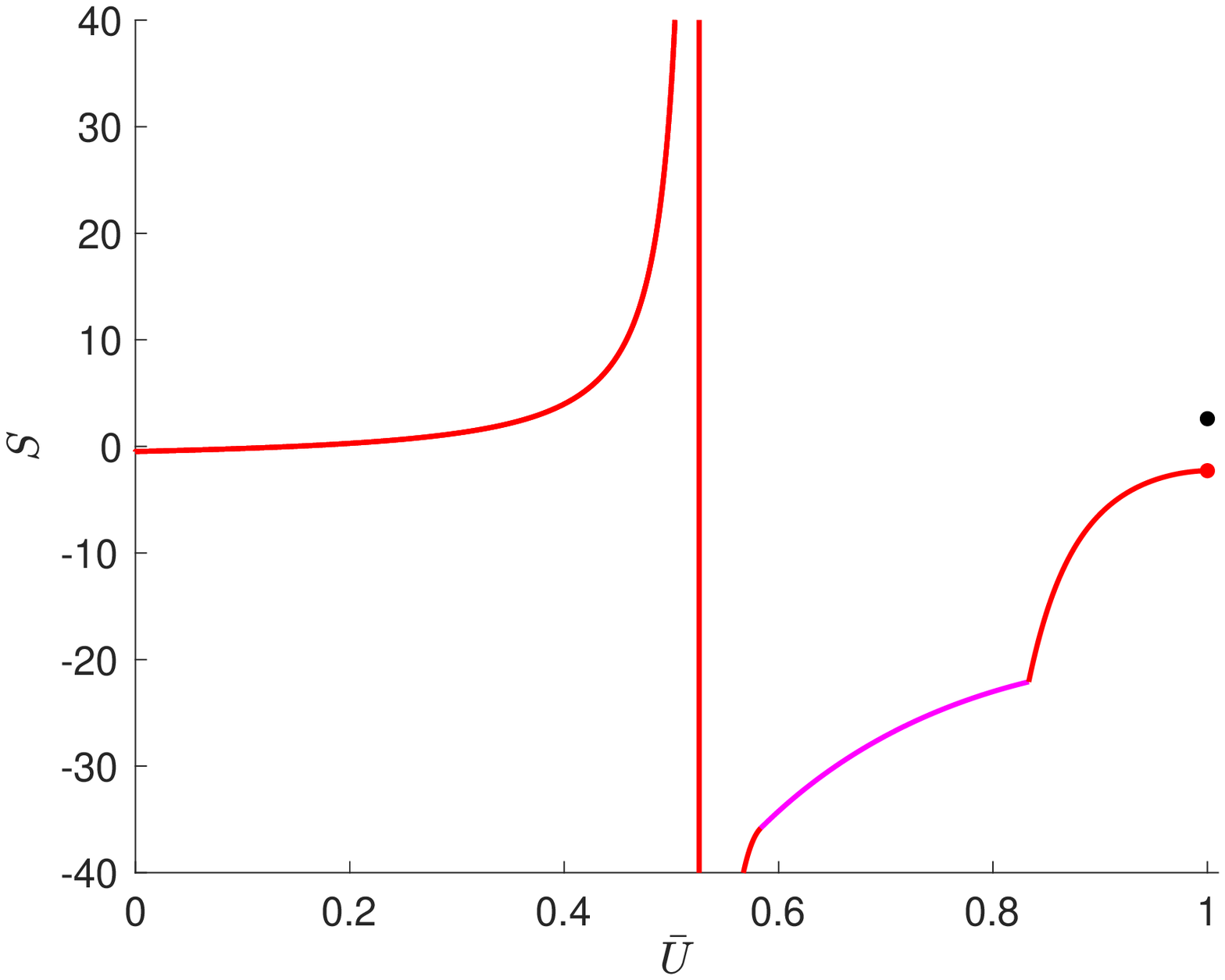}\\
(e)\includegraphics[width=0.37\textwidth]{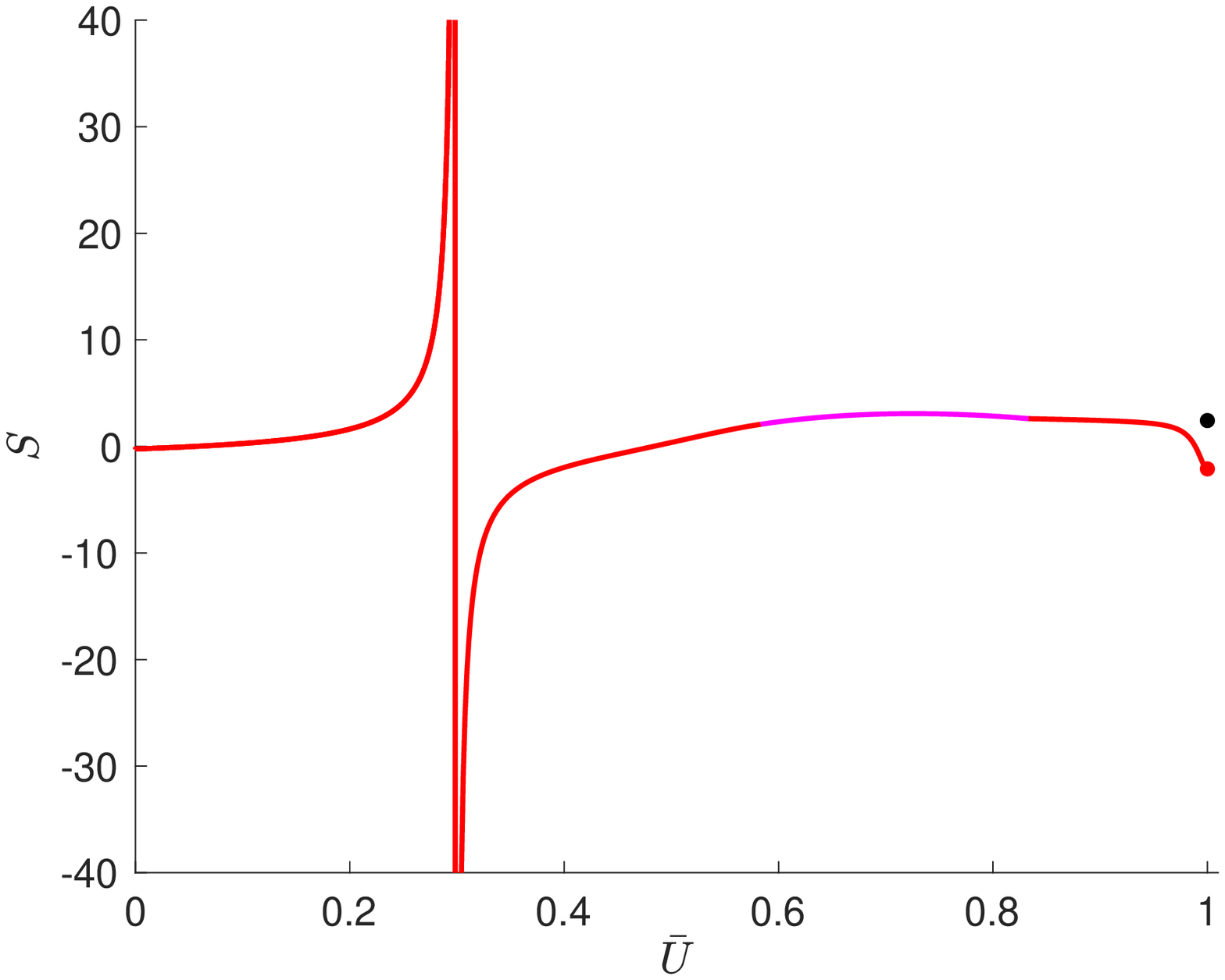}
(f)\includegraphics[width=0.37\textwidth]{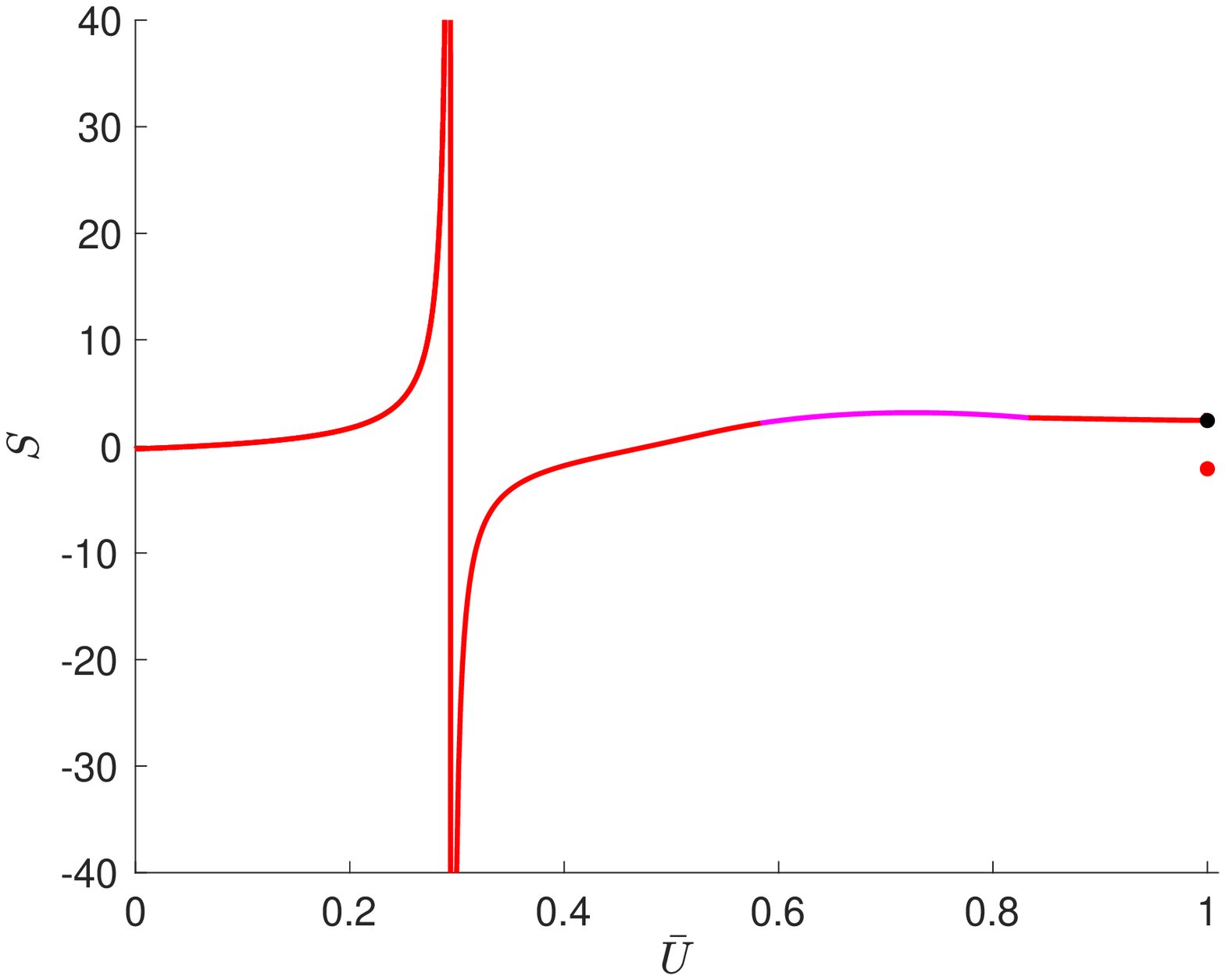}\\
(g)\includegraphics[width=0.37\textwidth]{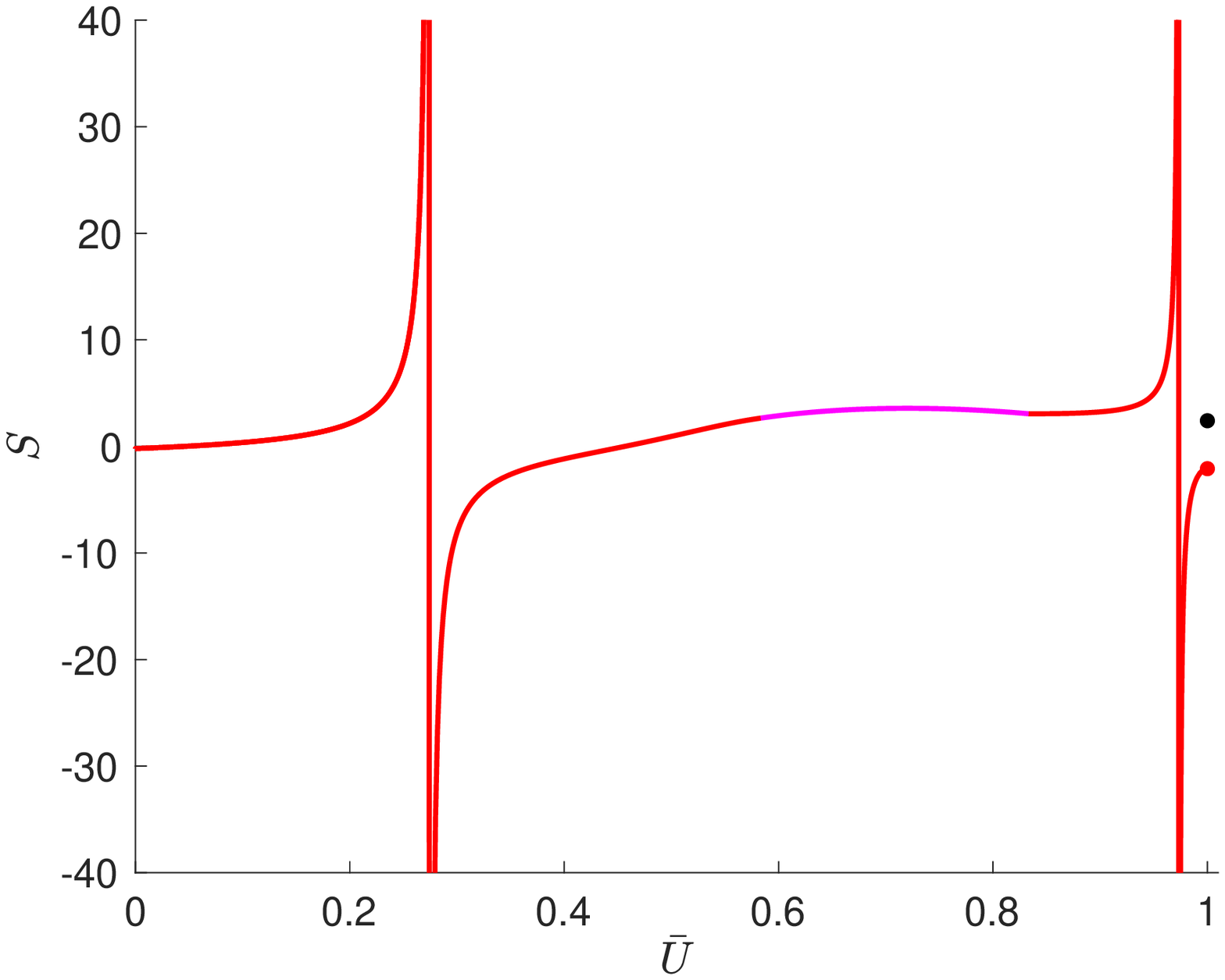}
\caption{Numerical integration of \eqref{eq:XYproj} along $\{Y = 0\}$ for (a) $\lambda = 0.1$, (b) $\lambda = \lambda_0 =  0$, (c) $\lambda = -0.1$, (d) $\lambda = -0.3$, (e) $\lambda = -0.8$, (f) $\lambda =  \lambda_1 \approx -0.80925$, (g) $\lambda = -0.85$. Red resp. black dots: (projectivisation of) the unstable resp. stable eigenvectors of the saddle-point at $\bar{U} = 1$. Magenta segments: concatenated fast jumps defined by the map \eqref{eq:projectivejump}.}
\label{fig:unfoldingslow}
\end{figure}

 These numerics are depicted in Fig. \ref{fig:unfoldingslow} for seven values of $\lambda$ within the range $[-0.85,0.1]$. As shown in Fig. \ref{fig:unfoldingslow}(a)--(c) and (e)--(g), a connection is formed between the unstable eigendirection of the saddle-point at $\bar{U} = 0$ and the stable eigendirection of the saddle-point at $\bar{U} = 1$ for the distinguished values $\lambda = \lambda_0 = 0$ and $\lambda_1 \approx -0.80925$. Simultaneously, we observe the winds that are necessarily generated upon crossing these eigenvalues. These winds can be continued in the parameter $\lambda$---such coordinate singularities can of course be characterised as zeroes of the flow with respect to a suitably chosen chart. They propagate to the left and are preserved across the jump, as is shown in Fig. \ref{fig:unfoldingslow}(c)--(d). \\

{\it The simplicity of the eigenvalues in the reduced problem.}   We have now determined the existence of eigenfunctions as solutions to the Riccati formulation \eqref{eq:XYproj1} of the linearised problem \eqref{eq:linslow}. Recall that the jump conditions for $X$ are given in Def. \ref{def:sloweigs}, and boundary conditions for $X$ given in \eqref{eq:projafps}, i.e.

$$
\lim_{\zeta \to \pm \infty} X = \hat{u}^{-}(\mu), \hat{s}^{+}(\mu). 
$$

In terms of the original reduced problem \eqref{eq:linslow},

this means that we have a value $\lambda = \mu$, and a solution to \eqref{eq:linslow} satisfying the boundary conditions
$$
\lim_{\zeta \to \pm \infty} P, V = 0, 
$$
and appropriate jump conditions. 

Differentiating the equation for $V$ with respect to $\zeta$, and substituting in for $\dot{P}$, and multiplying through by 
$$
g(\zeta) : = \exp\left( - \int^{\zeta} (\log(D(\bar{U})))_{s} + c d s  \right) = \frac{e^{-c\zeta}}{D(\bar{U})}
$$
we are led to the `Sturm-Liouville' form of the reduced problem on the slow manifold: 

\begin{align}\label{eq:SLform}
e^{c\zeta}\frac{d}{d\zeta}\left( \frac{e^{-c \zeta}}{D(\bar{U})} \frac{d}{d\zeta} V \right) +\tilde{Q}(\zeta) V = \mu V. 
\end{align} 
where $\tilde{Q}(\zeta)$ is defined as: 
$$ 
\tilde{Q}(\zeta) := R'(\bar{U}) + \frac{c D'(\bar{U}) \bar{U}_{\zeta}}{D(\bar{U})^{2}} 
$$

Letting $L$ be the symmetric linear operator : 
$$
L[y] :=  e^{c\zeta}\frac{d}{d\zeta}\left( \frac{e^{-c \zeta}}{D(\bar{U})} \frac{d y}{d\zeta} \right) + \tilde{Q}(\zeta) y
$$
we are interested in establishing the simplicity of the eigenvalues of such an operator where $U(\zeta)$ is given as the solution on the slow manifold with the appropriate boundary conditions, and jump at $\zeta =0$. 
\begin{remk}
As $D(\bar{U})$ is discontinuous at $\zeta = 0$, the domain of the linear operator $L$ is a function space that will necessarily incorporate a jump condition which is compatible with the one from Def \ref{def:sloweigs}. Further, the weighting function $e^{c \zeta}$ will define the inner product space on which $L$ is symmetric and well-defined. Here we only aim to establish  simplicity of any eigenvalues/eigenfunctions, which {\em a priori} exist.
\end{remk}

First we note that because $L$ is symmetric on eigenfunctions, the Fredholm alternative means that the geometric multiplicity of any eigenvalues must be the same as the algebraic multiplicity, and in particular if $(L-\mu)[y_{1}]=0$, there can be no non-zero solutions to $(L-\mu)[y] =y_{1}$.

Now, suppose that $v_{1}$ and $v_{2}$ are eigenfunctions of $L$ with the same eigenvalue $\mu$, so $L[v_{j}] = \mu v_{j}$. On the one hand, we have
\begin{align}
v_{1}L[v_{2}] - v_{2}L[v_{1}] = \mu(v_{1}v_{2}-v_{2}v_{1}) = 0,
\end{align}
while on the other we have 
\begin{align}
v_{1}L[v_{2}] - v_{2}L[v_{1}] & = e^{c\zeta}\frac{d}{d\zeta}\left(g(\zeta)(v_{1}\frac{d v_{2}}{d\zeta}-\frac{d v_{1}}{d\zeta} v_{2})\right).
\end{align}

Thus, the quantity 
\begin{align}
g(\zeta)(v_{1}\frac{d v_{2}}{d\zeta}-\frac{d v_{1}}{d\zeta} v_{2}) =: \frac{e^{-c\zeta}}{D(\bar{U})} W(v_{1},v_{2})
\end{align}
is a constant. Evaluating at $\zeta = + \infty$ we have 
$$
\lim_{\zeta \to +\infty}\frac{e^{-c\zeta}}{D(\bar{U})} W(v_{1},v_{2}) = 0.
$$ 
We observe that
$$
 v_{1,2}(\zeta) \sim e^{\frac{1}{2}\left(c + \sqrt{c^{2} + 4 D(0)(\mu -R'(0))}\right)\zeta}
$$
as $\zeta < 0$ grows very large in magnitude, and so in particular 
$$
\lim_{\zeta \to -\infty} \frac{e^{-c\zeta}}{D(\bar{U})} W(v_{1},v_{2}) = 0
$$
when $\mu > R'(0)$. Thus $\frac{e^{-c\zeta}}{D(\bar{U})} W(v_{1},v_{2}) \equiv 0$ on both sides of the jump and hence we have that $W(v_{1},v_{2}) \equiv 0$. Since the Wronskian of the two eigenfunctions is identically zero, they must be linearly dependent. We conclude that the slow eigenvalues $\lambda_{0}$ and $\lambda_{1}$ are both simple.  $\Box$

~\\

  \begin{figure}[t!] 
\centering
\includegraphics[width=0.98\textwidth]{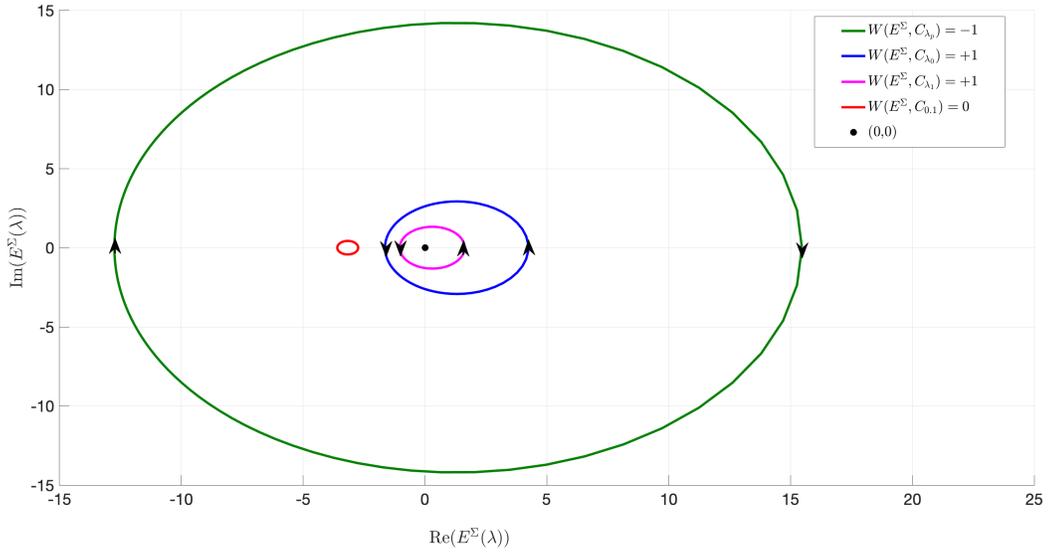}
\caption{Evaluation of the Riccati-Evans function \eqref{eq:riccati} along circles $C_z$ of radius 0.03  centered at the point $z \in \mathbb{C}$. The winding numbers along each contour are given in the legend and are also depicted by the arrowheads. Note that $\lambda_{0}$ and $\lambda_1$ denote the two eigenvalues and $\lambda_p =-0.08$ is the approximate location of the pole of $E^{\Sigma}(\lambda)$.}
\label{fig:riccatievans}
\end{figure}

We supplement our proofs with a numerical demonstration of the existence of the simple eigenvalues $\lambda_0,\,\lambda_1$, by means of a {\it Riccati-Evans} function. Define the cross section $\Sigma = \{\bar{U} = 0.95\}.$ The Riccati-Evans function associated with the Riccati equation \eqref{eq:XYproj} on the coordinate chart $\{V \neq 0\}$ is 
\begin{align} \label{eq:riccati}
E^{\Sigma}(\lambda) &:= s^{\Sigma}_1(\lambda)-u^{\Sigma}_0(\lambda),
\end{align}
where $s^{\Sigma}_1(\lambda)$ denotes the (first) intersection of the unique nontrivial solution $s_1(\zeta,\lambda)$ of \eqref{eq:XYproj} which converges asymptotically to the stable eigenvector of the  saddle point at $\bar{U} = 1$ as $\zeta\to+\infty$, and similarly for $u^{\Sigma}_0(\lambda)$, which connects to the unstable eigendirection of the saddle point as $\zeta \to -\infty$. Compare our definition to the general construction in \cite{harley}.\\

 We highlight a few key points about our function \eqref{eq:riccati}. The function is meromorphic (hence satisfying the argument principle), and it vanishes on eigenvalues $\lambda$. Its zeroes are intrinsic, depending neither on the choice of section nor on the choice of chart. On the other hand, as we can anticipate from the dynamics depicted in Fig. \ref{fig:unfoldingslow}, the poles arise for values of $\lambda$ where the solution leaves the chart by winding. In other words, the poles are an artifact of the choice of coordinate chart as well as the choice of section, and can usually  be moved (or even removed entirely) by a judicious selection of the chart and the section (with our choice of $\Sigma$, we locate a pole near $\lambda_p = -0.08$). See \cite{harley} for a general exposition.\\
 
  In this analysis we find it instructive to work with the `naive' chart $S = P/V$, $V \neq 0$, for all values of $\lambda$---the dimensionality of the slow eigenvalue problem is low enough that we can easily demonstrate the utility of the argument principle. As shown in Fig. \ref{fig:riccatievans}, the multiplicity of each eigenvalue can be readily computed by evaluating the corresponding winding number $W(E^{\Sigma},C)$ along simple closed contours $C$. Contours which surround poles, corresponding to blow-up of solutions off the chosen coordinate chart, will map under $E$ to contours with a negative winding number (corresponding to a clockwise orientation). Here we demonstrate that the poles are also simple. This completes our numerical verification of Theorem \ref{thm:sloweigscount}.

\section{Concluding remarks} \label{sec:conclusion}

We have given a complete characterisation of the spectral stability problem for shock-fronted travelling waves of the regularized system \eqref{eq:master}. But other types of high-order regularization terms can be applied to the underlying system that exhibits shock solutions; it can further be shown that these nonequivalent regularizations pick out distinct one-parameter families of shock-fronted travelling waves, which limit to singular solutions satisfying different rules. \\

Consider the following system which applies to both a viscous relaxation and a fourth-order {\it nonlocal} regularization, where the new parameter $a \geq 0$ characterises the relative weighting of the two regularizations. 

\begin{align} \label{eq:homotopy}
\frac{\del \bar{U}}{\del t} &= \frac{\del}{\del x} \left( D(\bar{U}) \frac{\del \bar{U}}{\del x} \right) + R(\bar{U}) + \eps a \frac{\del^3 \bar{U}}{\del x^2 \del t} - \eps^2 \frac{\del^4 \bar{U}}{ \del x^4}.
\end{align}

Setting $a = 0$ recovers a regularized system with `purely' nonlocal regularization. In this case, it can be shown that the corresponding one-parameter family of travelling waves solutions has a singular limit as a shock-fronted travelling wave which satisfies a so-called `equal area rule'-- the shock connection across disjoint branches of the potential function $F(\bar{U})$ is selected so that the area bounded above and below the shock height are exactly balanced.  Recent work \cite{proceedings22} has shown that families of shock-fronted waves persist robustly for $a > 0$, with each such family satisfying a generalized equal area rule in the singular limit. Furthermore, there is a finite value $a = a_V >0$ for which the shock connection is formed at the fold, again resulting in viscous-type shocks. The problem in this paper can be thought of as the regularization in the (scaled) limit as $a \to \infty$.  \\

\begin{figure}[h]
\begin{center}
\includegraphics[scale=0.4]{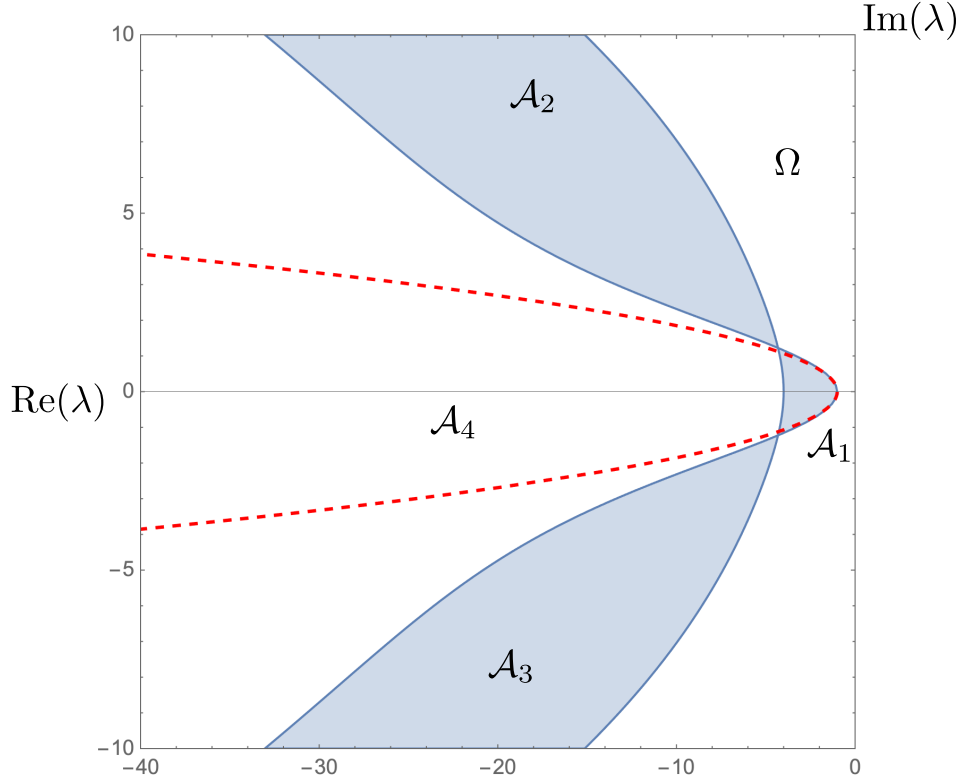}
\caption{Plot of the essential spectrum (shaded regions) of the of the linearised operator L in \eqref{eq:lin4}, which shows $\Omega$ and the $\cA_{i}$ partitioning the complex plane via the Fredholm borders (continuous spectrum). The shape of the essential spectrum means that the operator $L$, for non-zero $\ve$ is sectorial. The dashed line (red online) is one of the Fredholm borders  coming from the linearisation of the  reduced problem about the (constant) steady-state $\bar{U}^{-}$ c.f. the right figure in Fig \ref{fig:essential}. For the figure, $\ve =0.1$ and $c =a = 1$, with $D$ and $R$ as in \eqref{eq:diffusionterm} and \eqref{eq:reactionterm}. }\label{fig:4dessential}
\end{center}
\end{figure}

Thus, let $ a>0 $ be given and assume the existence of such a family of travelling waves for \eqref{eq:homotopy}. As before, we can find such waves as standing wave solutions $\bar{U}(z,t) \equiv \bar{U}(z)$ to:
\begin{align} \label{eq:travemix}
(\bar{U} - \ve a \bar{U}_{zz})_{t} = - \ve^{2}\bar{U}_{zzzz} -\ve a c \bar{U}_{zzz} + (F(\bar{U}))_{zz} +c\bar{U}_{z} + R(\bar{U}).
\end{align}
We now discuss how the stability problem changes in this more general case. Linearising about a standing wave solution $\bar{U}(z)$ to \eqref{eq:travemix} leads to the eigenvalue problem 
\begin{align}\label{eq:lin4}
\lambda p - \ve a \lambda p_{zz} = - \ve^{2}p_{zzzz} - \ve a c p_{zzz} + (D(\bar{U})p)_{zz}+ cp_{z}+R'(\bar{U}) p. 
\end{align}
Defining the variables 
\begin{align}
s & := \ve^{2}p_{zzz}+ \ve a c p_{zz} - ((D(U)+ \ve a \lambda) p)_{z}- cp \\ 
r & := \ve^{2}p_{zz} + \ve a c p_{z}- (D(U)+\ve a \lambda)p \\
q & = \ve p_{z}+ a c p,
\end{align}
we can write the closed system
\begin{align}
\begin{pmatrix} \ve p \\  \ve q \\ r \\ s 
\end{pmatrix}_{z} & =
\begin{pmatrix}
- ac & 1 & 0 & 0 \\ 
D(\bar{U}) + \ve a \lambda & 0 & 1 & 0 \\ 
c & 0 & 0 & 1 \\ 
R'(\bar{U}) - \lambda & 0 & 0 & 0 
\end{pmatrix} 
\begin{pmatrix} p \\ q \\ r \\ s 
\end{pmatrix},
\end{align}
with two fast and two slow variables. The calculations for the essential (and absolute) spectrum follow the same ideas as in Section \ref{sec:spectrum}. Presuming the existence of a wave $\bar{U}$ which exponentially approaches its end states $\bupm$, we have the dispersion relations are the pair of parametrised equations given by $k \in \R$: 
$$
\lambda_{\pm} = \frac{-\ve^{2} k^{4} - D(\bupm) k^{2} +R'(\bupm)}{1+ a \ve k^{2}} +i c k. 
$$

These are a pair of curves in the complex plane which are always opening leftward and intersecting the real axis at the points $(R'(\bupm),0)$ in the left half plane. The dispersion relations again form the Fredholm borders and make up the continuous spectrum. As before, $\C$  is partitioned into five regions: the region $\Omega$ which contains the right half plane, together with $\cA_{j}$ for $j = 1,2,3,4$. The essential spectrum in this case can again be determined by considering the signatures of the asymptotic matrices, and can be seen to consist of the regions $\cA_{1,2,3}$, which is the region ``between'' the Fredholm borders. We have

\begin{center}
\begin{tabular}{c|| c | c }
Region & $\sgn(A^{-})$ & $\sgn(A^{+})$ \\ 
\hline
$\Omega$ & $(-, -, ++)$ & $(-,-,+,+)$ \\ 
$\cA_{1}$ & $(-, -,-, +)$ & $(-,-,+,+)$ \\
$\cA_{2}$ & $(-, -, + +)$ & $(-,-,-,+)$ \\
$\cA_{3}$ & $(-, -, +, +)$ & $(-,-,-,+)$ \\
$\cA_{4}$ & $(-, -,-, +)$ & $(-,-,-,+)$ \\
\end{tabular}\label{tab:4dsigs}
\end{center}

Again, the Fredholm borders will be close to the Fredholm borders of the linearised reduced problem, linearised about the steady-states corresponding to the end states of the unperturbed problem. Eventually, as with the third order perturbation problem, the higher order modes in the dispersion relations will dominate, and the essential spectrum will diverge from the appropriate continuous spectrum of the unperturbed problem (see Fig: \ref{fig:4dessential}). In contrast to the case of `pure' viscous relaxation, however, the Fredholm borders do not asymptote to vertical lines in the complex plane and so the operator {\em is} sectorial, and spectral stability in this case would indeed imply linearised stability of the perturbation problem. \\

Having shown that the essential spectrum remains remarkably well-behaved in this more general context, let us now highlight key differences (and difficulties) in calculating the point spectrum. Here, the problem is now four-dimensional, and it can be verified with direct calculation that the two-dimensional slow manifolds on which lie saddle-type equilibria at $\bar{U} = 0$ and $\bar{U} = 1$ are now saddle-type (i.e. there is now both a fast stable and fast unstable direction). The augmented unstable bundle construction for $\eps > 0 $ now defines a complex 2-plane bundle over the sphere. \\

It is natural to ask whether the present problem is amenable to the splitting techniques introduced in \cite{GJ}: namely, is it possible to decompose this 2-plane bundle into `fast unstable' and `slow unstable' line bundles which are controlled by separated reduced eigenvalue problems?\\

 We conjecture that such a separation is not possible: as in the `pure' viscous case, it can be directly calculated that the eigenvalue problem again introduces the eigenvalue parameter $\lambda$ only `weakly,' i.e. through $\mathcal{O}(\eps \lambda)$ terms. Thus, the reduced fast eigenvalue problem again degenerates. As mentioned in the introduction, this poses issues for the construction of a fast elephant trunk over the entire wave; consequently, it does not appear to be feasible with the present techniques to uniformly separate a fast unstable line bundle from a slow unstable one. Furthermore, any such reduced line bundle over the fast layer is governed by the layer flow, which connects an unstable eigendirection to a {\it stable} one in the singular limit (i.e. we are not able to construct an {\it augmented} unstable line bundle for the reduced problem).  \\
 
Nonetheless, we assert that estimates similar to Lemma \ref{lem:exchange} can be brought to bear to regain control of the unstable 2-plane bundle as we track it across the fast shock layer.  It would be of interest to write down sharper (i.e. exponential closeness) estimates to the invariant manifold near the slow subbundle, as discussed in Remark \ref{remk:exchange}. The calculation in Appendix \ref{app:toy} gives a concrete example of such an exponential closeness result in an elementary problem. Such calculations appear to be very cumbersome in the general case. It would also be interesting to situate the construction of such invariant manifolds within the context of general invariant manifold results in nonautonomous dynamical systems (see e.g. \cite{rasmussen}). These are topics of ongoing work.

\section*{Acknowledgement}

Both authors acknowledge support from the Australian Research Council Discovery Project grant DP200102130. The authors thank  Martin Wechselberger at the University of Sydney for several insightful discussions. The first author thanks Martin Rasmussen at Imperial College London for a clarifying discussion about invariant manifolds of nonautonomous dynamical systems. 

\appendix 

\section{Construction of the fast and slow elephant trunks} \label{app:trunks}

{\it Proof of Lemma \ref{lem:fasttrunks}:} We follow the strategy to the proof of Lemma 4.2 in \cite{GJ}. The goal is to verify the following four conditions for the projectivization of the linearized system $\eqref{eq:visclinfast}$:
\begin{itemize}
\item ($(4.6)_{GJ}$ in Sec. IV of \cite{GJ}): There exists $\alpha > 0$ independent of $\eps$ and $\xi \in I$ so that the curve of critical points $\beta_0(\gamma,\eps)$ associated with the frozen family \eqref{eq:arbnonautfreeze} satisfies $$\text{Re}\,\sigma[G_{\beta}(\beta_0(\gamma,\eps),\gamma,\eps)]<-\alpha;$$
\item ($(4.7)_{GJ}$ in Sec. IV of \cite{GJ}): for any $d > 0$ there exists $\eps_0 = \eps_0(d)$ and a nested family of subintervals $I(d) \subset I$ (i.e. $I(d_1) \subset I(d_2)$ for $d_1 < d_2$) with $$\text{sup}_{(\beta,\xi,\eps) \in C} \{ |G|_{\xi} ,|| G_{\beta,\eps} ||, |\beta_{0\xi}|\} < d,$$ where the set $C$ can be taken as $$C = \{(\beta,\xi,\eps): |\beta-\beta_0(\xi,\eps)| < c_0, \xi \in I(d),0 < \eps \leq \eps_0\}$$ for some $c_0>0$;
\item ($(4.8)_{GJ}$ in Sec. IV of \cite{GJ}): we have $$K = \text{sup}_C ||D^2 G(\beta,\xi,\eps)|| < \infty; \text{~and}$$ 
\item ($(4.9)_{GJ}$ in Sec. IV of \cite{GJ}): let $A(\gamma) :=  G_{\beta}(\beta_0(\gamma,\eps),\gamma,\eps)$. Then there exists some invertible matrix $Y(\gamma)$ depending smoothly on $\gamma$, such that the following holds for some $ a> 0$ depending only on $\alpha$ (from condition $(4.6)_{GJ}$): 
\begin{align*}
\text{Re}(A(\gamma)\beta,\beta)_{\gamma} &< -a |\beta|_{\gamma}^2,\\
||Y(\gamma)||_{\gamma} &= 1,\\
||Y_{\gamma}(\gamma)||_{\gamma} &\leq d,
\end{align*}
where we define the inner product $(\beta_1,\beta_2)_{\gamma} := Y(\gamma) \beta_1 \cdot \overline{Y(\gamma) \beta_2}$ and the norm $|\beta|_{\gamma} = (\beta,\beta)^{1/2}_{\gamma}$.
\end{itemize}

Lemma 4.1 in \cite{GJ} then gives the result. By Lemma \ref{lem:asympt}, there exists $\bar{\xi} > 0$ such that the second bullet point in the Lemma above holds for all $|\xi| \geq \bar{\xi}$ and all sufficiently small $\eps > 0$. This verifies the condition $(4.6)_{GJ}$. Furthermore, the branch of equilibria $\beta_0(\xi,\lambda,\eps)$ for the corresponding frozen family is uniformly bounded for all $|\xi| > \bar{\xi}$  and all $\lambda \in K$, with $K \subset \Omega$ a fixed contour. Hence condition $(4.8)_{GJ}$ holds for some constant $K$, with $c_0$ in the definition equal to 1.\\

 Let us check the remaining conditions $(4.7)_{GJ}$ and $(4.9)_{GJ}$. Take the interval $I$ to be either $I_- = \{\xi \in \mathbb{R}: -\infty < \xi < -\bar{\xi}\}$ or $I_+ = \{\xi \in \mathbb{R}: \bar{\xi} < \xi < +\infty\}$. On either interval, the projectivized vector field $G(\beta,\xi,\lambda,\eps)$ depends on $\xi$ through components of the (nonautonomous) linearisation matrix, in particular through the diffusion and reaction terms (i.e. through $\bar{u}$). Hence for each constant $L >0$ there exists  $\eps_0 := \eps_0(L)$, so that for $0 < \eps \leq \eps_0$ we have $$\max \{ |G_{\xi}|, || G_{\xi,\beta}||, |\beta_{0,\xi} |\} < L \max\{|u'(\xi,\eps)|\}.$$

Now let $n(\delta)$ denote a $\delta$-tube around the singular limit of the travelling wave and assume that $I= I_-$. Resetting $\eps_0$ again if necessary, we can assume that the travelling wave $x(\xi,\eps)$ lies within $n(\delta)$ for each $\xi \leq \xi_{L}(\delta)$ . Then there exists a constant $M>0$ so that $|\bar{u}'(\xi,\eps)| \leq M\delta$ for all $\xi \leq \xi_L(\delta)$, since $\bar{u}' = \bar{v} - \Phi(\bar{u})$ is uniformly bounded along the entire wave.  Hence we may take $\delta(d) = d/M$ and define $$I_-(d) = \{\xi: \xi_L(\delta(d))\}.$$ Note that $\delta(d)$ depends smoothly on $d$ and that $\delta(0) = 0$, and furthermore that the definition for $I_-(d)$ provides nested intervals. This verifies condition $(4.7)_{GJ}$.\\

To verify the remaining condition $(4.9)_{GJ}$, it is enough to note that the linearization $A(\xi,\lambda,\eps) = A_0(\xi,\lambda,\eps) + \mathcal{O}(\delta) $ of $G$ (the projectivization of the linear system), evaluated along the strong unstable eigendirection, is a $2\times 2$ matrix which has two $\mathcal{O}(1)$ negative eigenvalues. This matrix  $A(\xi,\lambda,\eps)$, is necessarily negative definite for $I_L(d)$, and hence we can take $Y = A_0^{-1}$ such that $Y^{-1}AY$ produces a matrix with diagonal entries of uniformly negative real part and remaining entries of $\mathcal{O}(\delta)$. This immediately implies the needed inequalities. Applying Lemma 4.1 in \cite{GJ}, we have constructed the necessary elephant trunk $\Omega_-^f$ over $S^{a,-}_{\eps}$.\\

The preceding discussion applies identically for the case $I = I_+$ to produce a fast elephant trunk along $S^{a,+}_{\eps}$.  $\Box$\\

\section{Lemmas for estimates near the slow subbundle} \label{app:estimates}

\begin{lem} \label{lem:slowmanestimate}
Let $\delta,\, a > 0$ be fixed. There exists $\bar{\eps} > 0$ sufficiently small and $T > 0$, both depending only on $\delta$ and $a$, so that for each $0 < \eps \leq \bar{\eps}$ and $|\gamma| \geq a/\eps$, each solution $\hat{y}(s)$ of \eqref{eq:frozen} must satisfy at least one of the following:

\begin{itemize}
\item[(i)] $\hat{y}(0) \in N_{\delta} (\hat{\sigma}_s(\gamma,\lambda,\eps))$ \\
\item[(ii)] $\hat{y}(-T) \in N_{\delta} (\hat{f}_1(\gamma,\lambda,\eps))$. \\
\end{itemize}
\end{lem}

\begin{remk} \label{rem:gamma}
Fix some $\gamma$ as above and let $\hat{y}(s)$ denote a solution to the corresponding member of the family of projectivized frozen systems in \eqref{eq:frozen}, and let $y(s) \in \pi^{-1}(\hat{y}(s))$ be chosen so that $|y(0)|_{\infty} = 1$. The corresponding frozen linearized system in \eqref{eq:frozen} is autonomous and linear, and thus its solution can be written explicitly as a linear combination of eigenvectors
\begin{align} \label{eq:linearsol}
y(s) &= g_1 f_1e^{s\mu_{1,s}}+ g_2 f_2 e^{s\mu_{2,s}}+g_f f_f e^{s\mu_{f}}.
\end{align}
It thus follows from invariance of (generalized) eigenspaces that if (i) in Lemma \ref{lem:slowmanestimate} above does not hold for some fixed $\gamma$, then $g_f \neq 0$, and so there is some $T > 0$ so that (ii) holds. The point of this lemma is to show that for all sufficiently small $\eps$, such a $T > 0$ exists which can be chosen depending only on $\delta$ and $a$.
\end{remk}

{\it Proof of Lemma \ref{lem:slowmanestimate}}: Assume that item (i) does not hold for some $\gamma$ as given in the Lemma, so that there exists $T$ so that item (ii) holds for that value of $\gamma$. Using the normalization of $y(s)$ in Remark \ref{rem:gamma}, we note that $|g_i| \leq 1$. By Lemma \ref{lem:metbound}, there exists some $K > 0$, depending only on the metric $\rho$, so that
\begin{align} \label{eq:gfestimate}
|g_f| \geq K \delta.
\end{align}

From \eqref{eq:linearsol} we may write
\begin{align} \label{eq:factoredlin}
y(-T) &= e^{-\mu_{f}T}(g_f e_f + R(T)),
\end{align}
where
\begin{align*}
R(T) &= g_1 f_1 e^{-T(\mu_{1,s}-\mu_{f})} + g_2 f_2 e^{-T(\mu_{2,s}-\mu_{f})}.
\end{align*}

At this stage we remind the reader that the $\mu_i$ are still $\gamma$- and $\eps$-dependent. We now apply the asymptotic estimate in Lemma \ref{lem:asympt} to find $\bar{\eps} \geq 0$, and $\alpha > 0$ depending only on $\bar{\eps}$  and $a$, so that for each $\eps \in (0,\bar{\eps}]$ and $|\gamma|\geq a/\eps$, we extract the uniform bound $\min\{\text{Re}(\mu_{1,s} - \mu_f),\,\text{Re}(\mu_{2,s} - \mu_f)\} > \alpha$. It thus follows from the triangle inequality and the estimates $|g_i| \leq 1,\,|f_i|=1$ that
\begin{align*}
|R(T)| &< 2e^{-\alpha T}.
\end{align*}

For each $T' > T$, we have $|R(T')| < |R(T)|$. This estimate and \eqref{eq:gfestimate} applied to  the factored form \eqref{eq:factoredlin} then implies that 
\begin{align*}
\hat{y}(-T) \in N_{\delta}(\hat{f}_f(\gamma,\lambda,\eps))
\end{align*}

for $T$ depending only on $\alpha,\,\delta,\,a,\,\text{and }\rho$. $\Box$\\

We now turn to the dynamics on the slow timescale. Our objective here is to strengthen the estimate in Lemma \ref{lem:slowmanestimate} slightly, so that $a(\eps) \equiv a$ can be chosen independently of $\eps > 0$ such that the solution $\hat{E}$ remains uniformly near the slow subbundle for $|\zeta| \geq a$. This will be crucial in making the comparison to the linearized reduced flow defined on $S^{a,-}_{\varepsilon},S^{a,+}_{\varepsilon}$. \\

From the construction in Lemmas \ref{lem:fasttrunks} and \ref{lem:slowtrunks}, there exists $\delta_1>0$ such that for each sufficiently small $\eps$, 
\begin{align*}
N_{\delta_1}(\hat{f}_f(\xi,\lambda,\eps)) \subset \hat{\Omega}_f(\xi) &\text{~for~}\xi \leq -a/\eps \text{~and}\\
N_{\delta_1}(\hat{\sigma}_s(\xi,\lambda,\eps)) \subset \hat{\Omega}_s(\xi) &\text{~for~}\xi \geq a/\eps,\\
\end{align*}

where $\hat{\Omega}_f(\xi)$ and $\hat{\Omega}_s(\xi)$ denote the projectivizations of the slices of the corresponding elephant trunks within $\mathbb{CP}^2 \times \{\xi\}$.

\begin{lem} \label{lem:unifslowest}
Set $\delta >0$ with $\delta < \delta_1$, with $\delta_1$ as above, and fix $a > 0$.  Suppose that $\hat{Y}$ is a solution to the projectivization of the slow linearized equations
\begin{align}
\dot{\hat{Y}} &= \hat{A}(\hat{Y},\zeta,\lambda,\eps)
\end{align}
so that for each $\eps$ sufficiently small, there exists $A(\eps) > a$ with the property that
\begin{align}
\hat{Y}(\zeta,\lambda,\eps) \in N_{\delta}(\Sigma_s(\zeta,\lambda,\eps)). 
\end{align}

Then there exists $\bar{\eps} > 0$ such that for each $0 < \eps < \bar{\eps}$ and $|\zeta| \geq a$ we have
\begin{align} \label{eq:Yest}
\hat{Y}(\zeta,\lambda,\eps) \in  N_{\delta}(\Sigma_s(\zeta,\lambda,\eps)). 
\end{align}
\end{lem}

{\it Proof of Lemma \ref{lem:unifslowest}}:

Our strategy will be to verify the uniform closeness estimate separately over the slow branches $S^{a,-}_{\varepsilon}$ and $S^{a,+}_{\varepsilon}$. The argument for $S^{a,-}_{\varepsilon}$ follows the indirect proof in \cite{GJ} closely. With our modified exchange lemma in hand, the remaining closeness estimate over $S^{a,+}_{\varepsilon}$ is direct. \\

Suppose the lemma were false over $S^{a,-}_{\varepsilon}$; then there would exist sequences $\zeta_n < -a$ and $\eps_n \to 0$ so that \eqref{eq:Yest} fails to hold for all $n$. By passing to subsequences, we can assume that the following sequences converge simultaneously:

\begin{align*}
\zeta_n &\to \bar{\zeta} \text{ ~~~where~~ } a \leq -\bar{\zeta} < \infty,\\
X(\zeta_n,\eps_n) &\to \bar{x}~~~ \in S^{a,-}_{\varepsilon} , \text{ and}\\
\hat{Y}(\zeta_n,\lambda,\eps_n) = \hat{Y}_n &\to \hat{Y}_* ~~\in \mathbb{CP}^2.
\end{align*}

Let us first suppose that $\bar{x} \in S^{a,-}_{\varepsilon}$.  After centering the linearized equations near $\bar{\zeta}$ by making the change of variables $s = (\bar{\zeta} - \zeta_n)/\eps$ and $\gamma_n = \zeta_n/\eps_n$, we arrive at the following recentered equations on the fast timescale:

\begin{align*}
\frac{d\hat{z}}{ds} &= \hat{a}(\hat{z},\gamma_n + s, \lambda,\eps_n)\\
\hat{z}(0) &= \hat{Y}_n.
\end{align*}
Proceeding as before, we define an associated family of frozen systems. Let $\hat{z}_*(s,\lambda,n)$ denote the solution of the corresponding frozen system
\begin{align*}
\frac{d\hat{z}}{ds} &= \hat{a}(\hat{z}_*,\gamma_n,\lambda,\eps_n)\\
\hat{z}_*(0,\lambda,n) &= \hat{Y}_n.
\end{align*}

Now fix $\delta_2 > 0$ such that $\delta < \delta_2 < \delta_1$. By assumption, for each $n$ we have $\hat{Y}_n \notin N_{\delta} (\hat{\sigma}_s(\gamma_n, \lambda,\eps_n))$, and so it must be true by Lemma \ref{lem:slowmanestimate} that
\begin{align*}
\hat{z}(\gamma_n-T,\lambda,n) &\in N_{\delta_2}(\hat{f}_f(\gamma_n-T,\lambda,\eps_n)).
\end{align*}

Since $\delta_2 < \delta_1$, $\hat{y}$ must enter the fast elephant trunk about $\hat{f}_f$ over $S^{a,-}_{\varepsilon}$, which contradicts the assumed behavior of $\hat{y}$ as $\zeta \to -\infty$.\\

We now know that $\hat{Y}(\zeta,\lambda,\eps) \in N_{\delta}(\hat{\Sigma}(\zeta,\lambda,\eps))$ for each $\zeta \leq -a$. Then set $\delta > 0$ small enough that for each sufficiently small $\eps > 0$, Lemma \ref{lem:exchange} applies to $\hat{Y}(\zeta,\lambda,\eps)$  as the wave enters the vicinity of $\bar{x}:=X(a,\eps) \in S^{a,+}_{\varepsilon}$  from the fast layer. In particular, we have the asymptotics  that $\hat{Y}$ lies $\mathcal{O}(\eps)$-close to $\hat{\Sigma}$ at the time $\zeta = a$, with respect to the Fubini-Study metric, and so for each sufficiently small $\eps > 0$, we have $\hat{Y} \in N_{\delta}(\hat{\Sigma}_s(a,\lambda,\eps))$. Simultaneously, $\hat{Y}$ enters the slow elephant trunk over $S^{a,+}_{\varepsilon}$ since $\delta < \delta_1$, and thus remains $\delta$-close to the slow subbundle for each $\zeta > a$. This completes the proof.  $\Box$

\begin{coro} \label{coro:slowsubbundleapprox}
Fix $\lambda \in \mathbb{C}$ and $a,\delta > 0$.  Then there exists $\bar{\eps} > 0$ sufficiently small so that for each $\eps \in (0,\bar{\eps}]$, the unique solution $\hat{E}(\zeta,\lambda)$ of the projectivized linearized equations for which $\hat{E}(\zeta,\lambda) \to \hat{e}_{s,2}^-$ (the unstable eigenvector) as $\zeta \to -\infty$, also satisfies
\begin{align}
\hat{E}(\zeta,\lambda,\eps) \in N_{\delta}(\hat{\Sigma}_s(\zeta,\lambda,\eps)).
\end{align}

for $|\zeta| \geq a$.
\end{coro}

\section{Proof of Lemma \ref{lem:exchange}} \label{app:exchange}

In this proof we focus on clarifying the essential steps in the estimate, omitting unwieldy calculations while noting that they can be traced from \cite{jonestin}. The key step is to write down a slightly modified version of the estimate given in Proposition 8 in \cite{jonestin}; namely, there exist positive constants $B_1,\,B_2,\,B_3$, and $T_0$ depending only on the width $\Delta$ of the defining box \eqref{eq:boxdef}, so that for all $T \geq T_0$ and for each $q \in Q_T$, we have the upper bounds
\begin{equation} \label{eq:modgrowthestimate}
\begin{aligned}
|dy(q(t))| &\leq B_1 e^{\alpha t}\\
|db(q(t))| &\leq B_2 e^{(\alpha-\kappa) t} + B_3 \eps \lambda e^{\alpha t}
\end{aligned}
\end{equation}

for all $t \in [0,T]$, where $\alpha>0$ and $0 < \kappa< |\gamma_0|$ are growth rate constants characterizing respectively slow growth versus fast exponential contraction (see Proposition 8 and Corollary 2 in \cite{jonestin}). We highlight the essential modification: whereas in the standard exchange lemma, the component $db(t)$ of the tangent vector that is aligned along the fast fibers can be arranged to shrink exponentially quickly over long times, i.e. with growth rate $(\alpha-\kappa) < 0$  (in the case of a slow manifold with only attracting directions), there is now an extra, slowly-growing error term, arising from the new $\eps \lambda$ term in the eigenvalue problem. \\

The strategy to prove \eqref{eq:modgrowthestimate} proceeds essentially as in the proof of Proposition 8 in \cite{jonestin}. By smoothness of solutions to the ODEs, there exists some $T'>0$ such that the upper bound for $|db(t)|$ holds for $t \in [0,T']$ if we choose $B_2$ and $B_3$ sufficiently large. We seek to show that we can take $T' = T$ for choices $B_2,\,B_3$ which are independent of the initial condition. Assume that there is a maximal $T' >0 $ such that the error estimate is attained  (otherwise there is nothing to prove). We show that this leads to a contradiction if $T$ is large enough.\\

By a similar calculation to that in the proof of Proposition 8 in \cite{jonestin}, we find that
\begin{align*}
|dy(t)| &\leq \bar{B}_1 e^{\alpha t},
\end{align*}

where the new constant $\bar{B}_1 > 0$ depends upon the width of the defining box $\Delta > 0$, as well as the growth rate constants. There is also an extra term of the form $e^{(\alpha - \kappa)t}$ arising from applying the Duhamel principle (see Lemma 5.1 in \cite{jonestin}) to the extra term of the form $\eps \lambda M_2 db$, where $M_2$ is a smooth bounded function defined by compositions of the Fenichel straightening diffeomorphism and its inverse. We subsume this extra error term into the constant $\bar{B}_1$. \\

By the Duhamel principle and the bounding estimates on points in $Q_T$ together with the basic flow estimates for slowly-varying nonautonomous linear systems (see Proposition 2, Corollary 1, Proposition 6, and Proposition 7 in \cite{jonestin}), we then have
\begin{align*}
|db(t)| &\leq  \bar{K}_X e^{(\alpha - \kappa) t} \bar{M}_2 + \eps \lambda \int_0^{t}e^{-\mu(t-\zeta)} |M_2(b(\zeta),y(\zeta),\eps)| |dy(\zeta)| d\zeta\\
&\leq \bar{K}_X e^{(\alpha-\kappa)t} \bar{M}_2 + \eps \lambda M_{\Delta,\bar{\eps}}e^{-\mu t}\bar{B}_1 \frac{e^{(\mu+\alpha)t}}{\mu+\alpha}\\
&< B_2 e^{(\alpha-\kappa)t} + B_3 e^{\alpha t}
\end{align*}
 
 if $T$ is chosen large enough and for $B_2,\,B_3$ defined from the coefficients as in the above calculation. This contradicts the assumed maximality of $T'$, and so we conclude that $db(t)$ satisfies the required inequality for all sufficiently large $t$. \\
 
 The remainder of the proof follows the outlines of the main Theorem 6.5 in \cite{jonestin}.  By assumption we have that $|dy(0)| \geq M\eps$ for some $M > 0$, and let $\beta < 0$ be fixed. It follows from general considerations about slowly varying nonautonomous linear systems (see Proposition 2 in \cite{jonestin}) together with the Duhamel principle that  
 \begin{align*}
|dy(\cT_\eps)| \geq K_y \eps e^{\beta \mathcal{T}_q}.
\end{align*}

The result follows from applying the estimate for $|dy(t)|$ above and calculating the distance between $(db(\cT_\eps),dy(\cT_\eps))$ and $(0,dy(\cT_\eps))$ in the Fubini-Study metric. $\Box$

\begin{remk}
An elephant trunk lemma is used in Gardner and Jones to prove an analogous result. We point out that only partial elephant trunk-type estimates over $S^{a,-}_{\varepsilon}$ and $S^{a,+}_{\varepsilon}$ are required to provide an estimate of the type proven above.
\end{remk}

\section{Example for modified exchange lemma estimates} \label{app:toy}

We illustrate the assertions in Remark \ref{remk:exchange} with the following toy problem defined in a suitable box $\Delta$ in $\mathbb{R}^2 \times \mathbb{C}^2$, where $\lambda \in \mathbb{C}$ is taken to be a fixed constant:
\begin{equation}
\label{eq:toyexchange2}
\begin{aligned}
b' &= -b\\
y' &= \eps y\\
db' &= -db + \eps \lambda y\,dy\\
dy' &= \eps\, dy.
\end{aligned}
\end{equation} 

The critical manifold $\mathcal{M}_0$ of this system is given by $\{b=0,\,db=0\}$. This manifold is normally hyperbolic and attracting since the real parts of the nontrivial eigenvalues of the corresponding layer problem of \eqref{eq:toyexchange2} along $\mathcal{M}_0$ are both negative in $\Delta$, and hence there exists a one-parameter family $\mathcal{M}_{\eps}$ of attracting slow manifolds for sufficiently small values of $\eps>0$. Without loss of generality we think of this family as parametrized by the slow variables $(y,\,dy)$, and evidently $b = 0$ specifies one of the defining equations for $\mathcal{M}_{\eps}$. The slow subbundle along $\mathcal{M}_{\eps}$ in terms of this coordinate representation can be computed by calculating the eigenvector corresponding to the $\mathcal{O}(\eps)$ eigenvalue of the $2\times 2$ Jacobian of the latter two equations in  \eqref{eq:toyexchange2}. It is given by the span of the vector
\begin{align*}
\begin{pmatrix} \frac{\eps\lambda y}{1+\eps}\\1 \end{pmatrix}
\end{align*}

as $y$ varies. On the other hand, the system \eqref{eq:toyexchange2} can be explicitly solved for $(db(t),dy(t))$ given initial values $(b(0),y(0),db(0),dy(0)) = (b_0,y_0,db_0,dy_0)$ to find
\begin{align*}
db(t) &= \left( \frac{1}{1+2\eps} \right) \left( e^{-t}[db_0(1+2\eps) - \eps \lambda y_0 dy_0] + e^{2\eps t}\eps \lambda y_0 dy_0   \right)\\
dy(t) &= dy_0 e^{\eps t}. 
\end{align*}
From this calculation, it follows that $\mathcal{M}_{\eps}$ is given by $\{b=0,\,db=\eps \lambda y\,dy/(1+2\eps)\}$.  \\

For times $t = \mathcal{T} = \mathcal{O}(1/\eps)$, We find that the tangent vector $(db(t),dy(t))$ lies exponentially close to 
\begin{align*}
\begin{pmatrix} \frac{\eps \lambda y_0 dy_0}{1+2\eps} e^{2\eps \mathcal{T}}  \\ dy_0 e^{\eps \mathcal{T}} \end{pmatrix}.
\end{align*}
 
Evidently, this vector is $\mathcal{O}(\eps)$ close to the tangent bundle of the slow manifold, which has the local representation $$\begin{pmatrix} 0 \\ 1 \end{pmatrix}.$$  Further comparing the angle of this vector to the slow subbundle at the point $y = y_0 e^{\eps \mathcal{T}}$ by using the distance estimate Lemma 3.1 in \cite{jonestin}, we find that the angle is no larger than 
 \begin{align*}
\frac{2\eps^2 \lambda y_0}{(1+2\eps)(1+\eps)},
\end{align*}
i.e. the angle scales as $\mathcal{O}(\eps^2)$. This error estimate is sharper than that given by the tangent bundle of the slow manifold.

\end{document}